\documentstyle[11pt]{amsart}

\def\R{{\rm I}\hspace*{-0.8mm}{\rm R}}
\def\no{\nonumber}
\def\be{\begin{equation}}
\def\ee{\end{equation}}
\def\ba{\begin{eqnarray}}
\def\ea{\end{eqnarray}}
\newcommand{\beq}{\begin{equation}}
\newcommand{\eeq}{\end{equation}}
\newcommand{\ben}{\begin{eqnarray}}
\newcommand{\een}{\end{eqnarray}}
\newcommand{\beno}{\begin{eqnarray*}}
\newcommand{\eeno}{\end{eqnarray*}}
\def\tilde{\widetilde}

\def\e1{\epsilon}
\def\AAl{\mathcal{A}_{\lambda}}
\def\A0{\stackrel{\circ}{\AAl}}

\def\o1{\omega}
\def\01{\Omega}
\def\c1{\gamma}

\def\g1{\Sigma}

\def\l1{\Lambda}

\def\v1{\varphi}

\def\d1{\delta}
\def\part{\partial}

\def\f1{\frac}
\def\t1{\theta}
\def\s1{\sqrt{\e1}}
\def\b1{\beta}
\def\bar{\overline}
\def\bs{\begin{eqnarray*}}
\def\es{\end{eqnarray*}}
\def\m1{\Theta}
\def\w1{\wedge}
\def\mv{{\mathbf{v}}}
\def\mn{{\mathbf{n}}}
\def\md{{\mathbf{D}}}
\def\mp{{\mathbf{P}}}
\def\momega{{\mathbf{\Omega}}}
\def\msigma{{\mathbf{\sigma}}}
\def\mkappa{{\mathbf{\kappa}}}
\def\mcurl{{\mathbf{ curl}}}
\def\mN{{\mathbf{N}}}
\def\mh{{\mathbf{h}}}
\def\mve{{\mathbf{v}}_{\epsilon}}
\def\mne{{\mathbf{n}}_{\epsilon}}
\def\mde{{\mathbf{D}}_{\epsilon}}
\def\je{{\cal J}_{\epsilon}}
\def\momegae{{\mathbf{\Omega}_{\epsilon}}}

\def\mhe{{\mathbf{h}}_{\epsilon}}
\newcommand{\BS}{{\mathbb{S}}}

\numberwithin{equation}{section}

\newtheorem{theorem}{Theorem}[section]
\newtheorem{corollary}[theorem]{Corollary}
\newtheorem{remark}[theorem]{Remark}
\newtheorem{lemma}[theorem]{Lemma}
\newtheorem{proposition}[theorem]{Proposition}

\begin{document}

\title{Global existence of weak solution for the 2-D Ericksen-Leslie system}

\author{Meng WANG}

\address{Department of Mathematics, Zhejiang University, Hangzhou
310027, P. R. China}
\email{mathdreamcn@@zju.edu.cn}

\author{Wendong Wang}

\address{School of Mathematical Sciences, Dalian University of Technology, Dalian, 116024, P. R. China
and The Institute of Mathematical Sciences, The Chinese University of Hong Kong, Shatin, N.T., Hong Kong}
\email{wendong@@dlut.edu.cn}

\date{\today}

%\keywords{Global well-posedness, Ericksen-Leslie system}
%\vspace{.2in}

\begin{abstract}
We prove the global existence of weak solution for two dimensional Ericksen-Leslie system with the Leslie stress and general Ericksen stress
under the physical constrains on the Leslie coefficients. We also prove the local well-posedness of the Ericksen-Leslie system
in two and three spatial dimensions.
\end{abstract}

\maketitle

\section{Introduction}

\subsection{Ericksen-Leslie system}

The hydrodynamic theory of liquid crystals was established by  Ericksen \cite{Er, Er2} and Leslie \cite{Les} in the 1960's.
In this theory, the configuration of the liquid crystals is described by a director field $\mn=(n^1,n^2,n^3)\in \BS^2$.
The general Ericksen-Leslie system in $\R^3$ takes the form
\ben\label{eq:EL}
(EL)\left\{
\begin{array}{l}
\mv_{t}+\mv\cdot\nabla \mv=-\nabla p+\frac{\gamma}{Re}\Delta \mv+\frac{1-\gamma}{Re}\nabla\cdot\msigma,\\
\nabla\cdot\mv=0,\\
\mn\times (\mh-\gamma_1\mN-\gamma_2\md\cdot\mn)=0,
\end{array}\right.
\een
where $\mv=(v^1,v^2,v^3)$ is the velocity of the fluid, $p$ is the pressure, $Re$ is the Reynolds number and $\gamma\in(0,1)$.
We denote
\beno
&\mkappa=(\nabla\mv)^{T},\quad \md=\frac{1}{2}(\mkappa^{T}+\mkappa),\quad \momega=\frac{1}{2}(\mkappa^{T}-\mkappa),\\
&\mN=\mn_{t}+\mv\cdot\nabla\mn+\momega\cdot \mn,
\eeno
and obviously $\mN$ is vertical with the director field $\mn$.
The stress $\sigma$ is modeled by the phenomenological constitutive relation
$$\sigma=\sigma^{L}+\sigma^{E},$$
where $\msigma^{L}$ is the viscous (Leslie) stress defined by
\ben
\msigma^{L}=\alpha_1(\mn\mn:\md)\mn\mn+\alpha_2\mn\mN+\alpha_3\mN\mn
+\alpha_4\md+\alpha_5\mn\mn\cdot\md+\alpha_6\md\cdot\mn\mn,
\een
where the six constants $\alpha_1,...\alpha_6$ are called the Leslie coefficients, $\mn\mn:\md=\sum_{i,j}n^i\md_{ij}n^j$ and $\mn\mN=(n^i\mN^j)_{3\times3}$; $\msigma^{E}$ is the elastic (Ericksen) stress
defined by
\ben
\msigma^{E}=-\frac{\partial W}{\partial(\nabla \mn)}\cdot(\nabla \mn)^{T},
\een
where $W=W(\mn,\nabla\mn)$ is the Oseen-Frank density depending on the elastic constants $k_1,k_2,k_3,k_4$ with the form
\beno
W=k_1({\rm div}\mn)^2+k_2|\mn\times({\mathbf{\nabla}}\times\mn)|^2+k_3|\mn\cdot({\mathbf{\nabla}}\times\mn)|^2 +k_4\big(\textrm{tr}(\nabla\mn)^2-({\rm div}\mn)^2\big).
\eeno
As in \cite{GMS3,HX}, we rewrite $W$ as
\ba
W(\mn,\nabla\mn)=a|\nabla\mn|^2+V(\mn,\nabla\mn),
\ea
where $a=\min\{k_1,k_2,k_3\}$ and
$$
V(\mn,\nabla\mn)=(k_1-a)({\rm div}\mn)^2+(k_2-a)|\mn\times({\mathbf{\nabla}}\times\mn)|^2
+(k_3-a)|\mn\cdot(\mathbf{\nabla}\times\mn)|^2.
$$
The molecular field $\mh$ is given by
\ben
\mh=-\frac{\delta W}{\delta \mn}=(\nabla_{i}W_{p_{i}^{l}}-W_{n^{l}}),
\een
where $p_{i}^{l}=\nabla_i n^l$ and we adopt the standard summation convention.
Throughout this paper, we use the notations:
\beno
W_{n^i}=\frac {\partial W(\mn,\mp)} {\partial {n^i}},\quad W_{p_i^j}=\frac {\partial W(\mn, \mp)} {\partial {p_i^j}}.
\eeno

In order to ensure that the system (EL) has a basic energy law, the Leslie coefficients and the two constants $\gamma_1$, $\gamma_2$ should satisfy
\ben
&\alpha_2+\alpha_3=\alpha_6-\alpha_5,\label{Parodi}\\
&\gamma_1=\alpha_3-\alpha_2,\gamma_2=\alpha_6-\alpha_5,\label{Parodi-2}
\een
where (\ref{Parodi}) is called Parodi's relation.
We denote
$$
\beta_1=\alpha_1+\frac{\gamma_2^2}{\gamma_1},\quad \beta_2=\alpha_4,\quad \beta_3=\alpha_5+\alpha_6-\frac{\gamma_2^2}{\gamma_1}.
$$
A necessary and sufficient condition which ensures that the energy of (EL) in $\R^2$ or $\R^3$ is dissipated is
\begin{equation}\label{beta}
\left\{\begin{array}{ll}
\beta_2\ge 0,\quad 2\beta_2+\beta_3\ge 0,\quad \frac{3}{2}\beta_2+\beta_3+\beta_1\ge 0,\quad {\rm for }\,\,\R^3\\
\beta_2\ge 0,\beta_1+2\beta_2+\beta_3\ge 0,\beta_1< 0,\quad {\rm or}\quad
\beta_2\ge 0,2\beta_2+\beta_3\ge 0,\beta_1\ge 0,\quad {\rm for }\,\,\R^2
\end{array}\right.
\end{equation}
which
was introduced by Wang-Zhang-Zhang \cite{WZZ} for $\R^3$. Moreover, for (EL) in $\R^2$, the the condition is weaker and proved in Section 2 (see Remark 2.2 for details).

Let
$\mu_1=\frac{1}{\gamma_1}$ and $\mu_2=-\frac{\gamma_2}{\gamma_1}.$
The third equation of (\ref{eq:EL}) is equivalent to
\ben\label{eq:equal n}
\mn_t+\mv\cdot\nabla\mn+\mn\times((\momega\cdot\mn-\mu_1\mh-\mu_2\md\cdot\mn)\times \mn)=0.
\een

\subsection{Main results}

Most of earlier works treated the approximated or simplified system of (\ref{eq:EL}), since the general Ericksen-Leslie system is very complicated. Lin and Liu \cite{LL-ARMA} consider the  Ginzburg-Landau type approximation of (\ref{eq:EL}):
\begin{eqnarray}\label{eq:EL-pen}
\left\{
\begin{array}{l}
\mv_t+\mv\cdot\nabla\mv=-\nabla{p}+\frac{\gamma}{Re}\Delta\mv
+\frac{1-\gamma}{Re}\nabla\cdot\sigma,\\
\mn_t+\mv\cdot\nabla\mn+\momega\cdot\mn-\mu_1\Delta\mn-\mu_2\md\cdot\mn-\frac 1 {\varepsilon^2}(|\mn|^2-1)\mn=0,
\end{array}\right.
\end{eqnarray}
which is obtained by adding the penality term $\frac 1 {\varepsilon^2}(|\mn|^2-1)\mn$ in $W$.
The global existence of weak solution and the local existence and uniqueness of strong solution
of the system (\ref{eq:EL-pen}) were proved in \cite{LL-ARMA} under certain strong constrains on the Leslie coefficients.  However, whether the solution of (\ref{eq:EL-pen}) converges to that
of (\ref{eq:EL}) as $\varepsilon$ tends to zero remains open.

A simplest system preserving the basic energy law is
\begin{eqnarray}\label{eq:EL-simple}
\left\{
\begin{array}{l}
\mv_t+\mv\cdot\nabla\mv-\Delta\mv+\nabla{p}=-\nabla\cdot\big(\nabla \mn\odot \nabla \mn\big),\\
\mn_t+\mv\cdot\nabla\mn-\Delta \mn=|\nabla\mn|^2\mn,
\end{array}\right.
\end{eqnarray}
which is obtained by neglecting the Leslie stress and taking the elastic constants in $W$ as $k_1=k_2=k_3=1$.
In two dimensional case, the global existence of weak solution has been independently proved by Lin, Lin and Wang \cite{LLW} and Hong \cite{Hong},
where they construct a class of weak solution with at most a finite number of singular times. The uniqueness of weak solution is proved by Lin-Wang \cite{LW}
and Xu-Zhang \cite{XZ}. Recently, Hong and Xin \cite{HX} extended the result of \cite{LLW,Hong} to the Oseen-Frank model with general Ericksen stress. In three dimensional case, the global existence of weak solution of (\ref{eq:EL-simple}) is a challenging open problem.
In the case when $|\nabla\mn|^2\mn$ in (\ref{eq:EL-simple}) is replaced by $\frac 1 {\varepsilon^2}(|\mn|^2-1)\mn$, the global existence and partial regularity
of weak solution  were studied in \cite{Lin1, Lin2}.
We refer to \cite{Wang, WXL, HW} and references therein for more relevant results.

In a recent work \cite{WZZ}, Wang-Zhang-Zhang proved the local well-posedness of the Ericksen-Leslie system, and the global well-posednss for small
initial data under the physical constrain conditions (\ref{Parodi})-(\ref{beta}) on the Leslie coefficients. In \cite{WZZ}, they considered the Ericksen
stress with $k_1=k_2=k_3=1$. In this paper, we first extend their result to the case with general Ericksen stress,
which will be used in the proof of global existence of weak solution. It's worth mentioning that recently Hong-Li-Xin \cite{HLX} obtained the local well-posed results of the liquid crystal flow for the Oseen-Frank model without Leslie sress in $\R^3$ by Ginzburg-Landau approximation approach.

The first result of the paper is the local existence, uniqueness and blow-up criterion for strong solutions of the Ericksen-Leslie system (\ref{eq:EL}) with the Leslie stress and general Ericksen stress, which generalized the results in \cite{WZZ} and \cite{HLX}.
\begin{theorem}\label{thm:local}
Assume that the Leslie coefficients satisfy (\ref{Parodi})-(\ref{beta}). Let $s\ge 2$  be an integer,
the initial data $\nabla\mn_0\in H^{2s}(\R^d)$, $\mv_0\in H^{2s}(\R^d)$ for $d=2$ (or $d=3$), and ${\rm div}\mv_0=0$.
There exist $T>0$ and a unique solution $(\mv,\mn)$ of the Ericksen-Leslie system (\ref{eq:EL}) such that
\beno
\mv\in C\big([0,T];H^{2s}(\R^d)\big)\cap L^{2}(0,T;H^{2s+1}(\R^d)),\quad
\nabla\mn\in C\big([0,T];H^{2s}(\R^d)\big).
\eeno
Let  $T^{*}$ be the maximal existence time of the solution. If $T^{*}<+\infty$, then it is necessary to hold that
$$
\int_{0}^{T^{*}}\|\nabla\times \mv(t)\|_{L^{\infty}}+\|\nabla\mn(t)\|_{L^{\infty}}^{2}dt =+\infty.
$$
\end{theorem}

Our second main goal is to extend Hong-Xin's global existence result of weak solution in 2-D in \cite{HX} to the case with the Leslie stress. In the space $\R^2$, $(\mv,\mn)$ in (\ref{eq:EL}) satisfies $v^3=0$, $\partial_{x_3}\mv=0$ and $\partial_{x_3}\mn=0$. Moreover, we assume that $\nabla\cdot\msigma^{L}$ means
$\left(
                                           \begin{array}{ccc}
                                            1 & 0&0\\
                                             0&1& 0\\
                                              0&0 & 0\\
                                           \end{array}
                                         \right)\cdot(\nabla\cdot\msigma^{L}).$
Let $b\in \BS^2$ be a constant vector and we define
$$
H_{b}^{1}(\R^2;\BS^2)=\big\{u: u-b\in H^{1}(\R^2;\R^3), |u|=1 \,\,{\rm a.e. }\,\,\mbox{in} \,\, \R^2 \big\}.
$$

\begin{theorem}\label{thm:global}
 Assume that the Leslie coefficients satisfy (\ref{Parodi})-(\ref{beta}), and the initial data $(\mv_0,\mn_0)\in L^{2}(\R^2)\times H_{b}^{1}(\R^2;\BS^2)$ with ${\rm div}\mv_0=0$.
Then there exists a global weak solution $(\mv,\mn)$ of the Ericksen-Leslie system (\ref{eq:EL}),
which is smooth in $\R^2\times((0,+\infty)\setminus \{T_{l}\}_{l=1}^{L})$ for a finite number of times $\{T_l\}_{l=1}^{L}$. Moreover, there are two constants $\epsilon_0>0$ and $R_0>0$ such that each singular point
$(x_{i}^{l},T_{l})$ is characterized by the condition
$$
\limsup_{t\uparrow T_{l}}\int_{B_{R}(x_{i}^{l})}\big(|\nabla\mn|^2+|\mv|^2\big)(\cdot,t)dx> \epsilon_0
$$
for any $R>0$  with $R\le R_0$.
\end{theorem}
Compared with Hong-Xin's results in \cite{HX}, we consider Ericksen-Leslie system (\ref{eq:EL}) with the Leslie stress, for which we need to balance the interaction between the Leslie stress and the Oseen-Frank density $W(\mn,\nabla\mn)$ (for example, see Prop 2.1, Theorem 3.2, Lemma 4.2 and 4.4). And we also explore an important decomposition formula for $\mh$ with two nonlinear terms. To obtain global existence result, firstly we prove the local existence for strong solutions of Ericksen-Leslie system (\ref{eq:EL}) by using  Friedrich's approach and uniform energy estimates as in \cite{WZZ}, and the difficulty comes from the nonlinear terms of the molecular field $\mh$. We overcome it by the estimates for commutator operator and detailed analysis of $\mh$. For the simplified case $k_1=k_2=k_3$ ($h=2a\Delta\mn$), it's easy to obtain a positive higher order dissipated energy, for example
\beno
\int_{\R^d}\mn\times(\mh\times\mn)\cdot\mh dx\geq 4a^2\int_{\R^d}|\Delta\mn|^2dx+{\rm Lower \,Order\, Terms}
\eeno
However, when $k_1, k_2, k_3$ are different, it's important for us to obtain the above estimates. In fact, for the proof of local existence in Section 3.1, we get a weaker higher order dissipated term $\|\Delta^s\je\mhe\times\mne\|_{L^2(\R^d)}$, which is enough to control error terms. Moreover, for blow-up criterion of strong solutions in Section 3.3, we find a useful observation:
\ben\label{nablawn0}
&&\nabla_{\alpha}W_{p_{\alpha}^{l}}\cdot n^l\\
&=&-2k_2|\nabla\mn|^2-2(k_3-k_2)(\mn\cdot\mcurl\mn)^2-2(k_1-k_2)({\rm div}\mn)^2+2(k_1-k_2)\nabla_l( n^l {\rm div}\mn),\nonumber
\een
which make that the term $\int_{\R^d}|\Delta^s(\nabla_{l}(n^l{\rm div}\mn) \mn)\cdot\Delta^{s+1} \mn|dx$ can be written as the sum of commutator terms, which can be controlled, seeing (\ref{eq:dot}). For the arguments for global existence of weak solutions, we obtain certain local monotonicity inequalities and interior regularity estimates, and we follow the basic spirit of Struwe \cite{St1}, which is later developed by Hong-Xin in \cite{HX}.

\begin{remark}
After we finished this paper,  Professor Changyou Wang told us that they also obtained similar results as Theorem 1.2 in a recent joint work with Jinrui Huang and Fanghua Lin. Also we thank Prof. Changyou Wang for giving us some valuable suggestions to improve this paper, especially higher regularity estimates in Section 4, which are different with the arguments in \cite{HX}.
\end{remark}

This paper is organized as follows: In section 2, we introduce the basic energy law of the Ericksen-Leslie system (\ref{eq:EL}) and the decomposition formula for $\mh$; In section 3, we prove the local existence, uniqueness and blow-up criterion for strong solutions of the Ericksen-Leslie system (\ref{eq:EL}) with the Leslie stress and general Ericksen stress by using Friedrich's approach and energy estimates, where the special structure of $\mh$ is frequently exploited and used. Section 4 is devoted to the proof of global existence of weak solutions.

%-----------------------------------------section 2-----------------------------------------------------------------
\section{Basic energy-dissipation law}

In this section, we derive the basic energy law of the Ericksen-Leslie system (\ref{eq:EL}) under the conditions (\ref{Parodi})-(\ref{beta})
on the Leslie coefficients. We consider the solution $(\mv,\mn)$ in $\R^d$ with $d=2,3$.

%--------------------------------------------------lemma1----------------------------------------------------------------
\begin{proposition}\label{energy}
 Let $(\mv,\mn)$ be a smooth solution of (\ref{eq:EL}) with the initial values $(\mv_0, \mn_0)$. Then it holds that
\ba\label{energye}
&&\int_{\R^d}e(\mv(\cdot,t),\mn(\cdot,t))dx+\int_{0}^{t}\int_{\R^d}\Big(\frac{\gamma}{1-\gamma}|\nabla \mv|^2+\frac{1}{\gamma_1}|\mn\times \mh|^2\Big)dxds\nonumber\\
&&\quad+\int_{0}^{t}\int_{\R^d}\Big((\alpha_1+\frac{\gamma_2^2}{\gamma_1})|\mn\mn:\md|^2+(\alpha_5
+\alpha_6-\frac{\gamma_{2}^{2}}{\gamma_1})|\md\cdot\mn|^2
+\alpha_4(\md:\md)\Big)dxds\nonumber\\
&&=\int_{\R^d}e(\mv_0,\mn_0)dx,
\ea
where $e(\mv,\mn)$ is defined by
$$e(\mv,\mn)=W(\mn,\nabla\mn)+\frac{Re}{2(1-\gamma)}|\mv|^2.$$
\end{proposition}

\begin{remark}
If the Leslie coefficients satisfy the first relation of $(\ref{beta})$, Wang-Zhang-Zhang \cite{WZZ} proved that for any symmetric trace free $3\times3$ matrix $\md$ and $\mn\in\BS^2$,
\beno
(\alpha_1+\frac{\gamma_2^2}{\gamma_1})|\mn\mn:\md|^2+(\alpha_5
+\alpha_6-\frac{\gamma_{2}^{2}}{\gamma_1})|\md\cdot\mn|^2
+\alpha_4(\md:\md)\ge 0.
\eeno
 Thus, the energy is dissipated in this case. Moreover, for any symmetric trace free $3\times3$ matrix $\md$ with $D_{i3}=0$ for $i=1,2,3$, the energy is dissipated if and only if
\beno\label{condition:2d}
\beta_2\ge 0,\beta_1+2\beta_2+\beta_3\ge 0,\beta_1< 0,\quad {\rm or}\quad
\beta_2\ge 0,2\beta_2+\beta_3\ge 0,\beta_1\ge 0,
\eeno
which is weaker than $(\ref{beta})_1$.
\end{remark}
\noindent{\bf Proof of Remark 2.2:} For any $\mn=(n_1,n_2,n_3)\in S^2$ and $\md=\left(
                                           \begin{array}{ccc}
                                             x & y & 0 \\
                                             y & -x & 0\\
                                             0 & 0 & 0\\
                                           \end{array}
                                         \right)$,
\begin{equation}\label{sn}
    \beta_1(\mn\mn:\md)^2+\beta_2(\md:\md)+\beta_3|\md\cdot\mn|^2\ge 0,
\end{equation}
which is equivalent to that
\begin{equation*}
    \beta_1\left((n_1^2-n_2^2)x+2n_1n_2y\right)^2+\left(2\beta_2+\beta_3(n_1^2+n_2^2)\right)(x^2+y^2)\ge 0.
\end{equation*}
Normalize $(x,y)$ such that $x=\cos\alpha,y=\sin\alpha$. Then
\begin{equation*}
    \beta_1\left((n_1^2-n_2^2)\cos\alpha+2n_1n_2\sin\alpha\right)^2+\left(2\beta_2+\beta_3(n_1^2+n_2^2)\right)\ge 0
\end{equation*}
for all $\alpha\in [0,2\pi]$ and $\mn=(n_1,n_2,n_3)\in S^2$.

For $\mn\in S^2$, we can set $n_1=t\cos\psi,n_2=t\sin\psi$, where $0\le t\le 1$, and $\psi\in[0,2\pi]$. Thus (\ref{sn}) holds if and only if that
\begin{equation*}
    \beta_1 \cos^{2}(2\psi-\alpha)+\beta_3 \frac{1}{s}+2\beta_2\frac{1}{s^2}\ge 0,
\end{equation*}
for all $s\in [0,1]$, and $\psi,\alpha\in [0,2\pi]$. In other words,
\begin{equation*}
\beta_2\ge 0,\beta_1+2\beta_2+\beta_3\ge 0,\beta_1< 0,\quad {\rm or}\quad
\beta_2\ge 0,2\beta_2+\beta_3\ge 0,\beta_1\ge 0.
\end{equation*}

\noindent{\bf Proof of Proposition 2.1.}\,Multiplying the first equation of (\ref{eq:EL}) by $\mv$ and using the fact $\nabla\cdot\mv=0$, we get
\ba\label{2.1}
&&\frac{1}{2}\frac d {dt}\int_{\R^d}|\mv|^2dx+\frac{\gamma}{Re}\int_{\R^d}|\nabla \mv|^2dx\nonumber\\
&&=-\frac{1-\gamma}{Re}\int_{\R^d}\msigma^{L}:\nabla\mv dx
+\frac{1-\gamma}{Re}\int_{\R^d}W_{p_{j}^{k}}(\mn,\nabla\mn)\nabla_{i}n^k \nabla_{j}v^{i}dx,
\ea
and it follows from the definitions of $\msigma^{L}, \md, \momega$, and (\ref{Parodi})-(\ref{Parodi-2}) that
\beno
&&\int_{\R^d}\msigma^{L}:\nabla\mv dx\nonumber\\
&&=\int_{\R^d}(\alpha_1(\mn\mn:\md)\mn\mn+\alpha_2\mn\mN+\alpha_3\mN\mn+\alpha_4\md+\alpha_5\mn\mn\cdot\md+\alpha_6\md\cdot\mn\mn):(\md+\momega)dx\nonumber\\
&&=\int_{\R^d}\big\{\alpha_1(\mn\mn:\md)^2+(\alpha_2+\alpha_3)\mn\cdot(\md\cdot\mN)+(\alpha_5+\alpha_6)|\md\cdot\mn|^2+\alpha_4\md:\md\nonumber\\
&&\quad+(\alpha_2-\alpha_3)\mn\cdot(\momega\cdot\mN)+(\alpha_6-\alpha_5)(\md\cdot\mn)\cdot(\momega\cdot\mn)\big\}dx\nonumber\\
&&=\int_{\R^d}\alpha_1(\mn\mn:\md)^2+(\alpha_5+\alpha_6)|\md\cdot\mn|^2+\alpha_4\md:\md dx\nonumber\\
&&\quad+\int_{\R^d}\gamma_1\mN\cdot(\momega\cdot\mn)+\gamma_2(\md\cdot\mn)\cdot(\momega\cdot\mn)+\gamma_2\mn\cdot(\md\cdot\mN)dx.\nonumber
\eeno
For the last term of the above equality, recall that the vector triple product formula $a\times(b\times c)=b(a\cdot c)-c(a\cdot b)$. Using the equation (\ref{eq:equal n}) and the  antisymmetry of $\momega$ ($\mn\cdot \momega\cdot\mn=0$), we get
\beno
&&\int_{\R^d}\gamma_1\mN\cdot(\momega\cdot\mn)+\gamma_2(\md\cdot\mn)\cdot(\momega\cdot\mn)+\gamma_2\mn\cdot(\md\cdot\mN)dx\nonumber\\
&&=\int_{\R^d}(\frac{1}{\gamma_1}\mn\times(\mh\times\mn)-\frac{\gamma_2}{\gamma_1}\mn\times(\md\cdot\mn\times\mn))
\cdot(\gamma_1\momega\cdot\mn+\gamma_2\md\cdot\mn)dx\\
&&\quad+\int_{\R^d}\gamma_2(\md\cdot\mn)\cdot(\momega\cdot\mn)dx\nonumber\\
&&=\int_{\R^d}\mh\cdot\momega\cdot\mn+\frac{\gamma_2}{\gamma_1}\int_{\R^d}\mn\times(\mh\times\mn)
\cdot(\md\cdot\mn)dx-\frac{\gamma_2^2}{\gamma_1}\int_{\R^d}(|\md\cdot\mn|^2-|\mn\cdot\md\cdot\mn|^2)dx.\nonumber
\eeno
Thus, we have
\ba\label{2.2}
-\int_{\R^d}\msigma^{L}:\nabla\mv dx&=&-\int_{\R^d}(\alpha_1+\frac{\gamma_2^2}{\gamma_1})(\mn\mn:\md)^2+
(\alpha_5+\alpha_6-\frac{\gamma_2^2}{\gamma_1})|\md\cdot\mn|^2+\alpha_4\md:\md dx\nonumber\\
&&-\int_{\R^d}\mh\cdot\momega\cdot\mn-\frac{\gamma_2}{\gamma_1}\int_{\R^d}\mn\times(\mh\times\mn)
\cdot(\md\cdot\mn)dx.
\ea
On the other hand, for the functional $W(\mn,\nabla\mn)$ we have
\ba\label{2.3}
\frac{d}{dt}\int_{\R^d}W(\mn,\nabla\mn)dx
&=&\int_{\R^{d}}W_{n^{l}}n^{l}_{t}+W_{p_{i}^{l}}\partial_t\nabla_i n^{l}dx\nonumber\\
&=&\int_{\R^d}\big(W_{n^{l}}-\nabla_{i}W_{p_{i}^{l}}\big)\big(n_{t}^{l}+\mv\cdot\nabla n^{l}-\mv\cdot\nabla n^{l}\big)dx.
\ea
Due to $\nabla\cdot \mv=0$, we get
\ba\label{2.4}
&&-\int_{\R^d}(W_{n^l}-\nabla_{i}W_{p_{i}^{l}})v^{k}\nabla_kn^{l}dx\nonumber\\
&&=-\int_{\R^d}W_{n^l}v^{k}\nabla_kn^{l}-\int_{\R^D}W_{p_{i}^{l}}v^k\nabla_{ik}^{2}n^ldx
-\int_{\R^d}W_{p_{i}^{l}}\nabla_{i}v^k\nabla_{k}n^ldx\nonumber\\
&&=-\int_{\R^d}v^k\cdot\nabla_{k}W-\int_{\R^D}W_{p_{i}^{l}}\nabla_{i}v^k\nabla_{k}n^ldx\nonumber\\
&&=-\int_{\R^d}W_{p_{i}^{l}}\nabla_{i}v^k\nabla_{k}n^ldx,
\ea
while by scalar triple product formula $a\cdot(b\times c)=b\cdot(c\times a)=c\cdot(a\times b)$, it's easy to derive that
\ba\label{2.5}
&&\int_{\R^d}\big(W_{n^{l}}-\nabla_{i}W_{p_{i}^{l}}\big)\cdot\big(n_{t}^{l}+\mv\cdot\nabla n^{l}\big)dx\nonumber\\
&&=\int_{\R^d}\big(-\mh\big)\cdot\big(-\momega\cdot\mn
+\frac{1}{\gamma_1}\mn\times(\mh\times\mn)-\frac{\gamma_2}{\gamma_1}\mn\times(\md\cdot\mn\times\mn)\big)dx\nonumber\\
&&=\int_{\R^d}\mh\cdot\momega\cdot\mn dx-\frac{1}{\gamma_1}\int_{\R^d}|\mh\times\mn|^2dx+
\frac{\gamma_2}{\gamma_1}\int_{\R^d}\mn\times(\mh\times\mn)
\cdot(\md\cdot\mn)dx.
\ea
Summing up for (\ref{2.1})-(\ref{2.5}), we deduce that
\beno\label{snow}
&&\frac{d}{dt}\int_{\R^d}\frac{Re}{2(1-\gamma)}|v|^2+W(\mn,\nabla\mn)dx+\frac{\gamma}{1-\gamma}\int_{\R^d}|\nabla v|^2dx+\frac{1}{\gamma_1}\int_{\R^d}|\mh\times\mn|^2dx\nonumber\\
&&=-\int_{\R^d}(\alpha_1+\frac{\gamma_2^2}{\gamma_1})(\mn\mn:\md)^2+
(\alpha_5+\alpha_6-\frac{\gamma_2^2}{\gamma_1})|\md\cdot\mn|^2+\alpha_4\md:\md dx\nonumber
\eeno
Then the Proposition follows by integrating on the time.\endproof

Now we will give a decomposition formula for $\mh$, which plays an important role in our proof.
\begin{lemma}\label{decomposition} For the terms $\nabla_{\alpha}W_{p_{\alpha}^{l}},$ $W_{n^{l}}$ and $\mh$, we have the following representation:
\ba\label{nabla w}
(\nabla_{\alpha}W_{p_{\alpha}^{l}})&=&2a\Delta\mn+2(k_1-a)\nabla{\rm div}\mn
-2(k_2-a){\mathbf curl}{(\mn\times(\mathbf curl\mn\times\mn))}\nonumber\\
&&-2(k_3-a){\mathbf curl}({\mathbf curl}\mn\cdot \mn\mn),\nonumber\\
&=&2a\Delta\mn+2(k_1-a)\nabla{\rm div}\mn -2(k_2-a){\mathbf curl}({\mcurl\mn})\nonumber\\
&&-2(k_3-k_2){\mathbf curl}({\mathbf curl}\mn\cdot \mn\mn)
\ea
\ben
(W_{n^{l}})=2(k_3-k_2)({\mcurl}\mn\cdot\mn)({\mcurl\mn}),
\een
\begin{eqnarray}\label{h}
\mh
&=&2a\Delta\mn+2(k_1-a)\nabla{\rm div}\mn -2(k_2-a){\mathbf curl}({\mcurl\mn})\nonumber\\
&&-2(k_3-k_2){\mathbf curl}({\mathbf curl}\mn\cdot \mn\mn) -2(k_3-k_2)({\mcurl}\mn\cdot\mn)({\mcurl\mn})
\end{eqnarray}
%Especially, when $k_2=k_3$,
% $$
% \mh=(\nabla_{\alpha}W_{p_{\alpha}^{l}})=2a\Delta\mn+2(k_1-a)\nabla{\rm div}\mn-2(k_2-a){\mathbf curl}{\mathbf curl}\mn.
% $$
 \end{lemma}
{\it Proof:}   The proof is direct. Note that $(b\times c)\cdot (b\times c)=|b|^2|c|^2-(b\cdot c)^2,$ then
$$
\begin{array}{ll}
(W_{p_{\alpha}^{i}})_{i}^{\alpha}
&= 2a\left(\begin{array}{lll}\partial_1 n^1&\partial_2 n^1&\partial_3 n^1\\
                            \partial_1 n^2& \partial_2 n^2&\partial_3 n^2\\
                            \partial_1 n^3& \partial_2 n^3&\partial_3 n^3        \end{array}\right)
+2(k_1-a)\left(\begin{array}{lll}{\rm div}\mn&0&0\\
                                   0&{\rm div}\mn&0\\
                                  0&0&{\rm div}\mn         \end{array}\right)\\
                                  &\quad+2(k_2-a)\left(\begin{array}{lll}0&-(\partial_1 n^2-\partial_2 n^1 )&(\partial_3n^1-\partial_1n^3)\\
                                                                (\partial_1 n^2-\partial_2 n^1 )& 0&-(\partial_2n^{3}-\partial_3 n^2)\\
                                                               -(\partial_3n^1-\partial_1n^3)& (\partial_2 n^3-\partial_3 n^2)&0         \end{array}\right)\\
                                   &\quad +2(k_3-k_2)\left(\begin{array}{lll}0&-n^3(\mn\cdot {\rm curl}\mn)&n^2(\mn\cdot {\rm curl}\mn)\\
                                                                            n^3(\mn\cdot {\rm curl}\mn)  &0& -n^1(\mn\cdot {\rm curl}\mn)\\
                                                                            -n^2(\mn\cdot {\rm curl}\mn)&n^1 (\mn\cdot {\rm curl}\mn)&0
                                                                            \end{array}\right)\\
                                                                            \end{array}
$$
Then  it's easy to obtain (\ref{nabla w}) by making $\nabla^{\alpha}$ on both sides of the above equality.

%--------------------------------------------local solution--------------------------------------------------
\section{Local well-posedness, uniqueness and blow-up criterion}
Throughout this paper, we denote that $C$ is a constant depending on $d$, $\alpha_1,\cdots,\alpha_6$, $k_1$, $k_2$, $k_3$, $\gamma$, $Re$ and independent of the solution $(\mv,\mn)$, and different from line to line. The symbol $\langle\cdot,\cdot\rangle$ denotes the integral in $\R^d$ with $d=2,3$. Moreover, ${\mathcal P}(\cdot,\cdots,\cdot)$ denotes the polynomial depending on the variable quantities in the bracket whose order, for example, is less than $10$.
\subsection{Existence}
Firstly we use the classical  Friedrich's method to construct the approximate solutions of (\ref{eq:EL}) as in \cite{WZZ}. One of the main difference is that the representation formula of $\mh$ owns three different positive coefficients $k_1, k_2$ and $k_3$ and $\mh$ is nonlinear with respect to $\mn$.

We will frequently use the following lemma for the commutator; for example see \cite{bcd}.
\begin{lemma}\label{lem:galiado}
For $\alpha,\beta\in N^3$, it holds that
$$
\|D^{\alpha}(fg)\|_{L^2}\le C\sum_{|\gamma|=|\alpha|}\big(\|f\|_{L^{\infty}}\|D^{\gamma}g\|_{L^2}+\|g\|_{L^{\infty}}\|D^{\gamma}f\|_{L^2}\big),
$$
$$
\|[D^{\alpha},f]D^{\beta}g\|_{L^2}\le C\left(\sum_{|\gamma|=|\alpha|+|\beta|}\|D^{\gamma}f\|_{L^2}\|g\|_{L^{\infty}}+\sum_{|\gamma|=|\alpha|+|\beta|-1}\|\nabla f\|_{L^{\infty}}\|D^{\gamma}g\|_{L^2}\right).
$$
\end{lemma}

The local existence of (\ref{eq:EL}) is split into two steps.

{\bf Step 1.} Construction of the approximated solutions.

In order to construct an approximated system preserving the energy-dissipation law,  Wang-Zhang-Zhang \cite{WZZ} introduced the following equivalent system of (\ref{eq:EL})
\begin{equation}\label{eq:EL-m}
\left\{\begin{array}{l}
\mv_t+\mv\cdot\nabla\mv=-\nabla p+\frac{\gamma}{Re}\Delta\mv+\frac{1-\gamma}{Re}\nabla\cdot(\tilde{\sigma}^{L}+\sigma^{E}),\\
{\rm div}\mv=0,\\
\mn_t+\mv\cdot\nabla\mn+\mn\times\left((\momega\cdot\mn-\mu_1\mh-\mu_2\md\cdot\mn)\times\mn\right)=0,
\end{array}\right.
\end{equation}
where $\tilde{\sigma}^{L}=\sigma_{1}(\mv,\mn)+\sigma_2(\mn)$ with
\beno
&&\sigma_{1}(\mv,\mn)=\beta_1(\mn\mn:\md)\mn\mn+\beta_2|\mn|^4\md+\frac{\beta_3}{2}|\mn|^2(\mn\md\cdot\mn+\md\cdot\mn\mn),\\
&&\sigma_2(\mn)=\frac{1}{2}(-1-\mu_2)\mn\left(\mn\times(\mh\times\mn)\right)
+\frac{1}{2}(1-\mu_2)\left(\mn\times(\mh\times\mn)\right)\mn.
\eeno
It's easy to check that the above system is just (\ref{eq:EL}) when $|\mn|=1$.

Let
$$
{\cal J}_{\epsilon}f={\cal F}^{-1}(\phi( \frac{\xi}{\epsilon}){\cal F}f),
$$
where $\cal F$ is usual Fourier transform and $\phi(\xi)$ is a smooth cut-off function with $\phi=1$ in $B_1$ and $\phi=0$ outside of $B_2$. Let $\mathbf P$ be an operator which projects a vector field to its solenoidal part.
We construct the approximate system of (\ref{eq:EL-m}):
\begin{equation}\label{approximate}
\left\{
\begin{array}{llll}
\displaystyle\frac{\partial\mv_{\epsilon}}{\partial t}+{\cal J}_{\epsilon}{\mathbf P}({\cal J}_{\epsilon}\mv_{\epsilon}\cdot\nabla {\cal J}_{\epsilon}\mv_{\epsilon})\\
\quad\quad=\frac{\gamma}{Re}\je\Delta{\cal J}_{\epsilon}\mv_{\epsilon}+\frac{1-\gamma}{Re}\nabla\cdot{\cal J}_{\epsilon}{\mathbf P}\left(\sigma_1({\cal J}_{\epsilon}\mv_{\epsilon},{\cal J}_{\epsilon}{\mn_{\epsilon}})+\sigma_2({\cal J}_{\epsilon}\mn_{\epsilon})+\sigma^{E}({\cal J}_{\epsilon}\mn_{\epsilon})\right),\\
{\rm div}\mv^{\epsilon}=0,\\
\displaystyle\frac{\partial \mn_{\epsilon}}{\partial t}+{\cal J}_{\epsilon}\left({\cal J}_{\epsilon}\mv_{\epsilon}\cdot\nabla{\cal J}_{\epsilon}\mn_{\epsilon}+{\cal J}_{\epsilon}\mn_{\epsilon}\times\big(({\cal  J}_{\epsilon}\momega_{\epsilon}\cdot{\cal J}_{\epsilon}\mn_{\epsilon}-\mu_1{\cal J}_{\epsilon}\mh_{\epsilon}-\mu_2{\cal J }_{\epsilon}\md_{\epsilon}\cdot{\cal J}_{\epsilon}\mn_{\epsilon})\times {\cal J}_{\epsilon}\mn_{\epsilon}\big)\right)=0,\\
(\mv_{\epsilon},\mn_{\epsilon})|_{t=0}=({\cal J}_{\epsilon}\mv_{0},{\cal J}_{\epsilon}\mn_0).
\end{array}
\right.
\end{equation}
where
\begin{eqnarray*}{\cal J}_{\epsilon}h_{\epsilon}&=&2a\Delta{\cal J}_{\epsilon}\mne+2(k_1-a)\nabla{\rm div}{\cal J}_{\epsilon}\mne-2(k_2-a)\mcurl(\je\mne\times(\mcurl\je\mne\times\je\mne))\\
&&-2(k_3-a)\mcurl({\cal J}_{\epsilon}\mne\cdot\mcurl{\cal J}_{\epsilon}\mne{\cal J}_{\epsilon}\mne)-2(k_3-k_2)({\cal J}_{\epsilon}\mne\cdot\mcurl{\cal J}_{\epsilon}\mne)(\mcurl{\cal J}_{\epsilon}\mne)\end{eqnarray*}

By Cauchy-Lipschitz theorem, we know that there exist a strictly maximal time $T_{\epsilon}$ and
a unique solution $(\mve, \mn_{\epsilon})\in C([0,T_{\epsilon});H^{k}(\R^d))$ for any $k\ge 0$.
It's worth to mention that the choosing of $\je$ is different from that in \cite{WZZ}. Since $\mh$ is nonlinear, we need to use the uniform integration of $\je$ to overcome the difficulty from the commutator terms, for example, the Lie bracket $[\je, f]$ in (\ref{commutetor}), which needs much regularity of the cut-off function $\phi$.

{\bf Step 2.} Uniform energy estimates.

We introduce the following energy functional, which is related with the formula of $\mh$,
\begin{eqnarray*}
E_{s}(\mve,\mne)&=&\|\mne-\mn_0\|_{L^2}^2+\frac{Re}{2(1-\gamma)}\|\mve\|_{L^2}^{2}+\int_{\R^d}W(\mne,\nabla\mne)dx\\
&&+a\|\Delta^s\nabla\mne\|_{L^2}^2
+(k_1-a)\|\Delta^s{\rm div}\mne\|_{L^2}^2+(k_2-a)\|\je\mne\times\Delta^s\mcurl\mne\|_{L^2}^{2}\\
&&+(k_3-a)\|\je\mne\cdot\Delta^s\mcurl\mne\|_{L^2}^{2}
+\frac{Re}{2(1-\gamma)}\|\Delta^s\mve\|_{L^2}^{2}.
\end{eqnarray*}

{\bf Step 2.1.} Estimates of lower order terms in $E_{s}(\mve,\mne)$.
Similar arguments hold for $\tilde{\sigma}^{L}$ as in Proposition \ref{energy}, and other terms are directly estimated. In fact, the approximate system has the following energy estimate:
\ba\label{energye2}
&&\frac{d}{dt}\int_{\R^d}\frac{Re}{2(1-\gamma)}|\mv_{\epsilon}|^2+W(\mn_{\epsilon},\nabla\mn_{\epsilon})dx
+\frac{\gamma}{1-\gamma}\int_{\R^d}|\nabla\je\mv_{\epsilon}|^2dx\nonumber\\
&&=-\int_{\R^d} \beta_1|\je\md_{\epsilon}:\je\mn_{\epsilon}\je\mn_{\epsilon}|^2
+\beta_2|\je\mne|^4\je\md_{\epsilon}:\je\md_{\epsilon}+\beta_3|\je\mn_{\epsilon}|^2|\je\md_{\epsilon}\cdot\je\mn_{\epsilon}|^2dx
\nonumber\\
&&\quad+C{\mathcal P}\big(\|\mn_{\epsilon}\|_{L^{\infty}},\|\mv_{\epsilon}\|_{L^{\infty}},\|\nabla\mn_{\epsilon}\|_{L^{\infty}}\big)
\big(\|\nabla\mn_{\epsilon}\|_{H^{1}}^2+\|\nabla\mv_{\epsilon}\|_{L^{2}}^2\big)\nonumber\\
&&\le C{\mathcal P}\big(\|\mn_{\epsilon}\|_{L^{\infty}},\|\mv_{\epsilon}\|_{L^{\infty}},\|\nabla\mn_{\epsilon}\|_{L^{\infty}}\big)
\big(\|\nabla\mn_{\epsilon}\|_{H^{1}}^2+\|\nabla\mv_{\epsilon}\|_{L^{2}}^2\big)
\ea
Using (\ref{approximate}) and  the definition of $\je\mhe$, we have
\ba\label{n}
&&\frac{d}{dt}\|\mn_\epsilon-\mn_0\|_{L^2}^{2}=2\langle\partial_t\mn_{\epsilon},\mn_{\epsilon}-\mn_0\rangle\nonumber\\
&&\le C\left(\|\mv_{\epsilon}\|_{L^2}+\|\nabla\mv_{\epsilon}\|_{L^2}+\|\nabla\mn_{\epsilon}\|_{L^2}+\|\Delta\mn_{\epsilon}\|_{L^2}\right)
(1+\|\mn_{\epsilon}\|_{L^\infty}+\|\nabla\mn_{\epsilon}\|_{L^\infty})^{4}\|\mn_{\epsilon}-\mn_0\|_{L^2}\nonumber\\
&&\le C(1
+\|\mn_{\epsilon}\|_{L^\infty}+\|\nabla\mne\|_{L^\infty})^4E_{s}(\mv_{\epsilon},\mn_{\epsilon}).
\ea
{\bf Step 2.2.} We turn to the estimate of higher order derivatives of $\mne$ in $E_{s}(\mve,\mne)$.
\ba\label{highnn}
&&\frac12\frac{d}{dt}\langle\nabla\Delta^s\mne,\nabla\Delta^s\mne\rangle\nonumber\\
&&=-\langle\nabla\Delta^{s}(\je\mve\cdot\nabla\je\mne),\Delta^{s}\nabla\je\mne\rangle
+\langle\Delta^{s}\big(\je\mne\times\left((\je\momegae\cdot\je\mne)\times\je\mne\right)\big),\Delta^{s+1}\je\mne\rangle\nonumber\\
&&\quad-\mu_2\langle\Delta^{s}\big(\je\mne\times\left((\je\mde\cdot\je\mne)\times\je\mne\right)\big),\Delta^{s+1}\je\mne\rangle\nonumber\\
&&\quad-\mu_1\langle\Delta^{s}\big(\je\mne\times(\je\mhe\times\je\mne)\big),\Delta^{s+1}\je\mne\rangle\nonumber\\
&&\doteq I_1+I_2+I_3+I_4.
\ea
In the same way, we have
\ba\label{highdn}
\frac12\frac{d}{dt}\langle{\rm div}\Delta^s\mne,\rm div\Delta^s\mne\rangle
&=&-\langle\Delta^{s}{\rm div}(\je\mve\cdot\nabla\je\mne),\Delta^{s}{\rm div}\je\mne\rangle\nonumber\\
&&+\langle\Delta^{s}\big(\je\mne\times\left((\je\momegae\cdot\je\mne)\times\je\mne\right)\big),\Delta^{s}\nabla{\rm div}\je\mne\rangle\nonumber\\
&&-\mu_2\langle\Delta^{s}\big(\je\mne\times\left((\je\mde\cdot\je\mne)\times\je\mne\right)\big),\Delta^{s}\nabla{\rm div}\je\mne\rangle\nonumber\\
&&-\mu_1\langle\Delta^{s}\big(\je\mne\times(\je\mhe\times\je\mne)\big),\Delta^{s}\nabla{\rm div}\je\mne\rangle\nonumber\\
&\doteq& I_1'+I_2'+I_3'+I_4'.
\ea
%and recall a simple relation $\langle\nabla\times f,g\rangle=\langle f,\nabla\times g\rangle$, then
%\ba\label{highcn}
%\frac12\frac{d}{dt}\langle\Delta^s{\mcurl}\mne,\Delta^s\mcurl\mne\rangle
%&=&-\langle\mcurl \Delta^{s}(\je\mve\cdot\nabla\je\mne),\Delta^{s}\mcurl\je\mne\rangle\nonumber\\
%&&-\langle\Delta^{s}\big(\je\mne\times\left((\je\momegae\cdot\je\mne)\times\je\mne\right)\big),\Delta^{s}{\mcurl\mcurl}\je\mne\rangle\nonumber\\
%&& +\mu_2\langle\Delta^{s}\big(\je\mne\times\left((\je\mde\cdot\je\mne)\times\je\mne\right)\big),\Delta^{s}\mcurl\mcurl\je\mne\rangle\nonumber\\
%&&+\mu_1\langle\Delta^{s}\big(\je\mne\times(\je\mhe\times\je\mne)\big),\Delta^{s}\mcurl\mcurl\je\mne\rangle\nonumber\\
%&\doteq& I_1''+I_2''+I_3''+I_4''.
%\ea
For the nonlinear term $\frac{d}{dt}\langle\je\mne\times\Delta^s\mcurl\mne,\je\mne\times\Delta^s\mcurl\mne\rangle$, since ${\rm div}\mv^{\epsilon}=0$, (\ref{approximate}) and Lemma 3.1 yield that
\beno
&&\frac12\frac{d}{dt}\langle\je\mne\times\Delta^s\mcurl\mne,\je\mne\times\Delta^s\mcurl\mne\rangle\nonumber\\
&&=\langle\je\partial_t\mne\times\Delta^s\mcurl\mne-\je\mne\times\Delta^{s}\mcurl({\cal J}_{\epsilon}(\je\mve\cdot\nabla\je\mne)),
\je\mne\times\Delta^s\mcurl\mne\rangle\nonumber\\
&&\quad+\langle\je\mne\times\Delta^{s}\mcurl(\partial_t\mne+{\cal J}_{\epsilon}(\je\mve\cdot\nabla\je\mne)),
\je\mne\times\Delta^s\mcurl\mne\rangle\\
&&\le -\langle\je\mne\times\Delta^{s}\mcurl({\cal J}_{\epsilon}(\je\mve\cdot\nabla\je\mne)),
\je\mne\times\Delta^s\mcurl\mne\rangle\nonumber\\
&&\quad+\langle\je\mne\times\Delta^{s}\mcurl(\partial_t\mne+{\cal J}_{\epsilon}(\je\mve\cdot\nabla\je\mne)),
\je\mne\times\Delta^s\mcurl\mne\rangle\\
&&\quad+C(\|\nabla\mn_{\epsilon}\|_{L^{\infty}}\|\mv_{\epsilon}\|_{L^{\infty}}+\|\nabla\mv_{\epsilon}\|_{L^{\infty}}+
\|\nabla^2\mn_{\epsilon}\|_{L^{\infty}}+\|\nabla\mn_{\epsilon}\|_{L^{\infty}}^2)(1+\|\mn_{\epsilon}\|_{L^{\infty}})^5 \|\nabla\mn_{\epsilon}\|_{H^{2s}}^2\\
&&\le\langle[\je,\je\mne\times\je\mne\times]\Delta^s\mcurl\mne,\Delta^{s}\mcurl(\je\mve\cdot\nabla\je\mne)\rangle\nonumber\\
&&\quad+\langle\je\mne\times\je\mne\times\Delta^s\mcurl\je\mne,[\Delta^s\mcurl,\je\mve\cdot]\nabla\je\mne\rangle\nonumber\\
&&\quad+\langle\je\mne\times\Delta^{s}\mcurl(\partial_t\mne+{\cal J}_{\epsilon}(\je\mve\cdot\nabla\je\mne)),
\je\mne\times\Delta^s\mcurl\mne\rangle\\
&&\quad+C{\mathcal P}(\|\mn_{\epsilon}\|_{L^{\infty}},\|\nabla\mn_{\epsilon}\|_{L^{\infty}},\|\mv_{\epsilon}\|_{L^{\infty}},\|\nabla\mv_{\epsilon}\|_{L^{\infty}},
\|\nabla^2\mn_{\epsilon}\|_{L^{\infty}}) \|\nabla\mn_{\epsilon}\|_{H^{2s}}^2\\
&&\le\langle[\je,\je\mne\times\je\mne\times]\Delta^s\mcurl\mne,\Delta^{s}\mcurl(\je\mve\cdot\nabla\je\mne)\rangle\nonumber\\
&&\quad-\langle\je\mne\times\Delta^{s}\mcurl{\cal J}_{\epsilon}(\je\mne\times((\je\momegae\cdot\je\mne)\times\je\mne)),
\je\mne\times\Delta^s\mcurl\mne\rangle\\
&&\quad+\mu_2\langle\je\mne\times\Delta^{s}\mcurl{\cal J}_{\epsilon}(\je\mne\times((\je\mde\cdot\je\mne)\times\je\mne)),
\je\mne\times\Delta^s\mcurl\mne\rangle\\
&&\quad+\mu_1\langle\je\mne\times\Delta^{s}\mcurl{\cal J}_{\epsilon}(\je\mne\times(\je\mhe\times\je\mne)),
\je\mne\times\Delta^s\mcurl\mne\rangle\\
&&\quad+C(\delta){\mathcal P}(\|\mn_{\epsilon}\|_{L^{\infty}},\|\nabla\mn_{\epsilon}\|_{L^{\infty}},\|\mv_{\epsilon}\|_{L^{\infty}},\|\nabla\mv_{\epsilon}\|_{L^{\infty}},
\|\nabla^2\mn_{\epsilon}\|_{L^{\infty}}) \|\nabla\mn_{\epsilon}\|_{H^{2s}}^2+\delta\|\nabla\je\mv_{\epsilon}\|_{H^{2s}}^2\\
&&\doteq I_1''+I_2''+I_3''+I_4''+\delta\|\nabla\je\mv_{\epsilon}\|_{H^{2s}}^2\\
&&\quad+C(\delta){\mathcal P}(\|\mn_{\epsilon}\|_{L^{\infty}},\|\nabla\mn_{\epsilon}\|_{L^{\infty}},\|\mv_{\epsilon}\|_{L^{\infty}},\|\nabla\mv_{\epsilon}\|_{L^{\infty}},
\|\nabla^2\mn_{\epsilon}\|_{L^{\infty}}) \|\nabla\mn_{\epsilon}\|_{H^{2s}}^2
\eeno
where $\delta>0$, to be decided later. Similarly, we can obtain
\beno
&&\frac12\frac{d}{dt}\langle\je\mne\cdot\Delta^s\mcurl\mne,\je\mne\cdot\Delta^s\mcurl\mne\rangle\nonumber\\
&&\le-\langle[\je,\je\mne\je\mne\cdot]\Delta^s\mcurl\mne,\Delta^{s}\mcurl(\je\mve\cdot\nabla\je\mne)\rangle\nonumber\\
&&\quad-\langle\je\mne\cdot\Delta^{s}\mcurl{\cal J}_{\epsilon}(\je\mne\times((\je\momegae\cdot\je\mne)\times\je\mne)),
\je\mne\cdot\Delta^s\mcurl\mne\rangle\\
&&\quad+\mu_2\langle\je\mne\cdot\Delta^{s}\mcurl{\cal J}_{\epsilon}(\je\mne\times((\je\mde\cdot\je\mne)\times\je\mne)),
\je\mne\cdot\Delta^s\mcurl\mne\rangle\\
&&\quad+\mu_1\langle\je\mne\cdot\Delta^{s}\mcurl{\cal J}_{\epsilon}(\je\mne\times(\je\mhe\times\je\mne)),
\je\mne\cdot\Delta^s\mcurl\mne\rangle\\
&&\quad+C(\delta){\mathcal P}(\|\mn_{\epsilon}\|_{L^{\infty}},\|\nabla\mn_{\epsilon}\|_{L^{\infty}},\|\mv_{\epsilon}\|_{L^{\infty}},\|\nabla\mv_{\epsilon}\|_{L^{\infty}},
\|\nabla^2\mn_{\epsilon}\|_{L^{\infty}}) \|\nabla\mn_{\epsilon}\|_{H^{2s}}^2+\delta\|\nabla\je\mv_{\epsilon}\|_{H^{2s}}^2\\
&&\doteq I_1'''+I_2'''+I_3'''+I_4'''+\delta\|\nabla\je\mv_{\epsilon}\|_{H^{2s}}^2\\
&&\quad+C(\delta){\mathcal P}(\|\mn_{\epsilon}\|_{L^{\infty}},\|\nabla\mn_{\epsilon}\|_{L^{\infty}},\|\mv_{\epsilon}\|_{L^{\infty}},\|\nabla\mv_{\epsilon}\|_{L^{\infty}},
\|\nabla^2\mn_{\epsilon}\|_{L^{\infty}}) \|\nabla\mn_{\epsilon}\|_{H^{2s}}^2
\eeno

Firstly, we estimate the terms of $I_1$, $I_1'$, $I_1''$ and $I_1'''$.
Due to ${\rm div}\mv^{\epsilon}=0$ and Lemma 3.1, we have
\beno
&&|I_1|+|I_1'|\nonumber\\
&&\le |\langle[\nabla\Delta^s,\je\mv_{\epsilon}\cdot]\nabla\je\mn_{\epsilon},\Delta^{s}\nabla\je\mn_{\epsilon}\rangle|+
|\langle[{\rm div}\Delta^s,\je\mv_{\epsilon}\cdot]\nabla\je\mn_{\epsilon},\Delta^{s}{\rm div}\je\mn_{\epsilon}\rangle|\nonumber\\
&&\le C_{\delta}(\|\nabla\mve\|_{L^{\infty}}+\|\nabla\mne\|_{L^{\infty}}^{2})\|\nabla\mne\|_{H^{2s}}^{2}
+\delta\|\nabla\je\mve\|_{H^{2s}}^{2}
\eeno
As to $I_1''$ and $I_1'''$, we have the following estimate
\beno
[{\cal J}_{\epsilon},f]\nabla_jg
&=&\int_{\R^d}\phi_{\epsilon}(y)f(x-y)\nabla_jg(x-y)dy-\int_{R^{d}}\phi_{\epsilon}(y)f(x)\nabla_jg(x-y)dy\nonumber\\
&\leq&|\int_{\R^d}\nabla_j\phi_{\epsilon}(y)\int_{0}^{1}y\cdot\nabla f(x-\tau y)d\tau g(x-y)dy|\nonumber\\
&&+|\int_{\R^d}\phi_{\epsilon}(y)\nabla_j f(x-y) g(x-y)dy|,
\eeno
where $\phi_{\epsilon}(x)=\frac{1}{\epsilon^d}\phi(\frac{x}{\epsilon})$. Hence, for $1\leq p\leq\infty$, by Young inequality we get
\begin{eqnarray}\label{commutetor}
\|[{\cal J}_{\epsilon},f]\nabla_jg\|_{L^p}
&\leq&C\big(\|\nabla_j\phi_{\epsilon}(y)y\|_{L^1}\|\nabla f\|_{L^{\infty}}+\|\phi_{\epsilon}(y)\|_{L^1}\|\nabla f\|_{L^{\infty}}\big)\|g\|_{L^p}\nonumber\\
&\leq& C(1+\|\nabla f\|_{L^{\infty}})\|g\|_{L^p}
\end{eqnarray}
Applying (\ref{commutetor}) and Lemma \ref{lem:galiado} to $I_1''$, we obtain
\beno
|I_1''|&\leq&|\langle[\je,\je\mne\times\je\mne\times]\Delta^s\mcurl\nabla_j\mne,\Delta^{s-1}\nabla_j\mcurl(\je\mve\cdot\nabla\je\mne)\rangle|
+\delta\|\nabla\je\mv_{\epsilon}\|_{H^{2s}}^2\\
&&+C(\delta){\mathcal P}(\|\mn_{\epsilon}\|_{L^{\infty}},\|\nabla\mn_{\epsilon}\|_{L^{\infty}},\|\mv_{\epsilon}\|_{L^{\infty}},\|\nabla\mv_{\epsilon}\|_{L^{\infty}},
\|\nabla^2\mn_{\epsilon}\|_{L^{\infty}}) \|\nabla\mn_{\epsilon}\|_{H^{2s}}^2\\
&\le&2\delta\|\nabla\je\mv_{\epsilon}\|_{H^{2s}}^2+C(\delta){\mathcal P}(\|\mn_{\epsilon}\|_{L^{\infty}},\|\nabla\mn_{\epsilon}\|_{L^{\infty}},\|\mv_{\epsilon}\|_{L^{\infty}},\|\nabla\mv_{\epsilon}\|_{L^{\infty}},
\|\nabla^2\mn_{\epsilon}\|_{L^{\infty}}) \|\nabla\mn_{\epsilon}\|_{H^{2s}}^2
\eeno
Similar estimates hold for the term $I_1'''$,
thus we have
\ben\label{i1}
&&|I_1|+|I_1'|+|I_1''|+|I_1'''|\nonumber\\
&&\le 5\delta\|\nabla\je\mv_{\epsilon}\|_{H^{2s}}^2+C(\delta){\mathcal P}(\|\mn_{\epsilon}\|_{L^{\infty}},\|\nabla\mn_{\epsilon}\|_{L^{\infty}},\|\mv_{\epsilon}\|_{L^{\infty}},\|\nabla\mv_{\epsilon}\|_{L^{\infty}},
\|\nabla^2\mn_{\epsilon}\|_{L^{\infty}}) \|\nabla\mn_{\epsilon}\|_{H^{2s}}^2\nonumber\\
\een
Secondly, for the terms $I_2-I_2'''$, by the formula of $\je\mh$, Lemma \ref{lem:galiado} and (\ref{commutetor}) we have
\ben\label{i2}
&&2a I_2+2(k_1-a)I_2'+2(k_2-a)I_2''+2(k_3-a)I_2'''\nonumber\\
&&\le \langle \Delta^{s}\big(\je\mne\times\left((\je\momegae\cdot\je\mne)\times\je\mne\right)\big),\Delta^{s}\je\mhe   \rangle\nonumber\\
&&\quad+2(k_3-k_2)\langle \Delta^{s}\big(\je\mne\times\left((\je\momegae\cdot\je\mne)\times\je\mne\right)\big), \Delta^{s} ((\mcurl\je\mne\cdot\je\mne)\mcurl\je\mne)\rangle\nonumber\\
&&\quad+C|\langle[\Delta^s,\je\mne\times\je\mne\times]\mcurl\mcurl\je\mne,\Delta^{s}(\je\mne\times((\je\momegae\cdot\je\mne)\times\je\mne))\rangle|\nonumber\\
&&\quad+C|\langle[\Delta^s,\je\mne\je\mne\cdot]\mcurl\mcurl\je\mne,\Delta^{s}(\je\mne\times((\je\momegae\cdot\je\mne)\times\je\mne))\rangle|\nonumber\\
&&\quad+C|\langle[\je,\je\mne\times\je\mne\times]\Delta^s\mcurl\mcurl\mne,\Delta^{s}(\je\mne\times((\je\momegae\cdot\je\mne)\times\je\mne))\rangle|\nonumber\\
&&\quad+C|\langle[\je,\je\mne\je\mne\cdot]\Delta^s\mcurl\mcurl\mne,\Delta^{s}(\je\mne\times((\je\momegae\cdot\je\mne)\times\je\mne))\rangle|\nonumber\\
&&\quad+C(\delta){\mathcal P}(\|\mn_{\epsilon}\|_{L^{\infty}},\|\nabla\mn_{\epsilon}\|_{L^{\infty}},\|\mv_{\epsilon}\|_{L^{\infty}},\|\nabla\mv_{\epsilon}\|_{L^{\infty}},
\|\nabla^2\mn_{\epsilon}\|_{L^{\infty}}) \|\nabla\mn_{\epsilon}\|_{H^{2s}}^2
+\delta\|\nabla\je\mve\|_{H^{2s}}^{2}\nonumber\\
&&\le \langle \big(\je\mne\times\left((\je\Delta^{s}\momegae\cdot\je\mne)\times\je\mne\right)\big),\Delta^{s}\je\mhe   \rangle\nonumber\\
&&\quad+|\langle [\Delta^s,\je\mne\times(\je\mne\times(\je\mne\cdot))]\nabla_l\momegae,\Delta^{s-1}\nabla_l\je\mhe   \rangle|\nonumber\\
&&\quad+C(\delta){\mathcal P}(\|\mn_{\epsilon}\|_{L^{\infty}},\|\nabla\mn_{\epsilon}\|_{L^{\infty}},\|\mv_{\epsilon}\|_{L^{\infty}},\|\nabla\mv_{\epsilon}\|_{L^{\infty}},
\|\nabla^2\mn_{\epsilon}\|_{L^{\infty}}) \|\nabla\mn_{\epsilon}\|_{H^{2s}}^2
+2\delta\|\nabla\je\mve\|_{H^{2s}}^{2},\nonumber\\
&&\le \langle \big(\je\mne\times\left((\je\Delta^{s}\momegae\cdot\je\mne)\times\je\mne\right)\big),\Delta^{s}\je\mhe   \rangle\nonumber\\
&&\quad+C(\delta){\mathcal P}(\|\mn_{\epsilon}\|_{L^{\infty}},\|\nabla\mn_{\epsilon}\|_{L^{\infty}},\|\mv_{\epsilon}\|_{L^{\infty}},\|\nabla\mv_{\epsilon}\|_{L^{\infty}},
\|\nabla^2\mn_{\epsilon}\|_{L^{\infty}}) \|\nabla\mn_{\epsilon}\|_{H^{2s}}^2
+3\delta\|\nabla\je\mve\|_{H^{2s}}^{2}.\nonumber\\
\een
For the terms $I_3,\cdots,I_3'''$, similar arguments yield that
\ben\label{i3}
&&2a I_3+2(k_1-a)I_3'+2(k_2-a)I_3''+2(k_3-a)I_3'''\nonumber\\
&&\le -\mu_2\langle \big(\je\mne\times\left((\je\Delta^{s}\mde\cdot\je\mne)\times\je\mne\right)\big),\Delta^{s}\je\mhe   \rangle\nonumber\\
&&\quad+C(\delta){\mathcal P}(\|\mn_{\epsilon}\|_{L^{\infty}},\|\nabla\mn_{\epsilon}\|_{L^{\infty}},\|\mv_{\epsilon}\|_{L^{\infty}},\|\nabla\mv_{\epsilon}\|_{L^{\infty}},
\|\nabla^2\mn_{\epsilon}\|_{L^{\infty}}) \|\nabla\mn_{\epsilon}\|_{H^{2s}}^2
+\delta\|\nabla\je\mve\|_{H^{2s}}^{2}.\nonumber\\
\een
Moreover, since using Lemma \ref{lem:galiado} and integration by parts, we have
\ben\label{the term h}
&&\|\Delta^s\big(\je\mne\times(\je\mhe\times\je\mne)\big)\|_{L^2(\R^d)}\nonumber\\
&& \leq C {\mathcal P}(\|\mn_{\epsilon}\|_{L^{\infty}},\|\nabla\mn_{\epsilon}\|_{L^{\infty}},
\|\nabla^2\mn_{\epsilon}\|_{L^{\infty}}) (\|\nabla\mn_{\epsilon}\|_{H^{2s}}+\|\nabla^{2s}(\je\mhe\times\je\mne)\|_{L^2})\nonumber\\
&&\leq C {\mathcal P}(\|\mn_{\epsilon}\|_{L^{\infty}},\|\nabla\mn_{\epsilon}\|_{L^{\infty}},
\|\nabla^2\mn_{\epsilon}\|_{L^{\infty}}) \nonumber\\
&&\quad \cdot(\|\nabla\mn_{\epsilon}\|_{H^{2s}}+\|\Delta^{s}\je\mhe\times\je\mne\|_{L^2}+\|[\Delta^{s},\je\mne\times]\je\mhe\|_{L^2})\nonumber\\
&&\leq C {\mathcal P}(\|\mn_{\epsilon}\|_{L^{\infty}},\|\nabla\mn_{\epsilon}\|_{L^{\infty}},
\|\nabla^2\mn_{\epsilon}\|_{L^{\infty}})\big(\|\nabla\mn_{\epsilon}\|_{H^{2s}}+\|\Delta^{s}\je\mhe\times\je\mne\|_{L^2}\big).
\een
Hence, for the terms $I_4,\cdots, I_4'''$, again using Lemma \ref{lem:galiado}, by the formula of $\je\mhe$ we get
\ben\label{i4}
&&2a I_4+2(k_1-a)I_4'+2(k_2-a)I_4''+2(k_3-a)I_4'''\nonumber\\
&&\le-\mu_1\langle\Delta^s\big(\je\mne\times(\je\mhe\times\je\mne)\big),\Delta^{s}\je\mhe\rangle\nonumber\\
&&\quad+C(\delta){\mathcal P}(\|\mn_{\epsilon}\|_{L^{\infty}},\|\nabla\mn_{\epsilon}\|_{L^{\infty}},
\|\nabla^2\mn_{\epsilon}\|_{L^{\infty}}) \|\nabla\mn_{\epsilon}\|_{H^{2s}}^2+\delta\|\Delta^s\je\mhe\times\je\mne\|_{L^2}^{2}\nonumber\\
&&\le-\mu_1\langle\Delta^s\je\mhe\times\je\mne,\Delta^s\je\mhe\times\je\mne\rangle\nonumber\\
&&\quad-\mu_1\langle[\Delta^{s},\je\mne\times](\je\mhe\times\je\mne),\Delta^{s}\je\mhe\rangle+\mu_1\langle[\Delta^{s},\je\mne\times]\je\mhe,\Delta^{s}\je\mhe\times\je\mne\rangle\nonumber\\
&&\quad+C(\delta){\mathcal P}(\|\mn_{\epsilon}\|_{L^{\infty}},\|\nabla\mn_{\epsilon}\|_{L^{\infty}},
\|\nabla^2\mn_{\epsilon}\|_{L^{\infty}}) \|\nabla\mn_{\epsilon}\|_{H^{2s}}^2+2\delta\|\Delta^s\je\mhe\times\je\mne\|_{L^2}^{2}\nonumber\\
&&\doteq-\mu_1\langle\Delta^s\je\mhe\times\je\mne,\Delta^s\je\mhe\times\je\mne\rangle+I_{41}+I_{42}\nonumber\\
&&\quad+C(\delta){\mathcal P}(\|\mn_{\epsilon}\|_{L^{\infty}},\|\nabla\mn_{\epsilon}\|_{L^{\infty}},
\|\nabla^2\mn_{\epsilon}\|_{L^{\infty}}) \|\nabla\mn_{\epsilon}\|_{H^{2s}}^2+2\delta\|\Delta^s\je\mhe\times\je\mne\|_{L^2}^{2}.
\een
Obviously, $I_{42}$ is bounded by
\ben\label{42}
|I_{42}|&\le& C(\delta){\mathcal P}(\|\mn_{\epsilon}\|_{L^{\infty}},\|\nabla\mn_{\epsilon}\|_{L^{\infty}},
\|\nabla^2\mn_{\epsilon}\|_{L^{\infty}}) \|\nabla\mn_{\epsilon}\|_{H^{2s}}^2+\delta\|\Delta^s\je\mhe\times\je\mne\|_{L^2}^{2}.\nonumber\\
\een
The term $I_{41}$ is a little complicated. In fact, using the formula of $\je\mhe$ we can rewrite $I_{41}$ as
\beno
I_{41}
&\le&2a\mu_1\langle\nabla[\Delta^{s},\je\mne\times](\je\mhe\times\je\mne),\nabla\Delta^s\je\mne\rangle\\
&&+2(k_1-a)\mu_1\langle{\rm div}[\Delta^{s},\je\mne\times](\je\mhe\times\je\mne),{\rm div}\Delta^s\je\mne\rangle\nonumber\\
&&+2(k_2-a)\mu_1\langle{\mcurl}[\Delta^{s},\je\mne\times](\je\mhe\times\je\mne),\Delta^s(\je\mne\times({\mcurl}\je\mne\times\je\mne))\rangle\nonumber\\
&&+2(k_3-a)\mu_1\langle{\mcurl}[\Delta^{s},\je\mne\times](\je\mhe\times\je\mne),\Delta^s(\mcurl\je\mne\cdot\je\mne\je\mne)\rangle\nonumber\\
&&+2(k_3-k_2)\mu_1\langle[\Delta^{s},\je\mne\times](\je\mhe\times\je\mne),\Delta^s\big((\mcurl\je\mne\cdot\je\mne)({\mcurl}\je\mne)\big)\rangle\nonumber\\
\eeno
For the first term of $I_{41}$, by Lemma \ref{lem:galiado} and integration by parts as in (\ref{the term h}), it's easy to derive that
\beno
&&\langle\nabla[\Delta^{s},\je\mne\times](\je\mhe\times\je\mne),\nabla\Delta^s\je\mne\rangle\\
&&= \langle[\Delta^s,\nabla\je\mne\times](\je\mhe\times\je\mne),\nabla\Delta^s\je\mne\rangle+
\langle[\Delta^s,\je\mne\times]\nabla(\je\mhe\times\je\mne),\nabla\Delta^s\je\mne\rangle\\
&&\leq  C{\mathcal P}(\|\mn_{\epsilon}\|_{L^{\infty}},\|\nabla\mn_{\epsilon}\|_{L^{\infty}},
\|\nabla^2\mn_{\epsilon}\|_{L^{\infty}}) \big(\|\nabla\mn_{\epsilon}\|_{H^{2s}}^2+ \|\nabla^{2s}(\je\mhe\times\je\mne)\|_{L^2}\|\nabla\mn_{\epsilon}\|_{H^{2s}}\big)\\
&&\leq  C(\delta){\mathcal P}(\|\mn_{\epsilon}\|_{L^{\infty}},\|\nabla\mn_{\epsilon}\|_{L^{\infty}},
\|\nabla^2\mn_{\epsilon}\|_{L^{\infty}}) \|\nabla\mn_{\epsilon}\|_{H^{2s}}^2+\delta\|(\Delta^s\je\mhe\times\je\mne)\|_{L^2}^2
\eeno
Other terms are similar to deal with, and finally it's easy to obtain
\ben\label{41}
|I_{41}|&\le& C(\delta){\mathcal P}(\|\mn_{\epsilon}\|_{L^{\infty}},\|\nabla\mn_{\epsilon}\|_{L^{\infty}},
\|\nabla^2\mn_{\epsilon}\|_{L^{\infty}}) \|\nabla\mn_{\epsilon}\|_{H^{2s}}^2+5\delta\|\Delta^s\je\mhe\times\je\mne\|_{L^2}^{2}.\nonumber\\
\een
It is concluded by (\ref{i4}), (\ref{42}) and (\ref{41}) that
\ben\label{4end}
&&2a I_4+2(k_1-a)I_4'+2(k_2-a)I_4''+2(k_3-a)I_4'''\nonumber\\
&&\le -\mu_1\langle\Delta^s\mhe\times\mne,\Delta^s\mhe\times\mne\rangle+C(\delta){\mathcal P}(\|\mn_{\epsilon}\|_{L^{\infty}},\|\nabla\mn_{\epsilon}\|_{L^{\infty}},
\|\nabla^2\mn_{\epsilon}\|_{L^{\infty}}) \|\nabla\mn_{\epsilon}\|_{H^{2s}}^2\nonumber\\
&&\quad+8\delta\|\Delta^s\je\mhe\times\je\mne\|_{L^2}^{2}.\nonumber\\
\een
Thus, by (\ref{energye2}), (\ref{n}), (\ref{i1}), (\ref{i2}), (\ref{i3}) and (\ref{4end}), we get
\ba\label{nnn}
&&\frac{d}{dt}\int_{\R^d}E_s(\mne,\mve)-\frac{Re}{2(1-\gamma)}\langle\Delta^s\mve,\Delta^s\mve\rangle dx+\mu_1\langle\Delta^s\je\mhe\times\je\mne,\Delta^s\je\mhe\times\je\mne\rangle\nonumber\\
&&\le \langle \je\mne\times\left((\Delta^{s}\je\momegae\cdot\je\mne)\times\je\mne\right),\Delta^{s}\je\mhe   \rangle\nonumber\\
&&\quad-\mu_2\langle \je\mne\times\left((\Delta^{s}\je\mde\cdot\je\mne)\times\je\mne\right),\Delta^{s}\je\mhe   \rangle\nonumber\\
&&\quad+ C({\delta}){\mathcal P}(\|\mne\|_{L^{\infty}},\|\nabla\mne\|_{L^{\infty}},\|\nabla^2\mne\|_{L^{\infty}},\|\mve\|_{L^{\infty}},\|\nabla\mve\|_{L^{\infty}})\|\nabla\mne\|_{H^{2s}}^2\nonumber\\
&&\quad+30\delta(\|\nabla\je\mve\|_{H^{2s}}^2
+\|\Delta^s\je\mhe\times\je\mne\|_{L^2}^2).
\ea

{\bf Step 2.3.} We consider the estimate of higher order derivatives for $\mve$ in $E_{s}(\mve,\mne)$. Acting the inner product
$\frac{Re}{(1-\gamma)}\Delta^s\mve$ on both sides of $(\ref{approximate})_1$, by the antisymmetry of $\je\momegae$ we get
\beno
&&\frac{Re}{2(1-\gamma)}\frac{d}{dt}\langle\Delta^s\mve,\Delta^s\mve\rangle
+\frac{\gamma}{1-\gamma}\langle\nabla\Delta^s\je\mve,\nabla\Delta^s\je\mve\rangle\nonumber\\
&&=-\frac{Re}{1-\gamma}\langle\Delta^s(\je\mve\cdot\nabla\je\mve),\Delta^s\je\mve\rangle
+\langle\Delta^s(W_{p_{j}^k}(\je\mne,\je\nabla\mne)\nabla_{i}\je\mne^k),\Delta^s\nabla_j\je\mve^i\rangle\nonumber\\
&&\quad-\langle\Delta^s\{\beta_1(\je\mne\je\mne:\je\mde)\je\mne\je\mne
+\frac{\beta_3}{2}|\je\mne|^2(\je\mne\je\mde\cdot\je\mne+\je\mde\cdot\je\mne\je\mne)\nonumber\\
&&\quad+\beta_2|\je\mne|^4\je\mde\},\Delta^s\je\mde
\rangle
+\mu_2\langle\Delta^s\left(\je\mne\times(\je\mhe\times\je\mne)\je\mne\right),\Delta^s\je\mde\rangle\nonumber\\
&&\quad-\langle\Delta^s\left(\je\mne\times(\je\mhe\times\je\mne)\je\mne\right),\Delta^s\je\momegae\rangle\nonumber\\
&&\doteq III_1+III_2+III_3+III_4+III_5.
\eeno
Then by Lemma 3.1 and ${\rm div}\je\mve=0$ we have
\beno
III_1\le C\|\nabla\mve\|_{L^{\infty}}\|\mve\|_{H^{2s}}^{2},
\eeno
\beno
III_2\le C({\delta}){\mathcal P}(\|\mne\|_{L^{\infty}},\|\nabla\mne\|_{L^{\infty}})\|\nabla\mne\|_{H^{2s}}^{2}+\delta\|\nabla\je\mve\|_{H^{2s}}^{2},
\eeno
\beno
III_4&\le& \mu_2\langle\je\mne\times(\Delta^s\je\mhe\times\je\mne),\Delta^s\je\mde\cdot\je\mne\rangle\nonumber\\
&&+C({\delta}){\mathcal P}(\|\mne\|_{L^{\infty}},\|\nabla\mne\|_{L^{\infty}},\|\nabla^2\mne\|_{L^{\infty}})\|\nabla\mne\|_{H^{2s}}^2
+\delta\|\nabla\je\mve\|_{H^{2s}}^{2},
\eeno
\beno
III_5&\le& -\langle\je\mne\times(\Delta^s\je\mhe\times\je\mne),\Delta^s\je\momegae\cdot\je\mne\rangle\nonumber\\
&&+C({\delta}){\mathcal P}(\|\mne\|_{L^{\infty}},\|\nabla\mne\|_{L^{\infty}},\|\nabla^2\mne\|_{L^{\infty}})\|\nabla\mne\|_{H^{2s}}^2
+\delta\|\nabla\je\mve\|_{H^{2s}}^{2},
\eeno
and by Remark 2.2,
\beno
III_3&\le&-\langle\beta_1(\je\mne\je\mne:\Delta^s\je\mde)\je\mne\je\mne+\beta_2|\je\mne|^4\Delta^s\je\mde,\Delta^s\je\mde\rangle
\\
&&-\frac{\beta_3}{2}|\je\mne|^2\langle(\je\mne\Delta^s\je\mde\cdot\je\mne+\Delta^s\je\mde\cdot\je\mne\je\mne),\Delta^s\je\mde\rangle\\
&&+C({\delta}){\mathcal P}(\|\mne\|_{L^{\infty}},\|\nabla\mne\|_{L^{\infty}},\|\nabla\mve\|_{L^{\infty}})(\|\je\mve\|_{H^{2s}}^{2}+\|\nabla\mne\|_{H^{2s}}^{2})
+\delta\|\nabla\je\mve\|_{H^{2s}}^2\\
&\leq& C({\delta}){\mathcal P}(\|\mne\|_{L^{\infty}},\|\nabla\mne\|_{L^{\infty}},\|\nabla\mve\|_{L^{\infty}})(\|\mve\|_{H^{2s}}^{2}+\|\nabla\mne\|_{H^{2s}}^{2})
+\delta\|\nabla\je\mve\|_{H^{2s}}^2.
\eeno
Summing up for $III_1-III_5$, we have
\ben\label{highv}
&&\frac{Re}{2(1-\gamma)}\frac{d}{dt}\langle\Delta^s\mve,\Delta^s\mve\rangle
+\frac{\gamma}{1-\gamma}\langle\Delta^s\nabla\je\mve,\Delta^s\nabla\je\mve\rangle\nonumber\\
&&\le \mu_2\langle\je\mne\times(\Delta^s\je\mhe\times\je\mne),\Delta^s\je\mde\cdot\je\mne\rangle\nonumber\\
&&\quad-\langle\je\mne\times(\Delta^s\je\mhe\times\je\mne),\Delta^s\je\momegae\cdot\je\mne\rangle
+5\delta\|\nabla\je\mve\|_{H^{2s}}^2\nonumber\\
&&\quad+C({\delta}){\mathcal P}(\|\mne\|_{L^{\infty}},\|\nabla\mne\|_{L^{\infty}},\|\nabla^2\mne\|_{L^{\infty}},
\|\mve\|_{L^{\infty}},\|\nabla\mve\|_{L^{\infty}})(\|\nabla\mne\|_{H^{2s}}^{2}+\|\mve\|_{H^{2s}}^{2})\nonumber\\
\een
Hence, choose $\delta$ small enough and by (\ref{nnn}), (\ref{highv}) we may show that
$$
\frac{d}{dt}E_{s}(\mve,\mne)+\frac{\gamma}{2(1-\gamma)}\|\Delta^s\nabla\je\mve\|_{L^2}^2
\le {\cal F}(E_{s}(\mve,\mne)),
$$
where $\cal F$ is an increasing function with ${\cal F}(0)=0$. It means that there exists a $T>0$ depending only on $E_{s}(\mv_0,\mn_0)$ such that for $\forall t\in [0,\min(T,T_{\epsilon})]$,
$$
E_{s}(\mve,\mne)(t)+\frac{\gamma}{2(1-\gamma)}\int_0^t\|\Delta^s\nabla\je\mve\|_{L^2}^2d\tau
\le 2E_{s}(\mv_0,\mn_0),
$$
which implies that $T_{\epsilon}\ge T$ by a continuous argument. Then the uniform estimates for the approximate solutions $(\mve,\mne)$ on $[0,T]$ hold, which yield that there exists a local solution $(\mv,\mn)$ of (\ref{eq:EL}) by the standard compactness arguments. Moreover, if $|\mn_0|=1$, multiply $\cdot\mn$ on both sides of (\ref{eq:equal n}) and we can obtain that $|\mn|=1$. Hence the proof is complete.

\subsection{Uniqueness}

The section is devoted to  the proof of uniqueness for strong solutions of (\ref{eq:EL}).

\begin{theorem}\label{uniqueness}
Assume that the Leslie coefficients satisfy (\ref{Parodi})-(\ref{beta}), and the initial data $\nabla\mn_0\in H^{2s}(\R^d)$, $\mv_0\in H^{2s}(\R^d)$ with ${\rm div}\mv_0=0$, where $s\geq 2$ and $d=2,3$. Then there exists a unique strong solution $(\mv,\mn)$ of the Ericksen-Leslie system (\ref{eq:EL}) in $\R^d\times (0,T)$ with the indicated data.
\end{theorem}

{\it Proof:}
For $(\mv_0,\mn_0)$ with $\nabla\mn_0\in H^{2s}(\R^d)$ and $\mv_0\in H^{2s}(\R^d)$ ($s\geq 2$), we may assume that there exist two strong solutions  $(\mv_1,\mn_1)$ and  $(\mv_2,\mn_2)$ in $\R^d\times (0,T)$ with the initial data, which satisfy
\ben\label{eq:bound}
\sup_{(x,t)\in (\R^d\times (0,T))}\sum_{i=1,2}\big(|\nabla^3 \mn_i|+|\nabla^2 \mn_i|+|\nabla \mn_i|+|\nabla^2 \mv_i|+|\nabla \mv_i|+|\mv_i|\big)\leq C.
\een

Firstly, we estimate $\|\mn_1-\mn_2\|_{L^2}$.
By the equation (\ref{eq:equal n}), we have
\begin{eqnarray*}
&&\frac{1}{2}\frac{d}{dt}\langle\mn_1-\mn_2,\mn_1-\mn_2\rangle\\
&&=\langle\mv_2\cdot\nabla\mn_2-\mv_1\cdot\nabla\mn_1,\mn_1-\mn_2\rangle+\langle\momega_2\cdot\mn_2-\momega_1\cdot\mn_1,\mn_1-\mn_2\rangle\\
&&\quad +\mu_1\langle\mn_1\times(\mh_1\times\mn_1)-\mn_2\times(\mh_2\times\mn_2),\mn_1-\mn_2\rangle\\
&&\quad +\mu_2\langle\mn_1\times(\md_1\cdot\mn_1\times\mn_1)-\mn_2\times(\md_2\cdot\mn_2\times\mn_2),\mn_1-\mn_2\rangle\\
&&\doteq A_1+A_2+A_3+A_4,
\end{eqnarray*}
where $\momega_i, \mh_i, \md_i$ represent the functions of $\mv_i,\mn_i$ for $i=1,2$, respectively. For simplicity's sake, we denote
$$
\delta\mn=\mn_1-\mn_2, \delta\mv=\mv_1-\mv_2, \delta\mh=\mh_1-\mh_2,$$
$$\delta\md=\md_1-\md_2,\delta\momega=\momega_1-\momega_2,\delta\mN=\mN_1-\mN_2.
$$
For the terms $A_1$ and $A_2$, by the assumption (\ref{eq:bound}) and integration by parts we have
\ben\label{eq:A12}
|A_1|+|A_2|\leq C\left(\|\nabla\delta\mn\|_{L^2}^2+\|\delta\mv\|_{L^2}^2+\|\delta\mn\|_{L^2}^2\right).
\een
Similarly, for the term $A_4$, it's easier to obtain
\ben\label{eq:A4}
A_{4}\leq C\left(\|\nabla\delta\mn\|_{L^2}^2+\|\delta\mv\|_{L^2}^2+\|\delta\mn\|_{L^2}^2\right).
\een
For the term $A_3$, by the above assumptions and the formula of $\mh$ (\ref{h}) we have the following estimates
\begin{eqnarray*}
A_3&\leq&\mu_1\langle\mn_1\times(\delta_{\mh}\times\mn_1),\delta\mn\rangle+C\|\delta\mn\|_{L^2}^{2}\\
&\leq &\mu_1\langle\delta\mh,\delta\mn\rangle-\mu_1\langle(\delta{\mn}\cdot\mn_1)\mn_1,\delta\mn\rangle+C\|\delta\mn\|_{L^2}^{2}\\
&\doteq&A_{31}+A_{32}+C\|\delta\mn\|_{L^2}^2,
\end{eqnarray*}
where
\begin{eqnarray*}
\delta\mh&=&2a(\Delta\delta\mn)+2(k_1-a)(\nabla{\rm div}\delta\mn)\\
&& -2(k_2-a)({\mathbf curl}(\mn_1\times({\mcurl\mn_1\times\mn_1}))-{\mathbf curl}(\mn_2\times({\mcurl\mn_2\times\mn_2})))\\
&&-2(k_3-a)({\mathbf curl}({\mathbf curl}\mn_1\cdot \mn_1\mn_1)-{\mathbf curl}({\mathbf curl}\mn_2\cdot \mn_2\mn_2))\\ &&-2(k_3-k_2)\big(({\mcurl}\mn_1\cdot\mn_1)({\mcurl\mn_1})-({\mcurl}\mn_2\cdot\mn_2)({\mcurl\mn_2})\big).
\end{eqnarray*}
Hence for $A_{31}$,
\begin{eqnarray}
A_{31}&\leq&-2a\mu_1\|\nabla\delta\mn\|_{L^2}^{2}-2\mu_1(k_1-a)\|{\rm div}\delta\mn\|_{L^2}^{2}-2\mu_1(k_2-a)\|\mcurl\delta{\mn}\times\mn_1\|_{L^2}^{2}\nonumber\\
&&-2\mu_1(k_3-a)\|\mcurl\delta\mn\cdot\mn_1\|_{L^2}^{2}+C(a,\mu_1)\|\delta\mn\|_{L^2}^{2}+a\mu_1\|\nabla\delta\mn\|_{L^2}^{2}\nonumber\\
&\leq&C\|\delta\mn\|_{L^2}^{2}.\label{eq:I31}
\end{eqnarray}
For $A_{32}$, similar arguments yield that
\begin{eqnarray}\label{eq:A32}
A_{32}\leq C\left(\|\delta\mn\|_{L^2}^2+\|\nabla\delta\mn\|_{L^2}^2\right).
\end{eqnarray}
Then, combining the above estimates we have
\begin{eqnarray}\label{eq:n1-n2}
\frac12\partial_t\|\delta\mn\|_{L^2}^2
\leq C\left(\|\nabla\delta\mn\|_{L^2}^{2}+\|\delta\mn\|_{L^2}^{2}+\|\delta\mv\|_{L^2}^{2}\right).
\end{eqnarray}
Secondly, we need to estimate $\|\delta\mv\|_{L^2}$. Multiplying $(\mv_1-\mv_2)$ on both sides of the equation $(\ref{eq:EL})_1$ and integrating we have
\begin{eqnarray*}
&&\frac12\partial_t\langle\delta\mv,\delta\mv\rangle+\frac{\gamma}{Re}\langle\nabla\delta\mv,\nabla\delta\mv\rangle\\
&&=\langle\mv_2\cdot\nabla\mv_2-\mv_1\cdot\nabla\mv_1,\delta\mv\rangle+\frac{1-\gamma}{Re}\langle\sigma^L_2-\sigma^L_1,\nabla\delta\mv\rangle
+\frac{1-\gamma}{Re}\langle\sigma^E_2-\sigma^E_1,\nabla\delta\mv\rangle\\
&&\doteq A_5+\frac{1-\gamma}{Re}A_6+A_7,
\end{eqnarray*}
Obviously,
\begin{eqnarray*}
A_{5}\leq C\|\delta\mv\|_{L^2}^{2}+\frac{\gamma}{16Re}\|\nabla\delta\mv\|_{L^2}^{2},
\end{eqnarray*}
and
\begin{eqnarray*}
A_{7}\leq C\big(\|\delta\mn\|_{L^2}^{2}+\|\nabla\delta\mn\|_{L^2}^{2}\big)+\frac{\gamma}{16Re}\|\nabla\delta\mv\|_{L^2}^2.
\end{eqnarray*}
The estimate of $A_6$ is difficult. First we note that, for $i=1,2$
\begin{eqnarray*}
\sigma^L_i =\alpha_1(\mn_i\mn_i:\md_i)\mn_i\mn_i+\alpha_2\mn_i\mN_i+\alpha_3\mN_i\mn_i
+\alpha_4\md_i+\alpha_5\mn_i\mn_i\cdot\md_i+\alpha_6\md_i\cdot\mn_i\mn_i,
\end{eqnarray*}
and by the assumptions $|\md_i|+|\mN_i|+|(\mn_i)_t|\leq C$ in  $\R^d\times [0,T]$ with $i=1,2$. Hence,  we have
\begin{eqnarray*}
A_{6}
&\leq& -\big\{\big\langle\alpha_1(\mn_1\mn_1:\delta\md)\mn_1\mn_1+\alpha_2\mn_1(\delta\mN)+\alpha_3(\delta\mN)\mn_1
+\alpha_4\delta\md,\delta\md+\delta\momega\big\rangle\big.+\\
&&\big.\langle\alpha_5\mn_1\mn_1\cdot\delta\md+\alpha_6\delta\md\cdot\mn_1\mn_1,\delta\md+\delta\momega\rangle\big\}+\frac{\gamma}{16Re}\|\nabla\delta\mv\|_{L^2}^2+C\|\delta\mn\|_{L^2}^2\\
&\doteq&A_8+\frac{\gamma}{16(1-\gamma)}\|\nabla\delta\mv\|_{L^2}^2+C\|\delta\mn\|_{L^2}^2.
\end{eqnarray*}
Then by the same arguments as in Proposition 2.1 for (\ref{2.2})
\begin{eqnarray*}
A_{8}
&\leq& -(\alpha_1+\frac{\gamma_2^2}{\gamma_1})\|\mn_1\mn_1:\delta\md\|_{L^2}^{2}-(\alpha_5+\alpha_6-\frac{\gamma_2^2}{\gamma_1})\|\delta\md\cdot\mn_1\|_{L^2}^2-\alpha_4\|\delta\md\|_{L^2}^{2}\\
&& -\langle \delta\mh,\delta\momega\cdot\mn_1\rangle-\frac{\gamma_2}{\gamma_1}\langle\delta\mh\times\mn_1,(\delta\md\cdot\mn_1)\times\mn_1\rangle\\
&&+\frac{\gamma}{16(1-\gamma)}\|\nabla\delta\mv\|_{L^2}^2+C\|\delta\mn\|_{L^2}^2.
\end{eqnarray*}
Hence by Remark 2.2 we have
\ben\label{eq:v1-v2}
&&\frac12\frac{d}{dt}\|\delta\mv\|_{L^2}^{2}+\frac{3\gamma}{4Re}\|\nabla\delta\mv\|_{L^2}^{2}\nonumber\\
&&\leq -\frac{1-\gamma}{Re}\langle \delta\mh,\delta\momega\cdot\mn_1\rangle-\frac{\gamma_2}{\gamma_1}(\frac{1-\gamma}{Re})\langle\delta\mh\times\mn_1,(\delta\md\cdot\mn_1)\times\mn_1\rangle
+C\big(\|\delta\mn\|_{L^2}^2+\|\nabla\delta\mn\|_{L^2}^{2}\big).\nonumber\\
\een

At last, to control the first two terms of the above inequlity, we introduce the functional $\bar{W}(\mn_1,\nabla(\mn_1-\mn_2))$, and
\begin{eqnarray*}
&&\bar{W}(\mn_1,\nabla(\mn_1-\mn_2))\\
&&=a|\nabla(\mn_1-\mn_2)|^2+(k_1-a)({\rm div}(\mn_1-\mn_2))^2+(k_2-a)|\mn_1\times({\mathbf{\nabla}}\times(\mn_1-\mn_2))|^2
\\&&\quad+(k_3-a)|\mn_1\cdot(\mathbf{\nabla}\times(\mn_1-\mn_2))|^2\\
&&=a|\nabla\delta\mn|^2+(k_1-a)|{\rm div}\delta\mn|^2+(k_2-a)|\mn_1\times(\nabla\times\delta\mn)|^2+(k_3-a)|\mn_1\cdot(\nabla\times\delta\mn)|^2.
\end{eqnarray*}
Moreover, making the same computations as in Lemma 2.3 we get
\begin{eqnarray*}
\nabla_{\alpha}\bar{W}_{p_{\alpha}^{l}}-\bar{W}_{n^{l}}&=&2a\Delta\delta\mn+2(k_1-a)\nabla{\rm div}\delta\mn -2(k_2-a){\mathbf curl}(\mn_1\times({\mcurl\delta\mn}\times\mn_1)\nonumber\\
&&\quad-2(k_3-a){\mathbf curl}({\mathbf curl}\delta\mn\cdot \mn_1\mn_1) +A_9,
\end{eqnarray*}
where $A_9$ is the term without $\partial_{ij}\mn_1$, and $|A_9|\leq C(|\nabla\delta\mn|+|\delta\mn|).$ Hence
\ben\label{eq:h-w}
|\nabla_{\alpha}\bar{W}_{p_{\alpha}^{l}}-\bar{W}_{n^{l}}-(\mh_1-\mh_2)|\leq C(|\nabla\delta\mn|+|\delta\mn|).
\een
Then
\begin{eqnarray*}
&&\frac{d}{dt}\int_{\R^d}\bar{W}(\mn_1,\nabla(\mn_1-\mn_2))dx\\
&&=\int_{\R^{d}}\bar{W}_{\mn^{l}}(\mn^{l}_1-\mn^{l}_2)_{t}+\bar{W}_{p_{i}^{l}}\partial_t\nabla_i (\mn_1^l-\mn_2^l)dx+\int_{\R^{d}}\bar{W}_{\mn^{l}}(\mn^{l}_2)_{t}dx\nonumber\\
&&\leq\langle\bar{W}_{\mn^{l}}-\nabla_{i}\bar{W}_{p_{i}^{l}}, \delta\mn^l_{t}+\mv_1\cdot\nabla \mn^{l}_1-\mv_2\cdot\nabla \mn^{l}_2-(\mv_1\cdot\nabla \mn^{l}_1-\mv_2\cdot\nabla \mn^{l}_2)\rangle+C\|\nabla\delta\mn\|_{L^2}^{2}\\
&&\doteq A_{10}-\langle\bar{W}_{\mn^{l}}-\nabla_{i}\bar{W}_{p_{i}^{l}},\mv_1\cdot\nabla \mn^{l}_1-\mv_2\cdot\nabla \mn^{l}_2\rangle+C\|\nabla\delta\mn\|_{L^2}^2,
\end{eqnarray*}
while
\begin{eqnarray*}
&&-\langle\bar{W}_{\mn^{l}}-\nabla_{i}\bar{W}_{p_{i}^{l}}),\mv_1\cdot\nabla \mn^{l}_1-\mv_2\cdot\nabla \mn^{l}_2\rangle\\
&&\leq |\int_{\R^d}\bar{W}_{p_{i}^{l}}\nabla_{i}\mv_1^k\nabla_{k}(\mn^l_1-\mn^l_2)dx|+|\int_{\R^d}(\bar{W}_{\mn^{l}}-\nabla_{i}\bar{W}_{p_{i}^{l}})(\delta\mv)\cdot\nabla \mn^{l}_2dx|+C\|\nabla(\delta\mn)\|_{L^2}^2\\
&&\leq\frac{\gamma}{16(1-\gamma)}\|\nabla\delta\mv\|_{L^2}^2+C\big(\|\delta\mn\|_{L^2}^{2}+\|\delta\mv\|_{L^2}^2+\|\nabla\delta\mn\|_{L^2}\big).
\end{eqnarray*}
To estimate the term $A_{10}$,  since $|\nabla^2\mv|+|\nabla^3\mn_i|\leq C$ for $i=1,2$, by (\ref{eq:h-w}) we have
\begin{eqnarray*}
A_{10}
&\leq&\langle\bar{W}_{\mn^{l}}-\nabla_{i}\bar{W}_{p_{i}^{l}}, -\delta\Omega\cdot\mn_1+\mu_1\mn_1\times(\delta\mh\times\mn_1)
+\mu_2\mn_1\times(\delta\md\cdot\mn_1\times\mn_1)\rangle\\
&&\quad +C\big(\|\delta\mn\|_{L^2}^2+\|\nabla\delta\mn|_{L^2}^2\big)\\
&\leq&\langle\delta\mh, \delta\momega\cdot\mn_1\rangle-\mu_1\langle\delta\mh\times\mn_1,\delta\mh\times\mn_1\rangle
-\mu_2\langle\delta\mh\times\mn_1,\delta\md\cdot\mn_1\times\mn_1\rangle\\
&&+\frac{\gamma}{16(1-\gamma)}\|\nabla\delta\mv\|_{L^2}^{2}+C\big(\|\delta\mn\|_{L^2}^2+\|\nabla\delta\mn|_{L^2}^2\big).
\end{eqnarray*}
Thus,
\ben\label{eq:w1-w2}
&&\frac{d}{dt}\int_{\R^d}\bar{W}(\mn_1,\nabla(\mn_1-\mn_2))dx\nonumber\\
&&\leq\langle\delta\mh, \delta\momega\cdot\mn_1\rangle-\mu_1\langle\delta\mh\times\mn_1,\delta\mh\times\mn_1\rangle
-\mu_2\langle\delta\mh\times\mn_1,\delta\md\cdot\mn_1\times\mn_1\rangle\nonumber\\
&&+\frac{\gamma}{4(1-\gamma)}\|\nabla\delta\mv\|_{L^2}^{2}+C\big(\|\delta\mn\|_{L^2}^2+\|\delta\mv\|_{L^2}^2+\|\nabla\delta\mn|_{L^2}^2\big).
\een
Then, combine the above all estimates (\ref{eq:n1-n2}), (\ref{eq:v1-v2}) and (\ref{eq:w1-w2}), and noting that $\mu_2=-\frac{\gamma_2}{\gamma_1}$ we obtain
\begin{eqnarray*}
&&\frac{d}{dt}\big(\|\delta\mn\|_{L^2}^2+\|\delta\mv\|_{L^2}^{2}
+\frac{1-\gamma}{Re}\int_{\R^d}\bar{W}(\mn_1,\nabla(\mn_1-\mn_2))dx\big)\\
&&\leq C\big(\|\delta\mn\|_{L^2}^2+\|\delta\mv\|_{L^2}^2+\int_{\R^d}\bar{W}(\mn_1,\nabla(\mn_1-\mn_2))dx\big),\\
\end{eqnarray*}
which complete the proof by Gronwall's inequality.

 %------------------------------------------SECTION 3.3 blow up-------------------------------------------------------------
\subsection{Blow up criterion}
In this subsection, we will prove a blow up criterion for strong solution $(\mv,\mn)$  of (\ref{eq:EL}) constructed in Section 3.1 under the assumption (\ref{Parodi})-(\ref{beta}) with the data $\nabla\mn_0\in H^{2s}(\R^d)$ and $\mv_0\in H^{2s} (\R^d)(d=2\, \textrm{or}\, 3)$.
Let  $T^{*}$ be the maximal existence time of the solution. If $T^{*}<+\infty$, then it is necessary to hold that
$$
\int_{0}^{T^{*}}\|\nabla\times \mv(t)\|_{L^{\infty}}+\|\nabla\mn(t)\|_{L^{\infty}}^{2}dt =+\infty.
$$
Recall that $|\mn|=1$, $\mn\times(\Delta\mn\times\mn)=\Delta\mn+|\nabla\mn|^2\mn$, and we can obtain much precise a priori estimates than Section 3.1. In fact, we are aimed at the following energy estimates,
\ben\label{eq:blow-up estimate}
\frac{d}{dt}E_{s}(\mv,\mn)\le C(1+\|\nabla\mn\|_{L^{\infty}}^{2}+\|\nabla\mv\|_{L^{\infty}})E_{s}(\mv,\mn),
\een
where
\begin{eqnarray*}
E_{s}(\mv,\mn)&=&\|\mn-\mn_0\|_{L^2}^2+\frac{Re}{2(1-\gamma)}\|\mv\|_{L^2}^{2}+\int_{\R^d}W(\mn,\nabla\mn)dx\\
&&+a\|\Delta^s\nabla\mn\|_{L^2}^2
+(k_1-a)\|\Delta^s{\rm div}\mn\|_{L^2}^2+(k_2-a)\|\mn\times\Delta^s\mcurl\mn\|_{L^2}^{2}\\
&&+(k_3-a)\|\mn\cdot\Delta^s\mcurl\mn\|_{L^2}^{2}
+\frac{Re}{2(1-\gamma)}\|\Delta^s\mv\|_{L^2}^{2}.
\end{eqnarray*}
Then by Logarithmic Sobolev inequality in \cite{bkm}
$$
\|\nabla\mv\|_{L^{\infty}}\le C\big(1+\|\nabla\mv\|_{L^2}+\|\nabla\times\mv\|_{L^{\infty}}\log(3+\|\mv\|_{H^k})\big),
$$
for any $k\ge 3$, we have
$$
\frac{d}{dt}E_{s}(\mv,\mn)\le C(1+\|\nabla\mn\|_{L^{\infty}}^{2}+\|\nabla\mv\|_{L^2}+\|\nabla\times\mv\|_{L^{\infty}})\log(3+E_{s}(\mv,\mn))E_{s}(\mv,\mn).
$$
Applying Gronwall's inequality to the above inequality,
$$
E_{s}(\mv,\mn)\le (3+E_{s}(\mv_0,\mn_0))^{\exp\left(C\int_{0}^{t}(1+\|\nabla\mv\|_{L^2}
+\|\nabla\mn\|_{L^{\infty}}^{2}+\|\nabla\times\mv\|_{L^{\infty}})d\tau\right)},
$$
for any $t\in [0,T^{*})$. Hence we complete the proof if (\ref{eq:blow-up estimate}) holds.
In order to obtain (\ref{eq:blow-up estimate}), we sketch the proof since it's similar to the arguments in Section 3.1, and  we divide it into three steps.

{\bf Step 1.} Estimate the lower order terms.
As in Proposition \ref{energy}, the energy law holds
\ben\label{energy3.3}
&&\frac{d}{dt}\int_{\R^d}\frac{Re}{2(1-\gamma)}|\mv|^2+W(\mn,\nabla\mn)dx
+\frac{\gamma}{1-\gamma}\int_{\R^d}|\nabla\mv|^2dx+\frac{1}{\gamma_1}\int_{\R^d}|\mn\times \mh|^2dx\nonumber\\
&&=-\int_{\R^d} \beta_1|\md:\mn\mn|^2+\beta_2\md:\md+\beta_3|\md\cdot\mn|^2dx
\le 0.
\een
Using (\ref{eq:equal n}) and  (\ref{h}), by integration by parts we have
\ben\label{n3.3}
&&\frac{d}{dt}\|\mn-\mn_0\|_{L^2}^{2}=2\langle\partial_t\mn,\mn-\mn_0\rangle\nonumber\\
&&=-2\langle\mv\cdot\nabla\mn_{0}+\mn\times\big((\Omega\cdot\mn-\mu_1\mh-\mu_2\md\cdot\mn)\times \mn\big),\mn-\mn_0\rangle\nonumber\\
&&\quad-2\langle\mv\cdot\nabla (\mn-\mn_0), \mn-\mn_0\rangle\nonumber\\
&&\le C\left(\|\nabla\mn_0\|_{L^{2}}\|\mv\|_{L^2}+
\|\nabla\mv\|_{L^2}\|(\mn-\mn_0)\|_{L^2}+\|\nabla\mn\|_{L^{2}}^2+\|\nabla(\mn-\mn_0)\|_{L^2}\|\nabla\mn\|_{2}\right)\nonumber\\
&&\le CE_{s}(\mv,\mn),
\een
where we have used $\nabla\mn\in C\big([0,T^*);H^{2s}(\R^d)\big)$.

{\bf Step 2.} Estimates the higher order derivatives of $\mn$. As in Section 3.1, we have
\ben
&&\frac12\frac{d}{dt}\langle\nabla\Delta^s\mn,\nabla\Delta^s\mn\rangle\nonumber\\
&&=-\langle\nabla\Delta^{s}(\mv\cdot\nabla\mn),\Delta^{s}\nabla\mn\rangle
+\langle\Delta^{s}\left(\mn\times((\momega\cdot\mn)\times\mn)\right),\Delta^{s+1}\mn\rangle\nonumber\\
&&\quad-\mu_2\langle\Delta^{s}(\mn\times\left((\md\cdot\mn)\times\mn\right)),\Delta^{s+1}\mn\rangle
-\mu_1\langle\Delta^{s}(\mn\times(\mh\times\mn)),\Delta^{s+1}\mn\rangle\nonumber\\
&&\doteq II_1+II_2+II_3+II_4,\label{highnn3.3}
\een
and
\ben
&&\frac12\frac{d}{dt}\langle{\rm div}\Delta^s\mn,\rm div\Delta^s\mn\rangle\nonumber\\
&&=-\langle\Delta^{s}{\rm div}(\mv\cdot\nabla\mn),\Delta^{s}{\rm div}\mn\rangle+\langle\Delta^{s}(\mn\times\left((\momega\cdot\mn)\times\mn\right)),\Delta^{s}\nabla{\rm div}\mn\rangle\nonumber\\
&&\quad-\mu_2\langle\Delta^{s}(\mn\times\left((\md\cdot\mn)\times\mn\right)),\Delta^{s}\nabla{\rm div}\mn\rangle-\mu_1\langle\Delta^{s}(\mn\times(\mh\times\mn)),\Delta^{s}\nabla{\rm div}\mn\rangle\nonumber\\
&&\doteq II_1'+II_2'+II_3'+II_4'\label{highdn3.3}
\een
Moreover, $\langle\mn\times\Delta^s{\mcurl}\mn,\mn\times\Delta^s\mcurl\mn\rangle
=\langle\Delta^s{\mcurl}\mn,\Delta^s\mcurl\mn\rangle-\langle\mn\cdot\Delta^s\mcurl\mn,\mn\cdot\Delta^s\mcurl\mn\rangle$,
and
\ben
&&\frac12\frac{d}{dt}\langle\Delta^s{\mcurl}\mn,\Delta^s\mcurl\mn\rangle\nonumber\\
&&=-\langle\mcurl \Delta^{s}(\mv\cdot\nabla\mn),\Delta^{s}\mcurl\mn\rangle-\langle\Delta^{s}
\big(\mn\times((\momega\cdot\mn)\times\mn)\big),\Delta^{s}{\mcurl\mcurl}\mn\rangle\nonumber\\
&&+\mu_2\langle\Delta^{s}\big(\mn\times((\md\cdot\mn)\times\mn)\big),
\Delta^{s}\mcurl\mcurl\mn\rangle+\mu_1\langle\Delta^{s}(\mn\times(\mh\times\mn)),
\Delta^{s}\mcurl\mcurl\mn\rangle\nonumber\\
&&\doteq II_1''+II_2''+II_3''+II_4''\label{highcn3.3}
\een
On the other hand,
\beno
&&\frac12\frac{d}{dt}\langle\mn\cdot\Delta^s\mcurl\mn,\mn\cdot\Delta^s\mcurl\mn\rangle\nonumber\\
&&=\langle\partial_t\mn\cdot\Delta^s\mcurl\mn+\mn\cdot\Delta^s\mcurl\partial_t\mn, \mn\cdot\Delta^s\mcurl\mn\rangle\nonumber\\
&&=-\langle(\mv\cdot\nabla\mn)\cdot\Delta^s\mcurl\mn+\mn\cdot\Delta^{s}\mcurl(\mv\cdot\nabla\mn),
\mn\cdot\Delta^s\mcurl\mn\rangle\nonumber\\
&&\quad+\langle(\partial_t \mn+\mv\cdot\nabla\mn)\cdot\Delta^s\mcurl\mn+\mn\cdot\Delta^{s}\mcurl(\partial_t\mn+\mv\cdot\nabla\mn),
\mn\cdot\Delta^s\mcurl\mn\rangle\nonumber\\
&&=\langle\mn\cdot(\mv\cdot\nabla)(\Delta^s\mcurl\mn)-\mn\cdot\Delta^s\mcurl(\mv\cdot\nabla\mn), \mn\cdot\Delta^s\mcurl\mn\rangle\nonumber\\
&&\quad+\langle\Delta^s(\partial_t\mn+\mv\cdot\nabla\mn),\mcurl\left((\mn\cdot\Delta^s\mcurl\mn)\mn\right)\rangle\nonumber\\
&&\quad+\langle(\partial_t\mn+\mv\cdot\nabla\mn)\cdot\Delta^s\mcurl\mn,\mn\cdot\Delta^s\mcurl\mn\rangle\\
&&=-\langle[\Delta^{s}\mcurl,\mv\cdot]\nabla\mn,(\mn\cdot\Delta^s\mcurl\mn)\mn\rangle+\langle\Delta^s(\partial_t\mn+\mv\cdot\nabla\mn),\Delta^s\mcurl((\mn\cdot\mcurl\mn)\mn)\rangle\nonumber\\
&&\quad+\langle\Delta^s(\partial_t\mn+\mv\cdot\nabla\mn),
\mcurl\left((\mn\cdot\Delta^s\mcurl\mn)\mn\right)-\Delta^s\mcurl((\mn\cdot\mcurl\mn)\mn)\rangle\nonumber\\
&&\quad+\langle(\partial_t\mn+\mv\cdot\nabla\mn)\cdot\Delta^s\mcurl\mn,\mn\cdot\Delta^s\mcurl\mn\rangle
\eeno
which is equal to
\ben\label{highccn3.3}
&=&-\langle[\Delta^{s}\mcurl,\mv\cdot]\nabla\mn,\mn\cdot\Delta^s\mcurl\mn \mn\rangle-\langle\Delta^{s}\big(\mn\times((\momega\cdot\mn)\times\mn)\big), \Delta^s\mcurl(\mcurl\mn\cdot\mn \mn) \rangle\nonumber\\
&&+\mu_2\langle\Delta^{s}\left(\mn\times((\md\cdot\mn)\times\mn)\right),\Delta^s\mcurl(\mcurl\mn\cdot\mn \mn)\rangle\nonumber\\
&&+\mu_1
\langle\Delta^{s}(\mn\times(\mh\times\mn)),\Delta^s\mcurl(\mcurl\mn\cdot\mn \mn)\rangle\nonumber\\
&&+\langle\Delta^s(\partial_t\mn+\mv\cdot\nabla\mn),
(\mn\cdot\Delta^s\mcurl\mn)\mcurl\mn-\Delta^s((\mn\cdot\mcurl\mn)\mcurl\mn) \rangle\nonumber\\
&&+\langle\Delta^s(\partial_t\mn+\mv\cdot\nabla\mn),\mn\times\nabla\Delta^s(\mn\cdot\mcurl\mn)-\mn\times\nabla(\mn\cdot\Delta^s\mcurl\mn)\rangle\nonumber\\
&&+\langle\Delta^s(\partial_t\mn+\mv\cdot\nabla\mn),\Delta^s(\mn\times\nabla(\mn\cdot\mcurl\mn))-\mn\times\nabla\Delta^s(\mn\cdot\mcurl\mn)
\rangle\nonumber\\
&&+\langle(\partial_t\mn+\mv\cdot\nabla\mn)\cdot\Delta^s\mcurl\mn,\mn\cdot\Delta^s\mcurl\mn\rangle\nonumber\\
&\doteq& II_1'''+II_2'''+II_3'''+II_4'''+II_5'''+II_6'''+II_7'''+II_8''',
\een
where we have used the following relation, for a function $f$ and a vector field $u$,
$$
\mcurl(fu)=f\mcurl(u)+\nabla f\times u.
$$

Applying Lemma 3.1, we have
\ben\label{i13.3}
&&|II_1|+|II_1'|+|II_1''|+|II_1'''|\nonumber\\
&&\le C\big(\|\nabla\mv\|_{H^{2s}}\|\nabla\mn\|_{L^{\infty}}
+\|\nabla\mv\|_{L^{\infty}}\|\nabla\mn\|_{H^{2s}}
\big)\|\nabla\mn\|_{H^{2s}}\nonumber\\
&&\le C_{\delta}(\|\nabla\mv\|_{L^{\infty}}+\|\nabla\mn\|_{L^{\infty}}^{2})\|\nabla\mn\|_{H^{2s}}^{2}+\delta\|\nabla\mve\|_{H^{2s}}^{2},
\een
where $\delta>0$, to be decided later.

For the terms $II_2,\cdots,II_2'''$, we will use the following Gagliardo-Sobolev inequality on $\R^d$ (for example, see \cite{Ad}). Let $\alpha\in N$ and $\alpha\geq 2s-1$, then for $1\leq j\leq [\frac{\alpha}{2}]$, $[\frac{\alpha}{2}]+1\leq k\leq \alpha$ and $f\in H^{\alpha+1}(\R^d)$ we have
\ben\label{gagliardo-sobilev2}
\|\nabla^jf\|_{L^{\infty}}\leq C\|\nabla f\|_{H^{\alpha}}^{\frac{j}{\alpha+1-d/2}}\|f\|_{L^{\infty}}^{1-\frac{j}{\alpha+1-d/2}},\quad \|\nabla^kf\|_{L^{2}}\leq C\|\nabla f\|_{H^{\alpha}}^{\frac{k-d/2}{\alpha+1-d/2}}\|f\|_{L^{\infty}}^{1-\frac{k-d/2}{\alpha+1-d/2}}\nonumber\\
\een
Hence, for $\alpha\geq 2s-1$ with $s\geq 2$, the following inequalities hold
\ben\label{gagliardo-sobilev-n}
&&\|\nabla^{\alpha+1}\mn\|_{L^2}\|\nabla\mn\|_{L^{\infty}}+\|\nabla^{\alpha}\mn\|_{L^2}\|\nabla^2\mn\|_{L^{\infty}}
+\|\nabla^{\alpha}\mn\|_{L^2}\|\nabla\mn\|_{L^{\infty}}^2\leq C \|\nabla\mn\|_{H^{\alpha+1}},\nonumber\\
&&(\|\nabla^2 \mn\|_{L^{\infty}}+\|\nabla\mn\|_{L^{\infty}}^2)\|\nabla^{\alpha}\mn\|_{L^{2}}\leq C\|\nabla\mn\|_{L^{\infty}}\|\nabla^{\alpha+1}\mn\|_{L^2}.
\een
by Lemma 3.1, (\ref{h}) and the above inequalities (\ref{gagliardo-sobilev-n}) we have
\ben\label{i23.3}
&&2a II_2+2(k_1-a)II_2'+2(k_2-a)II_2''+2(k_3-k_2)II_2'''\nonumber\\
&&\le \langle \mn\times\left((\Delta^{s}\momega\cdot\mn)\times\mn\right),\Delta^{s}\mh\rangle
+|\langle\nabla[\Delta^s,\mn\times(\mn\times(\mn\cdot))]\momega,\Delta^{s-1}\nabla\mh\rangle|\nonumber\\
&&\quad +C|\langle\Delta^{s}\big(\mn\times((\momega\cdot\mn)\times\mn)\big), \Delta^s(\mcurl\mn\cdot\mn \mcurl\mn) \rangle|\nonumber\\
&&\le \langle \mn\times\left((\Delta^{s}\momega\cdot\mn)\times\mn\right),\Delta^{s}\mh\rangle\nonumber\\
&&\quad +C\big(\|\nabla^{2s}\mv\|_{L^2}\|\nabla^2\mn\|_{L^{\infty}}
+\|\nabla^{2s}\mv\|_{L^2}\|\nabla\mn\|_{L^{\infty}}^2+\|\nabla^{2s+1}\mv\|_{L^2}\|\nabla\mn\|_{L^{\infty}}\big)\|\nabla\mn\|_{H^{2s}}
\nonumber\\
&&\quad +C\big(\|\nabla^{2s+1}\mv\|_{L^2}
+\|\nabla^{2s}\mn\|_{L^2}\|\nabla\mv\|_{L^{\infty}}\big)\|\nabla\mn\|_{H^{2s}}
\|\nabla\mn\|_{L^{\infty}}\nonumber\\
&&\le \langle \mn\times\left((\Delta^{s}\momega\cdot\mn)\times\mn\right),\Delta^{s}\mh\rangle
+C_{\delta}\big(\|\nabla\mv\|_{L^{\infty}}
+\|\nabla\mn\|_{L^{\infty}}^{2}\big)(\|\nabla\mn\|_{H^{2s}}^{2}+\|\mv\|_{H^{2s}}^{2})\nonumber\\
&&\quad+\delta(\|\nabla\mv\|_{H^{2s}}^{2}+\|\Delta^{s+1}\mn\|_{L^{2}}^{2}).
\een
 Similar arguments hold for the terms $II_3,\cdots,II_3'''$,
\ben\label{i33}
&&2a II_3+2(k_1-a)II_3'+2(k_2-a)II_3''+2(k_3-k_2)II_3'''\nonumber\\
&&\le -\mu_2\langle \mn\times\left((\Delta^{s}\md\cdot\mn)\times\mn\right),\Delta^{s}\mh   \rangle+C_{\delta}\big(\|\nabla\mv\|_{L^{\infty}}
+\|\nabla\mn\|_{L^{\infty}}^{2}\big)(\|\nabla\mn\|_{H^{2s}}^{2}+\|\mv\|_{H^{2s}}^{2})\nonumber\\
&&\quad+\delta(\|\nabla\mv\|_{H^{2s}}^{2}+\|\Delta^{s+1}\mn\|_{L^{2}}^{2}).
\een

For the term $II_5'''$, by Lemma 3.1, (\ref{eq:equal n}) and Gagliardo-Sobolev inequality (\ref{gagliardo-sobilev-n}), we have
\ben\label{ii5'''}
&&|II_5'''|\nonumber\\
&\le& C\big(\|\Delta^s\nabla\mv\|_{L^2}+\|\Delta^s\mn\|_{L^2}\|\nabla\mv\|_{L^{\infty}}
+\|\Delta^s\mn\|_{L^2}(\|\nabla\mn\|_{L^{\infty}}^2+\|\nabla^2\mn\|_{L^{\infty}})+\|\Delta^s\mh\times\mn\|_{L^2}\big)\nonumber\\
&&\cdot
(\|\Delta^s\mn\|_{L^2}\|\nabla\mn\|_{L^{\infty}}^2+\|\nabla\mn\|_{L^{\infty}}\|\Delta^s\nabla\mn\|_{L^2}+\|\Delta^s\mn\|_{L^2}\|\nabla^2\mn\|_{L^{\infty}})\nonumber\\
&\le& C(\|\nabla\mn\|_{L^{\infty}}^2+\|\nabla\mv\|_{L^{\infty}})\|\Delta^s\nabla\mn\|_{L^2}^2+\delta (\|\Delta^s\nabla\mv\|_{L^2}+\|\Delta^{s+1}\mn\|_{L^2})^2,
\een
and similar arguments hold for $II_6'''-II_8'''$.

For the terms $II_4-II_4'''$, we have
\ba\label{ni4}
&&2aII_4+2(k_1-a)II_4'+2(k_2-a)II_4''+2(k_3-k_2)II_4'''\nonumber\\
&&=-\mu_1\langle\Delta^s(\mn\times(\mh\times\mn)),\Delta^{s}(\nabla_{\alpha}W_{p_{\alpha}^{l}})\rangle\nonumber\\
&&=-\mu_1\langle\Delta^s\left(\mn\times((\mh-\nabla_{\alpha}W_{p_{\alpha}^{k}})\times\mn\right),\Delta^s(\nabla_\alpha W_{p_{\alpha}^{l}})\rangle\nonumber\\
&&\quad-\mu_1\langle\Delta^s\left(\mn\times(\nabla_{\alpha}W_{p_{\alpha}^{k}}\times\mn\right),\Delta^{s}(\nabla_\alpha W_{p_{\alpha}^{l}})\rangle\nonumber\\
 &&\doteq II_{41}+II_{42}.
 \ea
Clearly,
 \begin{equation}\label{ii}
 II_{41} \le C_{\delta}(\|\nabla\mn\|_{L^{\infty}}^{2}+1)\|\Delta^s\nabla\mn\|_{L^2}^{2}+\delta\|\Delta^{s+1}\mn\|_{L^2}^2,
\end{equation}
and
\ba\label{i}
II_{42}&=&-\mu_1\langle\Delta^{s}\left(\mn\times(\nabla_{\alpha}W_{p_{\alpha}^{l}}\times\mn)\right),
\Delta^s(2a\Delta \mn)\rangle\nonumber\\
&&-\mu_1\langle\Delta^{s}
\left(\mn\times((\nabla_{\alpha}W_{p_{\alpha}^{l}}-2a\Delta \mn)\times\mn)\right),
\Delta^s(\nabla_{\alpha}W_{p_{\alpha}^{l}}-2a\Delta \mn)\rangle\nonumber\\
&&-2a\mu_1\langle\Delta^{s}
\left(\mn\times(\Delta \mn\times\mn)\right),
\Delta^s(\nabla_{\alpha}W_{p_{\alpha}^{l}}-2a\Delta \mn)\rangle\nonumber\\
&\doteq&II_{43}+II_{44}+II_{45}.
\ea
Note that $\nabla\times(\nabla\times\mn)=\nabla({\rm div}\mn)-\Delta\mn$, and direct calculation shows that
\ben\label{nablawn}
&&\nabla_{\alpha}W_{p_{\alpha}^{l}}\cdot n^l\\
&=&-2k_2|\nabla\mn|^2-2(k_3-k_2)(\mn\cdot\mcurl\mn)^2-2(k_1-k_2)({\rm div}\mn)^2+2(k_1-k_2)\nabla_l( n^l {\rm div}\mn).\nonumber
\een
Thus, by (\ref{nabla w}), Lemma 3.1 and (\ref{gagliardo-sobilev-n}) we infer that
\beno
\no II_{43}&=&-2a\mu_1\langle\Delta^s\nabla_{\alpha}W_{p_{\alpha}^{l}},
\Delta^{s+1}\mn\rangle+2a\mu_1\langle\Delta^s((\nabla_{\alpha}W_{p_{\alpha}^{l}}\cdot n^l)\mn),\Delta^{s+1}\mn\rangle\nonumber\\
\no &=&-4a^2\mu_1\langle\Delta^{s+1}\mn,\Delta^{s+1}\mn\rangle-4a(k_1-a)\mu_1\langle\nabla\Delta^s{\rm div}\mn,\nabla\Delta^s{\rm div}\mn\rangle\nonumber\\
\no &&-4a(k_2-a)\mu_1\langle\nabla\Delta^s\mcurl\mn,\nabla\Delta^s\mcurl\mn\rangle\nonumber\\
\no &&-4a(k_3-k_2)\mu_1\langle\nabla_l\Delta^s(\mn\cdot\mcurl\mn),\mn\cdot\nabla_l\Delta^{s}\mcurl\mn\rangle\nonumber\\
\no &&-4a(k_3-k_2)\mu_1\langle [\nabla_l\Delta^s,\mn](\mn\cdot\mcurl\mn),\nabla_l\Delta^{s}\mcurl\mn\rangle\nonumber\\
\no &&+2a\mu_1\langle\Delta^s((\nabla_{\alpha}W_{p_{\alpha}^{l}}\cdot n^l)\mn),\Delta^{s+1}\mn\rangle\nonumber\\
\no &\le& -4a^2\mu_1\langle\Delta^{s+1}\mn,\Delta^{s+1}\mn\rangle-4a(k_1-a)\mu_1\langle\nabla\Delta^s{\rm div}\mn,\nabla\Delta^s{\rm div}\mn\rangle\nonumber\\
\no &&-4a(k_2-a)\mu_1\langle\nabla\Delta^s\mcurl\mn,\nabla\Delta^s\mcurl\mn\rangle\nonumber\\
\no &&-4a(k_3-k_2)\mu_1\langle\mn\cdot\nabla_l\Delta^s(\mcurl\mn),\mn\cdot\nabla_l\Delta^{s}\mcurl\mn\rangle\nonumber\\
&&+C_{\delta}(\|\nabla\mn\|_{L^{\infty}}^{2}+1)\|\Delta^s\nabla\mn\|_{L^2}^{2}+2\delta\|\Delta^{s+1}\mn\|_{L^2}^2,
\eeno
where for the last term of the second equality we have used (\ref{nablawn}) and the following observation:
\ben\label{eq:dot}
&&\langle\Delta^s(\nabla_{l}(n^l{\rm div}\mn)\cdot n^k),\Delta^{s+1} n^k\rangle\nonumber\\
&&=\langle\Delta^s\nabla_{l}(n^l{\rm div}\mn),n^k\Delta^{s+1}n^{k}\rangle+\langle[\Delta^s,n^k]\nabla_{l}(n^l{\rm div}\mn),\Delta^{s+1}n^{k}\rangle\nonumber\\
&&=\langle\Delta^s\nabla_{l}(n^l{\rm div}\mn),\Delta^{s}(n^k \Delta n^{k})\rangle-\langle\Delta^s\nabla_{l}(n^l{\rm div}\mn),[\Delta^{s},n^k] \Delta n^{k}\rangle\nonumber\\
&&\quad+\langle[\Delta^s,n^k]\nabla_{l}(n^l{\rm div}\mn),\Delta^{s+1}n^{k}\rangle,
\een
which follows from Lemma 3.1, (\ref{gagliardo-sobilev-n}) and $n^k \Delta n^{k}=-|\nabla\mn|^2$ that
\beno
\langle\Delta^s(\nabla_{l}(n^l{\rm div}\mn)\cdot n^k),\Delta^{s+1} n^k\rangle\leq C_{\delta}(\|\nabla\mn\|_{L^{\infty}}^{2}+1)\|\Delta^s\nabla\mn\|_{L^2}^{2}+\delta\|\Delta^{s+1}\mn\|_{L^2}^2.
\eeno
Due to $a=\min\{k_1,k_2,k_3\}$ and $k_2-a\geq k_2-k_3$ when $k_2>k_3$, thus we have
\ben\label{ii1}
II_{43}\leq -4a^2\mu_1\langle\Delta^{s+1}\mn,\Delta^{s+1}\mn\rangle
+C_{\delta}(\|\nabla\mn\|_{L^{\infty}}^{2}+1)\|\Delta^s\nabla\mn\|_{L^2}^{2}+2\delta\|\Delta^{s+1}\mn\|_{L^2}^2.
\een
Similarly, we may obtain that
\begin{eqnarray}\label{ii4}
\no II_{44}
\no &\le&-\mu_1\langle \mn\times\Delta^s(\nabla_{\alpha}W_{p_{\alpha}^{l}}-2a\Delta \mn) ,\mn\times\Delta^s(\nabla_{\alpha}W_{p_{\alpha}^{k}}-2a\Delta \mn)\rangle\nonumber\\
&&+C_{\delta}(\|\nabla\mn\|_{L^{\infty}}^{2}+1)\|\Delta^s\nabla\mn\|_{L^2}^{2}+\delta\|\Delta^{s+1}\mn\|_{L^2}^2,
\end{eqnarray}
and
\begin{eqnarray}\label{ii3}
II_{45}&\le&-4a(k_1-a)\mu_1\langle\nabla\Delta^s{\rm div}\mn,\nabla\Delta^s{\rm div}\mn\rangle-4a(k_2-a)\mu_1\langle\nabla\Delta^s\mcurl\mn,\nabla\Delta^s\mcurl\mn\rangle\nonumber\\
&&-4a(k_3-k_2)\mu_1\langle\mn\cdot\nabla_l\Delta^s(\mcurl\mn),\mn\cdot\nabla_l\Delta^{s}\mcurl\mn\rangle\nonumber\\
&&+C_{\delta}(\|\nabla\mn\|_{L^{\infty}}^{2}+1)\|\Delta^s\nabla\mn\|_{L^2}^{2}+\delta\|\Delta^{s+1}\mn\|_{L^2}^2\nonumber\\
&\leq& C_{\delta}(\|\nabla\mn\|_{L^{\infty}}^{2}+1)\|\Delta^s\nabla\mn\|_{L^2}^{2}+\delta\|\Delta^{s+1}\mn\|_{L^2}^2.
\end{eqnarray}

{\bf Step 3.} Estimates of higher order derivatives of $\mv$. Obviously we only need to consider the terms $III_2-III_5$.
Especially, Lemma 3.1 and (\ref{gagliardo-sobilev-n}) yield that
\ben\label{iii2}
III_2\le C_{\delta}(1+\|\nabla\mn\|_{L^{\infty}}^{2})\|\nabla\mn\|_{H^{2s}}^{2}
+\delta\|\nabla\mv\|_{H^{2s}}^2,
\een
\begin{eqnarray}\label{iii4}
III_4\le \mu_2\langle\mn\times(\Delta^s\mh\times\mn),\Delta^s\md\cdot\mn\rangle+C_{\delta}(1+\|\nabla\mn\|_{L^{\infty}}^2)\|\nabla\mn\|_{H^{2s}}^2
+\delta\|\nabla\mv\|_{H^{2s}}^{2},
\end{eqnarray}
\begin{eqnarray}\label{iii5}
III_5\le -\langle\mn\times(\Delta^s\mh\times\mn),\Delta^s\momega\cdot\mn\rangle
+C_{\delta}(1+\|\nabla\mn\|_{L^{\infty}}^2)\|\nabla\mn\|_{H^{2s}}^2
+\delta\|\nabla\mv\|_{H^{2s}}^{2}.
\end{eqnarray}
Now we estimate one term of $III_3$. Since Lemma 3.1 and (\ref{gagliardo-sobilev-n}) imply that
\beno
&&|\langle\nabla_j[\Delta^s,\mn^j\mn\cdot]\nabla\mv,\Delta^s\mv\rangle|\\
&&\leq C\big(\|\nabla^2\mn\|_{L^{\infty}}\|\nabla^{2s}\mv\|_{L^2}+\|\nabla\mv\|_{L^{\infty}}\|\nabla\mn\|_{H^{2s}}
+\|\nabla\mn\|_{L^{\infty}}\|\nabla\mv\|_{H^{2s}}\big)\|\mv\|_{H^{2s}},
\eeno
while for the first term, by the inequality  (\ref{gagliardo-sobilev2}) we get
\beno
\|\nabla^2\mn\|_{L^{\infty}}\|\nabla^{2s}\mv\|_{L^2}&\leq & C(\|\nabla\mn\|_{L^{\infty}}\|\nabla^{2s+1}\mv\|_{L^2})^{1-\frac{1}{2s-\frac {d}{2}}}(\|\nabla\mv\|_{L^{\infty}}\|\nabla^{2s+1}\mn\|_{L^2})^{\frac{1}{2s-\frac{d}{2}}}\\
&\leq &C\big(\|\nabla\mn\|_{L^{\infty}}\|\nabla^{2s+1}\mv\|_{L^2}+\|\nabla\mv\|_{L^{\infty}}\|\nabla^{2s+1}\mn\|_{L^2}\big)
\eeno
Hence we have
\beno
&&|\langle\nabla_j[\Delta^s,\mn^j\mn\cdot]\nabla\mv,\Delta^s\mv\rangle|\\
&&\leq C_{\delta}(\|\nabla\mv\|_{L^{\infty}}+\|\nabla\mn\|_{L^{\infty}}^{2})(\|\mv\|_{H^{2s}}^{2}+\|\nabla\mn\|_{H^{2s}}^{2})
+\delta\|\nabla\mv\|_{H^{2s}}^2
\eeno
Consequently, it's not difficult to obtain the estimates of $III_3$,
\ben\label{iii3}
III_3\le C_{\delta}(\|\nabla\mv\|_{L^{\infty}}+\|\nabla\mn\|_{L^{\infty}}^{2})(\|\mv\|_{H^{2s}}^{2}+\|\nabla\mn\|_{H^{2s}}^{2})
+\delta\|\nabla\mv\|_{H^{2s}}^2.
\een

At last, by choosing $\delta$ sufficiently small,  it is concluded from (\ref{energy3.3})-(\ref{iii3}) that
$$
\frac{d}{dt}E_{s}(\mv,\mn)\le C(1+\|\nabla\mn\|_{L^{\infty}}^{2}+\|\nabla\mv\|_{L^{\infty}})E_{s}(\mv,\mn).
$$
Hence the proof is complete.
%On the other hand, if $k_3\leq k_2$, we consider the formula of $\nabla_{\alpha}W_{p_{\alpha}^{l}}$ as follows:
%\beno
%(\nabla_{\alpha}W_{p_{\alpha}^{l}})&=&2a\Delta\mn+2(k_1-a)\nabla{\rm div}\mn
%-2(k_2-k_3){\mathbf curl}{(\mn\times(\mathbf curl\mn\times\mn))}\nonumber\\
%&&-2(k_3-a){\mathbf curl}({\mathbf curl}\mn\cdot \mn\mn),\nonumber
%\eeno
%similar arguments yield that the bound of $II_{43}$ and $II_{45}$.

%------------------------------------------------Sec 4 Proof of Theorem 1.2-------------------------------------------------
\section{global existence of weak solution}
In this section, we prove global existence of weak solutions of (\ref{eq:EL}) in $\R^2$. Firstly we derive higher regularity estimates and local monotonicity inequality under the condition that local energy is uniformly small, where we follow the basic spirit of Struwe \cite{St1} which is later developed by Hong-Xin in \cite{HX}. Finally, we conclude the global existence by local existence in Section 3 and a priori estimates in this section.

For two constants $\tau$ and $T$ with $0\leq \tau<T$, we denote
$$
\begin{array}{ll}
V(\tau, T):&=\{\mn:\R^2\times [\tau,T]\rightarrow S^2|~\mn\mbox{ is measurable and satisfies }\\
&~~~~~~~{\rm{esssup}}_{\tau\le t\le T}\int_{\R^2}|\nabla \mn(\cdot,t)|^2dx+\int_{\tau}^{T}\int_{\R^2}|\nabla^2\mn|^2+|\partial_t \mn|^2dxdt<\infty\},
\end{array}
$$
and
$$
\begin{array}{ll}
H(\tau, T):&=\{\mv:\R^2\times [\tau,T]\rightarrow \R^2|~\mv\mbox{ is measurable and satisfies }\\
&~~~~~~~{\rm{esssup}}_{\tau\le t\le T}\int_{\R^2}|\mv(\cdot,t)|^2dx+\int_{\tau}^{T}\int_{\R^2}|\nabla\mv|^2dxdt<\infty\}.
\end{array}
$$

%--------------------------------------------------------------lemma 4.1---------------------------------------------
\subsection{A priori regularity estimates}
The following technical lemma could be found in \cite{St1}.
\begin{lemma}\label{struwe}
There are constants $C$ and $R_0$ such that for any $u\in V(0,T)$ and any $R\in(0,R_0]$, we have

\ba\label{eq:4.1}
\int_{\R^2\times[0,T]}|\nabla u|^4dxdt&\le & C{\rm esssup}_{0\le t\le T, x\in\R^2}\int_{B_{R}(x)}|\nabla u(\cdot,t)|^2dx\nonumber\\
&& \cdot(\int_{\R^2\times[0,T]}|\nabla^2 u|^2+R^{-2}\int_{\R^2\times[0,T]}|\nabla u|^2dxdt).\nonumber\\
\ea
\end{lemma}
By the same proof as in \cite{St1}, we can get that there exists a constant $C_1$
such that for any $f\in H(0,T)$ and any $R>0$, it holds that
\ba\label{eq:4.2}
\int_{\R^2\times[0,T]} |f|^4dxdt  &\le&
C_1{\rm{esssup}}_{0\le t\le T, x\in\R^2}\int_{B_{R}(x)}|f(\cdot,t)|^2dx\nonumber\\
&& \cdot(\int_{\R^2\times[0,T]}
|\nabla f|^2dxdt
+R^{-2}\int_{\R^2\times[0,T]}|f|^2dxdt).
\ea

%--------------------------------------------------------lemma4.2----------------------------------------------------------
\begin{lemma}\label{lem:4.2}Assume that the Leslie coefficients satisfy (\ref{Parodi})-(\ref{beta}). Let $(\mv,\mn)\in H(0,T)\times V(0,T)$ be a solution of (\ref{eq:EL}) with initial values
$\mv_0\in L^2$ and $\mn_0\in H_{b}^{1}$. Then $\exists$ $\epsilon_{1}>0$ and $R_0>0$ such that if
$$
\rm{esssup}_{0\le t\le T,x\in R^2}\int_{B_{R}(x)}|\mv(\cdot,t)|^2+|\nabla\mn(\cdot,t)|^2dx<\epsilon_1,
\quad\quad\forall R\in (0,R_0],
$$
then
\ba\label{eq:nabla2n}
\int_{R^2\times[0,T]}|\nabla^2\mn|^2+|\nabla\mv|^2dxdt\le C(1+TR^{-2})\int_{\R^2}e(\mv_0,\mn_0)dx,
\ea
\ba\label{eq:nablan4}
\int_{R^2\times[0,T]}|\nabla\mn|^4+|\mv|^4dxdt\le C\epsilon_{1}(1+TR^{-2})\int_{\R^2}e(\mv_0,\mn_0)dx,
\ea
where $e(\mv_0,\mn_0)$ denotes
$W(\mn_0,\nabla\mn_0)+\frac{Re}{2(1-\gamma)}|\mv_0|^2.$
\end{lemma}
{\it Proof:}
Since $$
\frac{d}{dt}\int_{\R^2}W(\mn,\nabla\mn) dx=\int_{\R^2}\big(W_{n^{l}}-\nabla_{i}W_{p_{i}^{l}}\big)\cdot n_{t}^{l}dx=-\int_{\R^2}\mh\cdot\mn_t dx,$$
multiplying (\ref{eq:equal n}) with $\mh$, using (\ref{2.4}) and the definition of $\mh$ we get
 \ba\label{nh}
&&\frac{d}{dt}\int_{\R^2}W(\mn,\nabla\mn) dx+\mu_1\int_{\R^2}(\mn\times(\mh\times\mn))\cdot\mh dx\nonumber\\
&&=\int_{\R^2}\mh\cdot(\momega\cdot\mn)dx-\mu_2\int_{\R^2}(\md\cdot\mn)\cdot(\mn\times(\mh\times\mn))dx
+\int_{\R^2}\left((\mv\cdot\nabla)\mn\right)\cdot\mh dx\nonumber\\
&&=\int_{\R^2}\mh\cdot(\momega\cdot\mn)dx-\mu_2\int_{\R^2}(\md\cdot\mn)\cdot(\mn\times(\mh\times\mn))dx- \int_{\R^2}W_{p_{j}^{k}}\nabla_i n^k\nabla_j v^idx,
\ea
where $\mu_1=\frac{1}{\gamma_1} $ and $\mu_2=-\frac{\gamma_2}{\gamma_1}$.
On the other hand, $(\ref{eq:EL})_3$ implies $\mN=\frac{1}{\gamma_1}\mn\times((\mh-\gamma_2\md\cdot\mn)\times\mn)$, hence by (\ref{2.2})  and Remark 2.2, we have
\ben\label{nv}
&&\frac{\gamma}{1-\gamma}\int_{\R^2\times[0,T]}|\nabla \mv|^2dxdt\nonumber\\
 &&\le \int_{\R^2\times[0,T]}W_{p_{j}^{k}}\nabla_i n^k\nabla_j v^idxdt+\frac{1-\gamma}{2Re}\int_{\R^2}|\mv_0|^2dx\nonumber\\
&&\quad-\int_{\R^2\times [0,T]}\big((\alpha_1+\frac{\gamma_2^2}{\gamma_1})(\mn\mn:\md)^{2}
 +(\alpha_5+\alpha_6-\frac{\gamma_2^2}{\gamma_1})|\md\cdot\mn|^{2}
 +\alpha_4\md:\md\big)dxdt\nonumber\\
&&\quad-\int_{\R^2\times[0,T]}\mh\cdot(\momega\cdot\mn)dxdt
 +\mu_2\int_{\R^2\times[0,T]}(\md\cdot\mn)\cdot(\mn\times(\mh\times\mn))dxdt\nonumber\\
 &&\le  \int_{\R^2\times[0,T]}W_{p_{j}^{k}}\nabla_i n^k\nabla_j v^idxdt+\frac{1-\gamma}{2Re}\int_{\R^2}|\mv_0|^2dx\nonumber\\
&&\quad-\int_{\R^2\times[0,T]}\mh\cdot(\momega\cdot\mn)dxdt+\mu_2\int_{\R^2\times[0,T]}(\md\cdot\mn)\cdot(\mn\times(\mh\times\mn))dxdt
\een
Consequently, it follows from (\ref{nh}) and (\ref{nv}) that
\ba\label{nvv}
   \quad\mu_1\int_{\R^{2}\times[0,T]}(\mn\times(\mh\times\mn))\cdot\mh dxdt+\frac{\gamma}{1-\gamma}\int_{\R^2\times[0,T]}|\nabla \mv|^2dx
   \le E(\mv_0,\mn_0).
\ea
Now we will estimate the term $\int_{\R^{2}\times[0,T]}(\mn\times(\mh\times\mn))\cdot\mh dxdt$.
Due to $
   (W_{n^{l}})=2(k_3-k_2)({\rm curl}\mn\cdot\mn)({\rm curl\mn})
   $ and $|\mn|=1$, we get
 \begin{eqnarray*}
(\mn\times(\mh\times\mn))\cdot\mh
   %&=(\nabla_{\alpha}W_{p_{\alpha}^{l}}-\nabla_{\alpha}W_{p_{\alpha}^{i}}n^jn^{l}) \cdot %h^l-(W_{n^{l}}-W_{n^{j}}n^{j}n^{l})\cdot h^{l}\\
   &=&\left(\mn\times(\nabla_{\alpha}W_{p_{\alpha}^{l}}\times\mn)\right)\cdot\mh
   -\left(\mn\times(W_{n^{l}}\times\mn)\right)\cdot\mh\nonumber\\
   &=&\mn\times(\nabla_{\alpha}W_{p_{\alpha}^{l}}\times\mn)\cdot\mh-\left(W_{n^{l}}-W_{n^{l}}\cdot\mn \mn\right)\cdot\mh\nonumber\\
   &=&\mn\times(\nabla_{\alpha}W_{p_{\alpha}^{l}}\times\mn)\cdot\mh-2(k_3-k_2)({\mcurl}\mn\cdot\mn)({\mcurl\mn}\cdot\mh)\nonumber\\     &&+2(k_3-k_2)({\mcurl}\mn\cdot\mn)^2\mn\cdot\mh.
\end{eqnarray*}

By Lemma \ref{decomposition} and integrating by parts we have the following estimates:
\begin{eqnarray*}
&&\int_{\R^2\times[0,T]}\Delta\mn\cdot(\nabla_{\alpha}W_{p_{\alpha}^{l}}-2a\Delta\mn)dxdt\\
&&\ge2\int_{\R^2\times[0,T]}\{(k_1-a)|\nabla{\rm div}\mn|^2+(k_2-a)|\nabla(\mcurl\mn\times\mn)|^2+(k_3-a)|\nabla(\mcurl\mn\cdot\mn)|^2\}dxdt\nonumber\\
  &&\quad-C\int_{\R^2\times[0,T]}|\nabla\mn|^2(|\nabla^2\mn|+|\nabla\mn|^2)dxdt
\end{eqnarray*}
Using Lemma \ref{decomposition} again, $\mn\cdot\Delta\mn=-|\nabla\mn|^2$ and the above estimates, we derived
   \ba\label{nhn}
   &&\int_{\R^2\times[0,T]}(\mn\times(\mh\times\mn))\cdot\mh dxdt\nonumber\\
   &&\ge \int_{R^2\times[0,T]}(\mn\times(\nabla_{\alpha}W_{p_{\alpha}^{l}}\times\mn))\nabla_{\alpha}W_{p_{\alpha}^{l}} dxdt
 -C\int_{\R^2\times[0,T]}|\nabla\mn|^2(|\nabla^2\mn|+|\nabla\mn|^2)dxdt\nonumber\\
&&\ge 2a\int_{\R^2\times[0,T]}\nabla_{\alpha}W_{p_{\alpha}^{l}} \cdot\Delta\mn dxdt +2a\int_{\R^2\times[0,T]}\Delta\mn\cdot(\nabla_{\alpha}W_{p_{\alpha}^{l}}-2a\Delta\mn)dxdt\nonumber\\
   &&\quad+\int_{\R^2\times[0,T]}(\mn\times((\nabla_{\alpha}W_{p_{\alpha}^{l}}-2a\Delta\mn)\times\mn)) \cdot(\nabla_{\alpha}W_{p_{\alpha}^{l}}-2a\Delta\mn) dxdt\nonumber\\
    &&\quad-C\int_{\R^2\times[0,T]}|\nabla\mn|^2(|\nabla^2\mn|+|\nabla\mn|^2)dxdt\nonumber\\
   & &\ge  4a\int_{\R^2\times[0,T]}\{a|\Delta\mn|^2+2(k_1-a)|\nabla{\rm div}\mn|^2+2(k_2-a)|\nabla(\mcurl\mn\times\mn)|^2\nonumber\\
  &&\quad+2(k_3-a)|\nabla(\mcurl\mn\cdot\mn)|^2\}dxdt-C\int_{\R^2\times[0,T]}|\nabla\mn|^2(|\nabla^2\mn|+|\nabla\mn|^2)dxdt\nonumber\\
   \ea
where we have used the following relationship
 \begin{eqnarray*}
    &&\int_{\R^2\times[0,T]}(\mn\times((\nabla_{\alpha}W_{p_{\alpha}^{l}}-2a\Delta\mn)\times\mn)) \cdot(\nabla_{\alpha}W_{p_{\alpha}^{l}}-2a\Delta\mn) dxdt\\
   &&= \int_{\R^2\times[0,T]}|(\nabla_{\alpha}W_{p_{\alpha}^{l}}-2a\Delta\mn)\times\mn)|^2 dxdt\geq 0.
\end{eqnarray*}
Hence by (\ref{nvv}) and (\ref{nhn}), we have
\ben\label{4.5}
&&a\mu_1 \int_{\R^2\times[0,T]}\{a|\Delta\mn|^2+(k_1-a)|\nabla{\rm div}\mn|^2+(k_2-a)|\nabla(\mcurl\mn\times\mn)|^2\}dxdt\nonumber\\
&&\quad+\int_{\R^2\times[0,T]}(k_3-a)|\nabla(\mcurl\mn\cdot\mn)|^2dxdt+\frac{\gamma}{1-\gamma}\int_{\R^2\times[0,T]}|\nabla \mv|^2dxdt\nonumber\\
&&\le C\int_{\R^2\times[0,T]}|\nabla\mn|^4 dxdt +E(\mv_0,\mn_0).
\een
Applying Lemma \ref{struwe}  and (\ref{eq:4.2}), we show that
 \begin{eqnarray*}
    &&\int_{\R^2\times[0,T]}|\nabla\mn|^4+|\mv|^4dxdt\\
   & &\le C\epsilon_1\int_{\R^2\times[0,T]}|\nabla^2\mn|^2+|\nabla\mv|^2dxdt
    +C\epsilon_1R^{-2}\int_{\R^2\times [0,T]}|\nabla\mn|^2+|\mv|^2 dxdt.
\end{eqnarray*}
Then (\ref{eq:nabla2n}) and (\ref{eq:nablan4}) follows from the above two estimates.
Thus the lemma is complete.

The following lemma is a local monotonicity inequality for strong solution $(\mv,\mn)$.
\begin{lemma}\label{lem:monotone}
Assume that the Leslie coefficients satisfy (\ref{Parodi})-(\ref{beta}). Let $(\mv,\mn)$ be a solution of (\ref{eq:EL}) with initial values $(\mv_{0},\mn_0)$ with
 $\mn_0\in V(0,T)$ and $\mv_0\in H(0,T)$. Assume that there exists  $\epsilon_1>0$ and $R_0>0$, such that
 $$
 {\rm esssup}_{x\in \R^2,0\le t\le T}\int_{B_{R_0}(x)}|\nabla\mn|^2+|\mv|^2dx<\epsilon_1.
 $$
 Then for all $s\in [0,T]$, $x_0\in \R^2$, and $R\le R_0$,
\begin{eqnarray*}
&&\int_{B_{R}(x_0)}e(\mv,\mn)(\cdot,s)dx
 +\frac{\gamma}{1-\gamma}\int_{0}^{s}\int_{B_{R}(x_0)}|\nabla \mv|^2dxdt
+\frac{1}{2\gamma_1}\int_{0}^{s}\int_{B_{R}(x_0)}|\mn\times\mh|^2dxdt\\
&& \le \int_{B_{2R}(x_0)}e(\mv_0,\mn_0)dx+C_2\frac{s^{1/2}}{R}(1+\frac{s}{R^2})^{1/2}\int_{\R^2}e(\mv_0,\mn_0)dx,
 \end{eqnarray*}
 where $C_2$ is a uniform constant.
 \end{lemma}

 {\it Proof:} Let $\phi\in C_{0}^{\infty}(B_{2R}(x_0))$  be a cut-off function with $\phi\equiv 1$ on $B_{R}(x_0)$ and
  $|\nabla\phi|\le \frac{C}{R}$, $|\nabla^2\phi|\le \frac{C}{R^2}$ for some $R\le R_0$.

  Multiply $(\ref{eq:EL})_1$ by $\phi^2\mv$, and integrate by parts
\ben\label{eq:local-v}
 && \frac{1}{2}\frac{d}{dt}\int_{\R^2}|\mv|^2\phi^2 dx+\frac{\gamma}{Re}\int_{\R^2}|\nabla\mv|^2\phi^2dx\nonumber\\
 & &=\int_{\R^2}(|\mv|^2+2p)\phi v^{l}\nabla_{l}\phi dx+\frac{\gamma}{Re}\int_{\R^2}|\mv|^2(|\nabla\phi|^2
+\phi\Delta\phi)dx\\
&&\quad+\frac{1-\gamma}{Re}\int_{\R^2}W_{p_{j}^{k}}\nabla_i n^k \nabla_jv^{i} \phi^{2}+W_{p_{j}^{k}}
  \nabla_i n^k v^{i} \nabla_j\phi^{2}dx-\frac{1-\gamma}{Re}\int_{\R^2}\msigma^{L}:\nabla(\mv\phi^{2})dx,\nonumber
\een
while similar to (4.6), we get
\beno
&&-\int_{\R^2}\msigma^{L}:\nabla(\mv\phi^2)dx+\int_{\R^2}\mv\cdot\msigma^{L}\cdot\nabla\phi^2 dx\\
&&\leq-\int_{\R^2}\big((\alpha_1+\frac{\gamma_2^2}{\gamma_1})(\mn\mn:\md)^{2}
 +(\alpha_5+\alpha_6-\frac{\gamma_2^2}{\gamma_1})|\md\cdot\mn|^{2}
 +\alpha_4\md:\md\big)\phi^2dx\nonumber\\
&&\quad-\int_{\R^2}\mh\cdot(\momega\cdot\mn)dx
 +\mu_2\int_{\R^2}(\md\cdot\mn)\cdot(\mn\times(\mh\times\mn))\phi^2dx\\
&&\leq-\int_{\R^2}\mh\cdot(\momega\cdot\mn)dx
 +\mu_2\int_{\R^2}(\md\cdot\mn)\cdot(\mn\times(\mh\times\mn))\phi^2dx
\eeno
Recall that
$$
\frac{d}{dt}\int_{\R^2}W(\mn,\nabla\mn)\phi^{2}dx=
\int_{\R^2}\phi^2\mh\cdot(\mv\cdot\nabla\mn-\mn_t-\mv\cdot\nabla\mn)dx-\int_{\R^2}W_{p_{j}^{k}}\partial_t n^k\nabla_j\phi^2dx,
$$
and as (\ref{2.4}), we can obtain
\begin{eqnarray*}
\int_{R^2}\phi^2\mh\cdot (\mv\cdot \nabla)\mn dx
=\int_{\R^2}\big(W\mv\cdot\nabla\phi^2-\phi^2W_{p_{j}^{k}}\nabla_{j}v^l\nabla_l n^k
  -W_{p_j^k}v^{l}\nabla_l n^k\nabla_j(\phi^2)\big)dx.
\end{eqnarray*}
Moreover, using (\ref{eq:equal n})
\beno
&&-\int_{\R^2}\phi^2\mh\cdot(\mn_t+\mv\cdot\nabla\mn)dx\\
&&=-\frac{1}{\gamma_1}\int_{\R^2}\mn\times((\mh-\gamma_2\md\cdot\mn)\times\mn)\cdot\mh \phi^2dx+\int_{\R^2}\phi^{2}\mh\cdot\Omega\cdot\mn dx.
\eeno
Then the above estimates together yield that
\ben\label{4.7}
&&\frac{1}{2}\frac{d}{dt}\int_{\R^2}\frac{Re}{1-\gamma}|\mv|^2\phi^2 dx+\frac{d}{dt}
   \int_{\R^2}W(\mn,\nabla\mn)\phi^2 dx\nonumber\\
&&\quad+\frac{\gamma}{1-\gamma}\int_{\R^2}|\nabla\mv|^2\phi^2 dx+\frac{1}{\gamma_1}\int_{\R^2}|\mn\times\mh|^2\phi^2\nonumber\\
&&\leq\frac{Re}{1-\gamma}\int_{\R^2}\phi(|\mv|^2+2|p|)|\mv||\nabla\phi| dx+\frac{\gamma}{1-\gamma}\int_{\R^2}|\mv|^{2}
   (|\nabla\phi|^2+\phi|\Delta\phi|)dx\nonumber\\
&&\quad-\int_{\R^2}\mv\cdot\msigma^{L}\cdot\nabla\phi^2 dx+\int_{\R^2}W\mv\cdot\nabla\phi^2dx-\int_{\R^2}W_{p_{j}^{k}}\partial_t n^k\nabla_j\phi^2dx\nonumber\\
&&\doteq B_1+B_2+B_3+B_4+B_5.
\een
Now, we estimate the following term
\begin{eqnarray*}
  B_3&=& -\int_{\R^2}\mv\cdot\msigma^{L}\cdot\nabla\phi^2 dx\\
  & =&-\int_{\R^2}\mv\nabla\phi^2:\big(\alpha_1(\mn\mn:\md)\mn\mn+\alpha_4\md+\alpha_5\mn\mn\cdot\md+\alpha_6\md\cdot\mn\mn\big) dx\\
  &&-\alpha_2\int_{\R^2}(\mv\cdot\mn)(\mN\cdot\nabla\phi^2)dx-\alpha_3\int_{\R^2} (\mv\cdot\mN)(\mn\cdot\nabla\phi^2) dx\\
   &\le& C\int_{\R^2}|\mv||\nabla\mv|\phi|\nabla\phi|dx
   -\alpha_2\int_{\R^2}(\mv\cdot\mn)(\mN\cdot\nabla\phi^2)dx-\alpha_3\int_{\R^2} (\mv\cdot\mN)(\mn\cdot\nabla\phi^2) dx
\end{eqnarray*}
Using (\ref{eq:equal n}) again,
\begin{eqnarray*}
&&\alpha_2\int_{\R^2}(\mv\cdot\mn)(\mN\cdot\nabla\phi^2)dx+\alpha_3\int_{\R^2} (\mv\cdot\mN)(\mn\cdot\nabla\phi^2) dx\\
&&= \mu_1\int_{\R^2}[\alpha_3((\mn\cdot\nabla)\phi^2)\mv+\alpha_2(\mv\cdot\mn)\nabla\phi^2]
     \cdot(\mn\times(\mh\times\mn))dx\\
&&\quad+\mu_2\int_{\R^2}[\alpha_3((\mn\cdot\nabla)\phi^2)\mv+\alpha_2(\mv\cdot\mn)\nabla\phi^2]
     \cdot(\mn\times((\md\cdot\mn)\times\mn))dx\\
&&\le C\int_{\R^2}|\mv|\phi|\nabla\phi||\mn\times\mh|dx+C\int_{\R^2} |\mv|\phi|\nabla\phi||\nabla \mv| dx.
\end{eqnarray*}
     Thus we have
\ba\label{4.8}
&&\frac{d}{dt}\int_{\R^2}\frac{Re}{2(1-\gamma)}|\mv|^2\phi^2+W(\mn,\nabla\mn)\phi^2dx
      +\frac{3\gamma}{4(1-\gamma)}\int_{\R^2}|\nabla\mv|^2\phi^2 dx+\frac{3}{4\gamma_1}\int_{\R^2}|\mn\times\mh|^2\phi^2dx\nonumber\\
     && \le C(B_1+B_2)+|B_4|+|B_5|.\nonumber
\ea
Obviously,
\beno
\int_0^s|B_2|dt\leq C\frac{s}{R^2}\int_{\R^2}e(\mv_0,\mn_0) dx.
\eeno
For the terms $B_4$
and $B_5,$ by the equation of $\mn$, (\ref{eq:nablan4}) and Proposition \ref{energy}, we have
\beno
&&\int_0^s|B_4|+|B_5|dt\\
&&\leq C\int_{\R^2\times(0,s)}\big(|\nabla\mn|^2|\mv|+|\nabla\mn||\nabla\mv|+|\nabla\mn||\mn\times\mh|\big)\phi|\nabla\phi|dxdt\\
&&\leq \frac{\gamma}{4(1-\gamma)}\int_{\R^2\times(0,s)}|\nabla\mv|^2\phi^2 dxdt+\frac{1}{4\gamma_1}\int_{\R^2\times(0,s)}|\mn\times\mh|^2\phi^2dxdt\\
&&\quad+C\frac{s}{R^2}\int_{\R^2}e(\mv_0,\mn_0) dx+C\epsilon_{1}^{1/2}\frac{s^{1/2}}{R}(1+\frac{s}{R^2})^{1/2}\int_{\R^2}e(\mv_0,\mn_0) dx
\eeno
Finally, for the first term of $B_1$, by (\ref{eq:4.2}) it's easy to obtain
\begin{eqnarray*}
      \int_{\R^2\times(0,s)}|\mv|^2|\mv||\phi||\nabla\phi| dxdt&\le& (\int_{0}^{s}\int_{\R^2}|\mv|^4dxdt)^{1/2}
      \cdot(\int_{0}^{s}\int_{\R^2}\frac{|\mv|^2}{R^2})^{1/2} \\
      &\le& C\epsilon_{1}^{1/2}\frac{s^{1/2}}{R}(1+\frac{s}{R^2})^{1/2}\int_{\R^2}e(\mv_0,\mn_0) dx
\end{eqnarray*}
      for $R\le R_0$.
Meanwhile,
$$
\int_{\R^2\times(0,s)}|p||\mv||\phi||\nabla\phi| dxdt\le (\int_{0}^{s}
\int_{\R^2}|p|^2dxdt)^{1/2}\cdot(\int_{0}^{s}\int_{\R^2}\frac{|\mv|^2}{R^2}dxdt)^{1/2}.
$$
We note that
$$
\Delta p=\frac{1-\gamma}{Re}\nabla\cdot\big(\nabla\cdot(\msigma^{E}+\msigma^{L})\big)-\partial_i\partial_j(v^iv^j).
$$
By Calder\'{o}n-Zygmund estimates, (\ref{eq:nabla2n}), (\ref{eq:nablan4}) and Proposition \ref{energy}, we have
\beno\label{4.9}
\int_{0}^{s}\int_{\R^2}|p|^2dxdt&   \le &\int_{0}^{s}\int_{\R^2}|\nabla\mn|^4
         +|\mv|^4+|\nabla\mv|^2+|\nabla^2\mn|^2dxdt\nonumber\\
         &\le& C(1+\frac{s}{R^2})\int_{\R^2}e(\mv_0,\mn_0)dx.
         \eeno
Hence the lemma is true.

%---------------------------------------------------lemma4------------------------------------------------------------
Next lemma is devoted to promote the regularity of the solution  $(\mv,\mn)$.
\begin{lemma} \label{lem:4.4}
Assume that the Leslie coefficients satisfy (\ref{Parodi})-(\ref{beta}). Let $(\mv,\mn)$  be a solution of (\ref{eq:EL})
 with the initial value $(\mv_0,\mn_0)\in L^2\times H_{b}^{1}(\R^2,S^2)$ and ${\rm div}\mv_0=0$. There are constants
 $\epsilon_1$ and $R_0>0$ such that
 $$
 {\rm esssup}_{0\le t\le T,x\in\R^2}\int_{B_{R}(x)}|\nabla\mn(\cdot,t)|^2+|\mv|^2(\cdot,t)dx<\epsilon_1,\quad\forall R\in(0,R_0].$$
Then, for all $t\in[\tau, T]$ with $\tau\in (0,T)$, it holds that
$$
\int_{\R^2}|\nabla^2\mn(\cdot,t)|^2+|\nabla\mv|^2(\cdot,t)dx+\int_{\tau}^t\int_{\R^2}|\nabla^3\mn(\cdot,s)|^2+|\nabla^2\mv|^2(\cdot,s)dxds \le C(\epsilon_1,E_0,\tau, T,\frac{T}{R^2}),
$$
where $E_0=E(\mv_0,\mn_0)=\int_{\R^2}e(\mv_0,\mn_0)dx$.
\end{lemma}

{\it Proof:} Multiplying $(\ref{eq:EL})_1$ with $\Delta \mv$, we have
\ba\label{5.1}
&&\frac{1}{2}\frac{d}{dt}\int_{\R^2}|\nabla\mv|^2dx+ \frac{\gamma}{Re}\int_{\R^2}|\Delta\mv|^2dx\nonumber\\
&&=\int_{\R^2}(\mv\cdot\nabla)\mv\cdot\Delta \mv dx-\frac{1-\gamma}{Re}\int_{\R^2}(\nabla\cdot\msigma^{E})\cdot\Delta\mv dx-\frac{1-\gamma}{Re}\int_{\R^2}(\nabla\cdot\msigma^{L})\cdot\Delta\mv dx\nonumber\\
&&\le \frac{1}{4}\frac{\gamma}{Re}\int_{\R^2}|\Delta\mv|^2dx+C\int_{\R^2}|\mv\cdot\nabla\mv|^2dx+C\int_{\R^2}(|\nabla^2\mn|^2+|\nabla\mn|^4)|\nabla\mn|^2dx\nonumber\\ &&\quad-\frac{1-\gamma}{Re}\int_{\R^2}(\nabla\cdot\msigma^{L})\cdot\Delta\mv dx.
\ea
Note that
\begin{eqnarray*}
&&-\frac{1-\gamma}{Re}\int_{\R^2}(\nabla\cdot\msigma^{L})\cdot\Delta\mv dx \\
&&=\frac{1-\gamma}{Re}\int_{\R^2}\msigma^{L}:\Delta(\md+\momega)dx\\
&&= \frac{1-\gamma}{Re}\int_{\R^2}[\alpha_1(\mn\mn:\md)\mn\mn:\Delta\md+(\alpha_2+\alpha_3)\mn\mN:\Delta\md
+(\alpha_2-\alpha_3)\mn\mN:\Delta\Omega\\
&&\quad+\alpha_4\md:\Delta\md+(\alpha_5+\alpha_6)(\mn\mn\cdot\md):\Delta\md
+(\alpha_5-\alpha_6)(\mn\mn\cdot\md):\Delta\momega] dx.
\end{eqnarray*}
Obviously, we have
$$
\alpha_1\int_{\R^2}(\mn\mn:\md)\mn\mn:\Delta\md dx\le -\alpha_1\int_{\R^2}|\mn\mn:\nabla\md|^2dx +C\int_{\R^2}|\nabla\mn||\nabla\mv||\nabla^2\mv|dx,
$$
$$
\int_{\R^2} (\mn\mn\cdot\md):\Delta\md dx\le -\int_{\R^2}|(\nabla \md\cdot\mn)|^2dx +C\int_{\R^2}|\nabla\mn||\nabla\mv||\nabla^2\mv|dx.
$$
Moreover, recall that $\gamma_1=\alpha_3-\alpha_2$ and $\gamma_2=\alpha_6-\alpha_5=\alpha_3+\alpha_2$, hence
by $(\ref{eq:EL})_3$ and the anti-symmetric property of $\Delta\momega$ we get
\begin{eqnarray*}
&&-\int_{\R^2}\gamma_1(\mn\mN:\Delta\momega)+\gamma_2(\mn\mn\cdot\md):\Delta\momega dx\\
&&=-\int_{\R^2}\mn\cdot(-\mh+\gamma_1\mN+\gamma_2\md\cdot\mn+\mh):\Delta\momega dx
=\int_{\R^2}\mh\cdot(\Delta\momega\cdot\mn)dx.
\end{eqnarray*}
Using $(\ref{eq:EL})_3$ again, we obtain
\begin{eqnarray*}
&&\gamma_2\int_{\R^2}\mn\mN:\Delta\md dx\\
&&=\frac{\gamma_2}{\gamma_1}\int_{\R^2}(\Delta\md\cdot\mn)\cdot(\mn\times(\mh\times\mn))dx-\frac{\gamma_2^2}{\gamma_1}
\int_{\R^2}(\Delta\md\cdot\mn)\cdot(\mn\times((\md\cdot\mn)\times\mn))dx\\
&&\leq \frac{\gamma_2}{\gamma_1}\int_{\R^2}(\Delta\md\cdot\mn)\cdot(\mn\times(\mh\times\mn))dx-
\frac{\gamma_2^2}{\gamma_1}\int_{\R^2}(\Delta\md\cdot\mn)\cdot(\md\cdot\mn))dx\\
&&\quad+\frac{\gamma_2^2}{\gamma_1}\int_{\R^2}(\Delta\md\cdot\mn)\cdot((\md\cdot\mn)\cdot\mn\mn)dx\\
&&\le   \frac{\gamma_2}{\gamma_1}\int_{\R^2}(\Delta\md\cdot\mn)\cdot(\mn\times(\mh\times\mn))dx
+\frac{\gamma_2^2}{\gamma_1}\int_{\R^2}(|\nabla\md\cdot\mn|^2-|\mn^i\nabla\md_{ij}\cdot\mn^j|^2)dx\\
&&\quad+C\int_{\R^2}|\nabla\mn||\nabla\mv||\nabla^2\mv|dx
\end{eqnarray*}

Thus, combining the above estimates and applying Remark 2.2, we get
\ba\label{4.13}
&&\frac{1}{2}\frac{d}{dt}\int_{\R^2}|\nabla\mv|^2dx+ \frac{\gamma}{2Re}\int_{\R^2}|\Delta\mv|^2dx\nonumber\\
&&\le C\int_{\R^2}(|\nabla^2\mn|^2+|\nabla\mv|^2
 +|\nabla\mn|^4)(|\nabla\mn|^2+|\mv|^2)dx\nonumber\\
&&\quad+\frac{1-\gamma}{Re}\int_{\R^2}\mh\cdot(\Delta\momega\cdot\mn)dx
-\frac{1-\gamma}{Re}\mu_2\int_{\R^2}(\mn\times(\mh\times\mn))\cdot(\Delta\md\cdot\mn)dx\nonumber\\
\ea
On the other hand, we can differentiate (\ref{eq:equal n}),  multiply it by $\nabla_\beta \mh$, and we get
\ba\label{4.14}
&&\frac{d}{dt}\int_{\R^2}a|\Delta \mn|^2+(k_1-a)|\nabla{\rm div}\mn|^2
+(k_2-a)|\nabla(\mn\times\mcurl\mn)|^2+(k_3-a)|\nabla(\mn\cdot\mcurl\mn)|^2dx\nonumber\\
&&\quad-\int_{\R^2}[(\nabla_\beta\mv\cdot\nabla)n^i
+(\mv\cdot\nabla)\nabla_{\beta}n^{i}]\cdot\nabla_{\beta}h^{i}dx\nonumber\\
&&\leq\int_{\R^2}\nabla_{\beta}(\momega\cdot\mn)\cdot\nabla_{\beta}\mh dx-\mu_1\int_{\R^2}\nabla_{\beta}(\mn\times(\mh\times\mn))\cdot\nabla_{\beta}\mh dx\nonumber\\
&&\quad-\mu_2 \int_{\R^2}\nabla_{\beta}(\mn\times((\md\cdot\mn)\times\mn))\cdot\nabla_{\beta}\mh dx+\int_{\R^2}\delta|\nabla\mn_t|^2+C(\delta)(|\nabla^2\mn|^2+|\nabla\mn|^4)|\nabla\mn|^2dx\nonumber\\
&&\doteq B_6+B_7+B_8+\int_{\R^2}\delta|\nabla\mn_t|^2+C(\delta)(|\nabla^2\mn|^2+|\nabla\mn|^4)|\nabla\mn|^2dx,
\ea
where $\delta>0$, to be decided later.

Direct calculation  shows
\ben
 B_6&=&-\int_{\R^2}(\Delta\momega\cdot\mn)\cdot\mh dx-\int_{\R^2}(\Delta(\momega\cdot\mn)-\Delta\momega\cdot\mn)\cdot\mh dx\nonumber\\
 &\le& -\int_{\R^2}\mh\cdot(\Delta\momega\cdot\mn) dx+\int_{\R^2}(|\nabla\mv||\nabla \mn||\nabla\mh|+|\mh||\nabla^2\mv||\nabla\mn|)dx.
\een
Note the fact that $|\mn\cdot\nabla_{\beta}\Delta\mn|\leq C|\nabla\mn||\nabla^2\mn|$, and similar estimates to (\ref{nhn}) for $B_7$ imply
\ben
B_7\le -\mu_1 a\int_{\R^2}|\nabla^3\mn|^2+C\int_{\R^2}|\nabla\mn|^2(|\nabla^2\mn|^2+|\nabla\mn|^4)dx.
\een
At last, we estimate $B_8$.
\ben
 B_8&=&-\mu_2\int_{\R^2}\{(\nabla_{\beta}\mn\times((\md\cdot\mn)\times\mn))\cdot\nabla_{\beta}\mh
 + (\mn\times((\md\cdot\mn)\times\nabla_{\beta}\mn))\cdot\nabla_{\beta}\mh\nonumber\\
 &&+(\mn\times((\md\cdot\nabla_{\beta}\mn)\times\mn))\cdot\nabla_{\beta}\mh\} dx-\mu_2\int_{\R^2} (\mn\times((\nabla_{\beta}\md\cdot\mn)\times\mn))\cdot\nabla_{\beta}\mh dx\nonumber\\
 &\le& C\int_{\R^2}|\nabla\mv||\nabla\mn||\nabla\mh|dx +\mu_2\int_{\R^2}\mn\times((\Delta\md\cdot\mn)\times\mn)\cdot\mh dx\nonumber\\
 &&+\mu_2\int_{\R^2}\{(\nabla_{\beta}\mn\times((\nabla_{\beta}\md\cdot\mn)\times\mn))
  +(\mn\times((\nabla_{\beta}\md\cdot\nabla_{\beta}\mn)\times\mn))\nonumber\\
 &&+(\mn\times((\nabla_{\beta}\md\cdot\mn)\times\nabla_{\beta}\mn))\}\cdot\mh dx\nonumber\\
  &\le& C\int_{\R^2}|\nabla\mv||\nabla\mn||\nabla\mh|+|\nabla\mn||\nabla^2\mv||\mh|dx
  +\mu_2\int_{\R^2}((\Delta\md\cdot\mn)\cdot(\mn\times(\mh \times\mn)) dx.\nonumber\\
 \een
The definition of $\mh$ and (\ref{eq:equal n}) tell us that
\ben\label{eq:4.18}
&&-\mn\cdot\Delta\mn=|\nabla\mn|^2,\quad |\mh|\leq C(|\nabla\mn|^2+|\nabla^2\mn|),\nonumber\\
&&|\nabla \mh|\leq C(|\nabla\mn|^3+|\nabla^2\mn||\nabla\mn|+|\nabla^3\mn|),
\een
and
\ben\label{eq:4.19}
 |\nabla \mn_t|^2\leq \max\{1,|\mu_1|,|\mu_2|\} (|\nabla^3 \mn|^2+|\nabla^2 \mv|^2)+C(|\mv|^2+|\nabla\mn|^2)(|\nabla\mv|^2+|\nabla^2\mn|^2).
\een
Hence combining these estimates from (\ref{4.13}) to (\ref{eq:4.19}) and choosing $\delta=\frac{a\mu_1}{2+2|\mu_1|+2|\mu_2|}$, by Gagliardo-Sobolev inequlity we obtain
 \beno
&&\frac{d}{dt} \int_{\R^2}\frac{Re}{2(1-\gamma)}|\nabla\mv|^2 +a|\Delta\mn|^2dx
  +\int_{\R^2}\frac{\gamma}{4(1-\gamma)}|\nabla^2\mv|^2 +\frac{a\mu_1}{4}|\nabla^3\mn|^2 dx\nonumber\\
&&\quad+\frac{d}{dt}\int_{\R^2}(k_1-a)|\nabla{\rm div}\mn|^2
+(k_2-a)|\nabla(\mn\times\mcurl\mn)|^2+(k_3-a)|\nabla(\mn\cdot\mcurl\mn)|^2dx\nonumber\\
&&\le C\int_{\R^2}(|\mv|^2+|\nabla\mn|^2)(|\nabla\mv|^2+|\nabla^2\mn|^2)dx,\nonumber\\
&&\le C\big(\int_{\R^2}(|\mv|^4+|\nabla\mn|^4)dx\big)^{1/2}\big(\int_{\R^2}(|\nabla\mv|^4+|\nabla^2\mn|^4)dx\big)^{1/2}\\
&&\le \min\{\frac{\gamma}{8(1-\gamma)},\frac{a\mu_1}{8}\}\int_{\R^2}|\nabla^2\mv|^2 +|\nabla^3\mn|^2 dx\\
&&\quad+C\big(\int_{\R^2}(|\mv|^4+|\nabla\mn|^4)dx\big)\big(\int_{\R^2}(|\nabla\mv|^2+|\nabla^2\mn|^2)dx\big)
\eeno
Due to (\ref{eq:nabla2n}), we know that for $\tau\in (0,T)$ there exists a $\tau'\in (0,\tau)$ such that $\int_{\R^2}(|\nabla\mv|^2+|\nabla^2\mn|^2)(\cdot,\tau')dx<C(\tau,E_0,\frac{T}{R^2}).$
Then the required inequality follows from Gronwall's inequality, (\ref{eq:nabla2n}) and (\ref{eq:nablan4}); see also, \cite[Lemma 3.4]{HX}.\endproof

%---------------------------------------------------lemma4.5------------------------------------------------------------
Similar computations to Lemma 4.4, we may obtain higher regularity for $(\mv,\mn)$ of (\ref{eq:EL}).
\begin{lemma} \label{lem:4.5}
Assume that the Leslie coefficients satisfy (\ref{Parodi})-(\ref{beta}). Let $(\mv,\mn)$  be a solution of (\ref{eq:EL})
 with the initial value $(\mv_0,\mn_0)\in L^2\times H_{b}^{1}(\R^2,S^2)$ and ${\rm div}\mv_0=0$. There are constants
 $\epsilon_1$ and $R_0>0$ such that
 $$
 {\rm esssup}_{0\le t\le T,x\in\R^2}\int_{B_{R}(x)}|\nabla\mn(\cdot,t)|^2+|\mv|^2(\cdot,t)dx<\epsilon_1,\quad\forall R\in(0,R_0].$$
Then, for all $t\in[\tau, T]$ with $\tau\in (0,T)$, it holds that
$$
\int_{\R^2}|\nabla^3\mn(\cdot,t)|^2+|\nabla^2\mv|^2(\cdot,t)dx +\int_{\tau}^t\int_{\R^2}|\nabla^4\mn(\cdot,s)|^2+|\nabla^3\mv|^2(\cdot,s)dxds\le C(\epsilon_1,E_0,\tau, T, \frac{T}{R^2}).
$$
\end{lemma}

{\it Proof:} Since $|\nabla\mn|^2=|\mn\cdot\Delta\mn|\leq |\nabla^2\mn|$. Multiplying $(\ref{eq:EL})_1$ with $\Delta^2 \mv$, we have
\ba\label{5.10}
&&\frac{1}{2}\frac{d}{dt}\int_{\R^2}|\nabla^2\mv|^2dx+ \frac{\gamma}{Re}\int_{\R^2}|\nabla\Delta\mv|^2dx\nonumber\\
&&=-\int_{\R^2}(\mv\cdot\nabla)\mv\cdot\Delta^2 \mv dx+\frac{1-\gamma}{Re}\int_{\R^2}(\nabla\cdot\msigma^{E})\cdot\Delta^2\mv dx+\frac{1-\gamma}{Re}\int_{\R^2}(\nabla\cdot\msigma^{L})\cdot\Delta^2\mv dx\nonumber\\
&&\le \frac{1}{4}\frac{\gamma}{Re}\int_{\R^2}|\nabla\Delta\mv|^2dx+C\int_{\R^2}|\nabla(\mv\cdot\nabla\mv)|^2dx+\frac{1-\gamma}{Re}\int_{\R^2}(\nabla\cdot\msigma^{L})\cdot\Delta^2\mv dx\nonumber\\
&&\quad +C\int_{\R^2}(|\nabla^3\mn|^2|\nabla\mn|^2+|\nabla^2\mn|^4)dx.
\ea
As the above lemma, we also have
\beno
&&\frac{1-\gamma}{Re}\int_{\R^2}(\nabla\cdot\msigma^{L})\cdot\Delta^2\mv dx \\
&&=-\frac{1-\gamma}{Re}\int_{\R^2}\msigma^{L}:\Delta^2(\md+\momega)dx\\
&&=- \frac{1-\gamma}{Re}\int_{\R^2}[\alpha_1(\mn\mn:\md)\mn\mn:\Delta^2\md+(\alpha_2+\alpha_3)\mn\mN:\Delta^2\md
+(\alpha_2-\alpha_3)\mn\mN:\Delta^2\Omega\\
&&\quad+\alpha_4\md:\Delta^2\md+(\alpha_5+\alpha_6)(\mn\mn\cdot\md):\Delta^2\md
+(\alpha_5-\alpha_6)(\mn\mn\cdot\md):\Delta^2\momega] dx,
\eeno
while
\beno
&&-\alpha_1\int_{\R^2}(\mn\mn:\md)\mn\mn:\Delta^2\md dx\\
&&\le -\alpha_1\int_{\R^2}|\mn\mn:\Delta\md|^2dx +C\int_{\R^2}|\nabla\mn|^2|\nabla\mv||\nabla^3\mv|+|\nabla^2\mv||\nabla^3\mv||\nabla\mn|dx,
\eeno
and
\beno
&&-\int_{\R^2} (\mn\mn\cdot\md):\Delta^2\md dx\\
&&\le -\int_{\R^2}|(\Delta \md\cdot\mn)|^2dx+C\int_{\R^2}|\nabla^2\mn||\nabla\mv||\nabla^3\mv|+|\nabla^2\mv||\nabla^3\mv||\nabla\mn|dx.
\eeno
Moreover,
by $(\ref{eq:EL})_3$ and the anti-symmetric property of $\Delta\momega$ we get
\begin{eqnarray*}
&&\int_{\R^2}\gamma_1(\mn\mN:\Delta^2\momega)+\gamma_2(\mn\mn\cdot\md):\Delta^2\momega dx\\
&&=\int_{\R^2}\mn\cdot\Delta^2\momega\cdot(-\mh+\gamma_1\mN+\gamma_2\md\cdot\mn+\mh)dx
=\int_{\R^2}\mn\cdot\Delta^2\momega\cdot\mh dx.
\end{eqnarray*}
Using $(\ref{eq:EL})_3$ again, we obtain
\begin{eqnarray*}
&&-\gamma_2\int_{\R^2}\mn\mN:\Delta^2\md dx\\
&&\leq -\frac{\gamma_2}{\gamma_1}\int_{\R^2}(\Delta^2\md\cdot\mn)\cdot(\mn\times(\mh\times\mn))dx+
\frac{\gamma_2^2}{\gamma_1}\int_{\R^2}(\Delta^2\md\cdot\mn)\cdot(\md\cdot\mn))dx\\
&&\quad-\frac{\gamma_2^2}{\gamma_1}\int_{\R^2}(\Delta^2\md\cdot\mn)\cdot((\md\cdot\mn)\cdot\mn\mn)dx\\
&&\le  - \frac{\gamma_2}{\gamma_1}\int_{\R^2}(\Delta\md\cdot\mn)\cdot(\mn\times(\Delta\mh\times\mn))dx
-\frac{\gamma_2^2}{\gamma_1}\int_{\R^2}(|\Delta\md\cdot\mn|^2-|\mn^i\Delta\md_{ij}\cdot\mn^j|^2)dx\\
&&\quad+C\int_{\R^2}|\nabla^3\mv|(|\nabla^2\mn||\nabla\mv|+
|\nabla^2\mn|^2+|\nabla^2\mv||\nabla\mn|+|\nabla^3\mn||\nabla\mn|)dx
\end{eqnarray*}
Thus, combining the above estimates and applying Remark 2.2, we get
\ba\label{4.130}
&&\frac{1}{2}\frac{d}{dt}\int_{\R^2}|\nabla^2\mv|^2dx+ \frac{\gamma}{2Re}\int_{\R^2}|\nabla\Delta\mv|^2dx\nonumber\\
&&\le C\int_{\R^2}|\nabla\mv|^4+|\nabla^2\mn|^4dx
+C\int_{\R^2}(|\nabla^3\mn|^2+|\nabla^2\mv|^2)(|\nabla\mn|^2+|\mv|^2)dx\nonumber\\
&&\quad+\frac{1-\gamma}{Re}\int_{\R^2}\mn\cdot\Delta\momega\cdot\Delta\mh dx
+\frac{1-\gamma}{Re}\mu_2\int_{\R^2}(\mn\times(\Delta\mh\times\mn))\cdot(\Delta\md\cdot\mn)dx\nonumber\\
\ea

On the other hand, we can differentiate $\nabla_\beta$ to (\ref{eq:equal n}),  multiply it by $\nabla_\beta\Delta \mh$, and we get
\ba\label{4.140}
&&\frac{d}{dt}\int_{\R^2}a|\nabla\Delta \mn|^2+(k_1-a)|\triangle{\rm div}\mn|^2
+(k_2-a)|\Delta(\mn\times\mcurl\mn)|^2+(k_3-a)|\Delta(\mn\cdot\mcurl\mn)|^2dx\nonumber\\
&&\quad+\int_{\R^2}[(\nabla_\beta\mv\cdot\nabla)n^i
+(\mv\cdot\nabla)\nabla_{\beta}n^{i}]\cdot\nabla_{\beta}\Delta h^{i}dx\nonumber\\
&&\leq-\int_{\R^2}\nabla_{\beta}(\momega\cdot\mn)\cdot\nabla_{\beta}\Delta\mh dx+\mu_1\int_{\R^2}\nabla_{\beta}(\mn\times(\mh\times\mn))\cdot\nabla_{\beta}\Delta\mh dx\nonumber\\
&&\quad+\mu_2 \int_{\R^2}\nabla_{\beta}(\mn\times((\md\cdot\mn)\times\mn))\cdot\nabla_{\beta}\Delta\mh dx\nonumber\\
&&\quad+\delta\int_{\R^2}|\nabla^2\mn_t|^2
+C(\delta)\int_{\R^2}(|\nabla^2\mn|^4+|\nabla^3\mn|^2|\nabla\mn|^2+|\nabla^2\mn|^2|\nabla\mn|^4)dx\nonumber\\
&&\doteq B_6'+B_7'+B_8'+\delta\int_{\R^2}|\nabla^2\mn_t|^2
+C(\delta)\int_{\R^2}(|\nabla^2\mn|^4+|\nabla^3\mn|^2|\nabla\mn|^2+|\nabla^2\mn|^2|\nabla\mn|^4)dx,\nonumber\\
\ea
where $\delta>0$, to be decided later.

Since $|\Delta\mh|\leq C(|\nabla^4\mn|+|\nabla^3\mn||\nabla\mn|+|\nabla^2\mn|^2)$, then direct calculation  shows
\ben
&&|\int_{\R^2}[(\nabla_\beta\mv\cdot\nabla)n^i
+(\mv\cdot\nabla)\nabla_{\beta}n^{i}]\cdot\nabla_{\beta}\Delta h^{i}dx|\nonumber\\
&&\leq \frac{1}{16}a\mu_1\int_{\R^2}|\nabla^4\mn|^2dx+C\int_{\R^2}|\nabla\mv|^4+|\nabla^2\mn|^4dx\nonumber\\
&&+C\int_{\R^2}(|\nabla^3\mn|^2+|\nabla^2\mv|^2)(|\nabla\mn|^2+|\mv|^2)dx,
\een
and
\ben
B_6'&=&\int_{\R^2}(\Delta\momega\cdot\mn)\cdot\Delta\mh dx+\int_{\R^2}(\Delta(\momega\cdot\mn)-\Delta\momega\cdot\mn)\cdot\Delta\mh dx\nonumber\\
&\le& \int_{\R^2}\Delta\mh\cdot(\Delta\momega\cdot\mn) dx+ \frac{1}{16}a\mu_1\int_{\R^2}|\nabla^4\mn|^2dx\nonumber\\
&&+C\int_{\R^2}|\nabla\mv|^4+|\nabla^2\mn|^4dx
+C\int_{\R^2}(|\nabla^3\mn|^2+|\nabla^2\mv|^2)(|\nabla\mn|^2+|\mv|^2)dx\nonumber\\
\een
Note the fact that $|\mn\cdot\Delta^2\mn|\leq C(|\nabla\mn||\nabla^3\mn|+|\nabla^2\mn|^2)$, and similar estimates to (\ref{nhn}) for $B_7'$ imply
\ben
B_7'\le -\mu_1 a\int_{\R^2}|\nabla^4\mn|^2+C\int_{\R^2}(|\nabla\mn|^2|\nabla^3\mn|^2+|\nabla^2\mn|^4)dx.
\een
At last, we estimate $B_8'$.
\ben
B_8'&=&\mu_2\int_{\R^2}\{(\nabla_{\beta}\mn\times((\md\cdot\mn)\times\mn))\cdot\nabla_{\beta}\Delta\mh
 + (\mn\times((\md\cdot\mn)\times\nabla_{\beta}\mn))\cdot\nabla_{\beta}\Delta\mh\nonumber\\
 &&+(\mn\times((\md\cdot\nabla_{\beta}\mn)\times\mn))\cdot\nabla_{\beta}\Delta\mh\} dx+\mu_2\int_{\R^2} (\mn\times((\nabla_{\beta}\md\cdot\mn)\times\mn))\cdot\nabla_{\beta}\Delta\mh dx\nonumber\\
 &\le& \frac{1}{16}a\mu_1\int_{\R^2}|\nabla^4\mn|^2dx+C\int_{\R^2}|\nabla\mv|^4+|\nabla^2\mn|^4dx\nonumber\\
&&+C\int_{\R^2}(|\nabla^3\mn|^2+|\nabla^2\mv|^2)(|\nabla\mn|^2+|\mv|^2)dx\nonumber\\ &&-\mu_2\int_{\R^2}\mn\times((\Delta\md\cdot\mn)\times\mn)\cdot\Delta\mh dx-\mu_2\int_{\R^2}(\mn\times((\nabla_{\beta}\md\cdot\mn)\times\nabla_{\beta}\mn))\cdot\Delta\mh dx\nonumber\\
 &&-\mu_2\int_{\R^2}\{(\nabla_{\beta}\mn\times((\nabla_{\beta}\md\cdot\mn)\times\mn))
  +(\mn\times((\nabla_{\beta}\md\cdot\nabla_{\beta}\mn)\times\mn))\}\Delta\mh dx\nonumber\\
  &\le&-\mu_2\int_{\R^2}((\Delta\md\cdot\mn)\cdot(\mn\times(\Delta\mh \times\mn)) dx+\frac{1}{8}a\mu_1\int_{\R^2}|\nabla^4\mn|^2dx\nonumber\\
 &&+C\int_{\R^2}|\nabla\mv|^4+|\nabla^2\mn|^4dx
+C\int_{\R^2}(|\nabla^3\mn|^2+|\nabla^2\mv|^2)(|\nabla\mn|^2+|\mv|^2)dx\nonumber\\
 \een
The definition of $\mh$ and (\ref{eq:equal n}) tell us that
\ben\label{eq:4.190}
 |\nabla^2 \mn_t|^2&\leq& \max\{1,|\mu_1|,|\mu_2|\} (|\nabla^4 \mn|^2+|\nabla^3 \mv|^2)+C|\nabla^2\mv|^2|\nabla\mn|^2\nonumber\\
 &&+C|\nabla^2\mn|^2(|\nabla\mv|^2+|\nabla^2\mn|^2).
\een
Hence combining these estimates from (\ref{4.130}) to (\ref{eq:4.190}) and choosing $\delta$ small enough, by Gagliardo-Sobolev inequlity we obtain
 \beno
&&\frac{d}{dt} \int_{\R^2}\frac{Re}{2(1-\gamma)}|\nabla^2\mv|^2 +a|\nabla\Delta\mn|^2dx
  +\int_{\R^2}\frac{\gamma}{4(1-\gamma)}|\nabla^3\mv|^2 +\frac{a\mu_1}{4}|\nabla^4\mn|^2 dx\nonumber\\
&&\quad+\frac{d}{dt}\int_{\R^2}(k_1-a)|\Delta{\rm div}\mn|^2
+(k_2-a)|\Delta(\mn\times\mcurl\mn)|^2+(k_3-a)|\Delta(\mn\cdot\mcurl\mn)|^2dx\nonumber\\
&&\le C\int_{\R^2}|\nabla\mv|^4+|\nabla^2\mn|^4dx
+C\int_{\R^2}(|\nabla^3\mn|^2+|\nabla^2\mv|^2)(|\nabla\mn|^2+|\mv|^2)dx,\nonumber\\
&&\le C\big(\int_{\R^2}(|\mv|^4+|\nabla\mn|^4)dx\big)^{1/2}\big(\int_{\R^2}(|\nabla^2\mv|^4+|\nabla^3\mn|^4)dx\big)^{1/2}
+C\int_{\R^2}|\nabla\mv|^4+|\nabla^2\mn|^4dx\\
&&\le \min\{\frac{\gamma}{8(1-\gamma)},\frac{a\mu_1}{8}\}\int_{\R^2}|\nabla^3\mv|^2 +|\nabla^4\mn|^2 dx\\
&&\quad+C\int_{\R^2}(|\mv|^4+|\nabla\mn|^4+|\nabla\mv|^2+|\nabla^2\mn|^2)dx\int_{\R^2}(|\nabla^2\mv|^2+|\nabla^3\mn|^2)dx
\eeno
Then the required inequality follows from Gronwall's inequality, (\ref{eq:nabla2n}), (\ref{eq:nablan4}) and Lemma 4.4.\endproof

%---------------------------------------------------cOR 4.6------------------------------------------------------------
Similar computations as the above two lemmas, we may obtain interior regularity for $(\mv,\mn)$ of (\ref{eq:EL}).
\begin{corollary} \label{lem:4.5}
Assume that the Leslie coefficients satisfy (\ref{Parodi})-(\ref{beta}). Let $(\mv,\mn)$  be a solution of (\ref{eq:EL})
 with the initial value $(\mv_0,\mn_0)\in L^2\times H_{b}^{1}(\R^2,S^2)$ and ${\rm div}\mv_0=0$. There are constants
 $\epsilon_1$ and $R_0>0$ such that
 $$
 {\rm esssup}_{0\le t\le T,x\in\R^2}\int_{B_{R}(x)}|\nabla\mn(\cdot,t)|^2+|\mv|^2(\cdot,t)dx<\epsilon_1,\quad\forall R\in(0,R_0].$$
Then, for all $t\in[\tau, T]$ with $\tau\in (0,T)$, for any $l\geq 1$ it holds that
\ben\label{eq:l regularity}
&&\int_{\R^2}|\nabla^{l+1}\mn(\cdot,t)|^2+|\nabla^l\mv|^2(\cdot,t)dx +\int_{\tau}^t\int_{\R^2}|\nabla^{l+2}\mn(\cdot,s)|^2+|\nabla^{l+1}\mv|^2(\cdot,s)dxds
\nonumber\\
&&\le C(l,\epsilon_1,E_0,\tau, T, \frac{T}{R^2}).
\een
Moreover, $\mn$ and $\mv$ are regular for all $t\in(0,T)$.
\end{corollary}

{\it Proof:} The proof is divided into four steps.

{\bf Step 1.} For $l=1,2$, by Lemma 4.4 and 4.5 we know that (\ref{eq:l regularity}) holds. For $l=3$, similar arguments as in Lemma 4.5, it's not difficult to obtain (\ref{eq:l regularity}). Now we assume that for $ l_0\geq 3$, the estimate (\ref{eq:l regularity}) holds, and we prove that the case $l=l_0+1$ is still true.

{\bf Step 2.} Firstly, we prove the case that $l_0$ is odd, i.e. for $l=1,\cdots,2k-1$ and  $k\geq 2$, (\ref{eq:l regularity}) holds, and we are aimed  to obtain the case $l=2k$. Especially, by the assumptions and Sobolev inequality we have the following estimates:
\ben\label{eq:bound of nabla v}
|\nabla\mn|+|\nabla^2\mn|+|\mv|+|\nabla\mv|\leq C(\epsilon_1,E_0,\tau, T,  \frac{T}{R^2}), \quad {\rm in}\,\,\R^2\times(\tau,T),
\een
and for any $t\in (\tau,T)$,
\ben\label{eq:bond of n 2k-1 }
&&\int_{\R^2}|\nabla^{2k}\mn(\cdot,t)|^2+|\nabla^{2k-1}\mv|^2(\cdot,t)dx+\int_{\tau}^t\int_{\R^2}|\nabla^{2k+1}\mn(\cdot,s)|^2+|\nabla^{2k}\mv|^2(\cdot,s)dxds
 \nonumber\\
&&\le C(k,\epsilon_1,E_0,\tau,  T, \frac{T}{R^2}),
\een
which implies that
\ben\label{eq:bound of n 2k-2}
|\nabla^{2k-2}\mn|+|\nabla^{2k-3}\mv|\leq C(k,\epsilon_1,E_0,\tau, T,  \frac{T}{R^2}), \quad {\rm in}\,\,\R^2\times(\tau,T).
\een
As the process in Lemma 4.5,
multiplying $(\ref{eq:EL})_1$ with $\Delta^{2k} \mv$, we have
\ba\label{5.100}
&&\frac{1}{2}\frac{d}{dt}\int_{\R^2}|\nabla^{2k}\mv|^2dx+ \frac{\gamma}{Re}\int_{\R^2}|\nabla\Delta^k\mv|^2dx\nonumber\\
&&=-\int_{\R^2}(\mv\cdot\nabla)\mv\cdot\Delta^{2k} \mv dx+\frac{1-\gamma}{Re}\int_{\R^2}(\nabla\cdot\msigma^{E})\cdot\Delta^{2k}\mv dx+\frac{1-\gamma}{Re}\int_{\R^2}(\nabla\cdot\msigma^{L})\cdot\Delta^{2k}\mv dx\nonumber\\
&&\le \frac{1}{4}\frac{\gamma}{Re}\int_{\R^2}|\nabla\Delta^k\mv|^2dx+C\int_{\R^2}|\nabla^{2k-1}(\mv\cdot\nabla\mv)|^2dx
+\frac{1-\gamma}{Re}\int_{\R^2}(\nabla\cdot\msigma^{L})\cdot\Delta^{2k}\mv dx\nonumber\\
&&\quad +C\int_{\R^2}|\nabla^{2k}\msigma^{E}|^2dx.
\ea
Recall that the form of $\msigma^{E}$ is like $\nabla\mn\otimes\nabla\mn+\nabla\mn\otimes\nabla\mn\otimes\mn\otimes\mn$, hence by (\ref{eq:bound of nabla v})-(\ref{eq:bound of n 2k-2}) we have
\beno
\int_{\R^2}|\nabla^{2k}\msigma^{E}|^2dx\leq C(k,\epsilon_1,E_0,\tau,  T, \frac{T}{R^2})\big(\int_{\R^2}|\nabla^{2k+1}\mn|^2dx+1\big);
\eeno
similarly,
\beno
\int_{\R^2}|\nabla^{2k-1}(\mv\cdot\nabla\mv)|^2dx\leq C(k,\epsilon_1,E_0,\tau, T,  \frac{T}{R^2})\big(\int_{\R^2}|\nabla^{2k}\mv|^2dx+1\big).
\eeno
Note that
\beno
&&\frac{1-\gamma}{Re}\int_{\R^2}(\nabla\cdot\msigma^{L})\cdot\Delta^{2k}\mv dx \\
&&=- \frac{1-\gamma}{Re}\int_{\R^2}[\alpha_1(\mn\mn:\md)\mn\mn:\Delta^{2k}\md+(\alpha_2+\alpha_3)\mn\mN:\Delta^{2k}\md
+(\alpha_2-\alpha_3)\mn\mN:\Delta^{2k}\Omega\\
&&\quad+\alpha_4\md:\Delta^{2k}\md+(\alpha_5+\alpha_6)(\mn\mn\cdot\md):\Delta^{2k}\md
+(\alpha_5-\alpha_6)(\mn\mn\cdot\md):\Delta^{2k}\momega] dx.
\eeno
By (\ref{eq:bound of nabla v})-(\ref{eq:bound of n 2k-2}), we get
\beno
&&-\alpha_1\int_{\R^2}(\mn\mn:\md)\mn\mn:\Delta^{2k}\md dx\\
&&\le -\alpha_1\int_{\R^2}|\mn\mn:\Delta^{k}\md|dx +\delta\int_{\R^2}|\nabla^{2k+1}\mv|^2dx+C(\delta,k,\epsilon_1,E_0,\tau, T,  \frac{T}{R^2})(\int_{\R^2}|\nabla^{2k}\mv|^2dx+1),
\eeno
where $\delta>0$, to be decided. Moreover,
\beno
&&-\int_{\R^2} (\mn\mn\cdot\md):\Delta^{2k}\md dx\\
&&\le -\int_{\R^2}|(\Delta^k \md\cdot\mn)|^2dx+\delta\int_{\R^2}|\nabla^{2k+1}\mv|^2dx+C(\delta,k,\epsilon_1,E_0,\tau, T,  \frac{T}{R^2})(\int_{\R^2}|\nabla^{2k}\mv|^2dx+1).
\eeno
Since formally $h$ is $\nabla^2\mn+\nabla^2\mn\otimes\mn\otimes\mn+\nabla\mn\otimes\nabla\mn\otimes\mn$,
by $(\ref{eq:EL})_3$, the anti-symmetric property of $\Delta^k\momega$ and (\ref{eq:bound of nabla v})-(\ref{eq:bound of n 2k-2}) we get
\beno
&&\int_{\R^2}\gamma_1(\mn\mN:\Delta^{2k}\momega)+\gamma_2(\mn\mn\cdot\md):\Delta^{2k}\momega dx\\
&&=\int_{\R^2}\mn\cdot\Delta^{2k}\momega\cdot(-\mh+\gamma_1\mN+\gamma_2\md\cdot\mn+\mh)dx
=\int_{\R^2}\mn\cdot\Delta^{2k}\momega\cdot\mh dx\\
&&\leq \int_{\R^2}\mn\cdot\Delta^{k}\momega\cdot\Delta^{k}\mh dx+\delta\int_{\R^2}|\nabla^{2k+1}\mv|^2dx\\
&&+C(\delta,k,\epsilon_1,E_0,\tau,  T, \frac{T}{R^2})(\int_{\R^2}|\nabla^{2k+1}\mn|^2dx+1).
\eeno
Using $(\ref{eq:EL})_3$ again, similarly we obtain
\beno
&&-\gamma_2\int_{\R^2}\mn\mN:\Delta^{2k}\md dx\\
&&\leq -\frac{\gamma_2}{\gamma_1}\int_{\R^2}(\Delta^{2k}\md\cdot\mn)\cdot(\mn\times(\mh\times\mn))dx+
\frac{\gamma_2^2}{\gamma_1}\int_{\R^2}(\Delta^{2k}\md\cdot\mn)\cdot(\md\cdot\mn))dx\\
&&\quad-\frac{\gamma_2^2}{\gamma_1}\int_{\R^2}(\Delta^{2k}\md\cdot\mn)\cdot((\md\cdot\mn)\cdot\mn\mn)dx\\
&&\le  - \frac{\gamma_2}{\gamma_1}\int_{\R^2}(\Delta^k\md\cdot\mn)\cdot(\mn\times(\Delta^k\mh\times\mn))dx
+\frac{\gamma_2^2}{\gamma_1}\int_{\R^2}(|\Delta^k\md\cdot\mn|^2-|\mn^i\Delta^k\md_{ij}\mn^j|^2)dx\\
&&\quad+\delta\int_{\R^2}|\nabla^{2k+1}\mv|^2dx+C(\delta,k,\epsilon_1,E_0,\tau,  T, \frac{T}{R^2})(\int_{\R^2}|\nabla^{2k}\mv|^2dx+\int_{\R^2}|\nabla^{2k+1}\mn|^2dx+1).
\eeno
Thus, combining the above estimates and applying Remark 2.2, we get
\ba\label{4.1300}
&&\frac{1}{2}\frac{d}{dt}\int_{\R^2}|\nabla^{2k}\mv|^2dx+ \frac{\gamma}{2Re}\int_{\R^2}|\nabla\Delta^k\mv|^2dx\nonumber\\
&&\le \frac{1-\gamma}{Re}\int_{\R^2}\mn\cdot\Delta^k\momega\cdot\Delta^k\mh dx
+\frac{1-\gamma}{Re}\mu_2\int_{\R^2}(\mn\times(\Delta^k\mh\times\mn))\cdot(\Delta^k\md\cdot\mn)dx\nonumber\\
&&\quad+\delta\int_{\R^2}|\nabla^{2k+1}\mv|^2dx+C(\delta,k,\epsilon_1,E_0,\tau,  T, \frac{T}{R^2})(\int_{\R^2}|\nabla^{2k}\mv|^2dx+\int_{\R^2}|\nabla^{2k+1}\mn|^2dx+1).\nonumber\\
\ea

Differentiate $\nabla_\beta$ to (\ref{eq:equal n}),  multiply it by $\nabla_\beta\Delta^{2k-1} \mh$, and by (\ref{eq:bound of nabla v})-(\ref{eq:bound of n 2k-2}) we get
\ba\label{4.1400}
&&\frac{d}{dt}\int_{\R^2}a|\nabla\Delta^k \mn|^2+(k_1-a)|\triangle^k{\rm div}\mn|^2
+(k_2-a)|\Delta^k(\mn\times\mcurl\mn)|^2\nonumber\\
&&\quad+\frac{d}{dt}\int_{\R^2}(k_3-a)|\Delta^k(\mn\cdot\mcurl\mn)|^2dx
+\int_{\R^2}[(\nabla_\beta\mv\cdot\nabla)n^i
+(\mv\cdot\nabla)\nabla_{\beta}n^{i}]\cdot\nabla_{\beta}\Delta^{2k-1} h^{i}dx\nonumber\\
&&\leq-\int_{\R^2}\nabla_{\beta}(\momega\cdot\mn)\cdot\nabla_{\beta}\Delta^{2k-1}\mh dx+\mu_1\int_{\R^2}\nabla_{\beta}(\mn\times(\mh\times\mn))\cdot\nabla_{\beta}\Delta^{2k-1}\mh dx\nonumber\\
&&\quad+\mu_2 \int_{\R^2}\nabla_{\beta}(\mn\times((\md\cdot\mn)\times\mn))\cdot\nabla_{\beta}\Delta^{2k-1}\mh dx\nonumber\\
&&\quad+\delta\int_{\R^2}|\nabla^{2k}\mn_t|^2
+C(\delta,k,\epsilon_1,E_0,\tau,  T, \frac{T}{R^2})(\int_{\R^2}|\nabla^{2k+1}\mn|^2dx+1)\nonumber\\
&&\doteq B_6''+B_7''+B_8''+\delta\int_{\R^2}|\nabla^{2k}\mn_t|^2
+C(\delta,k,\epsilon_1,E_0,\tau, T,  \frac{T}{R^2})(\int_{\R^2}|\nabla^{2k+1}\mn|^2dx+1)\nonumber\\
\ea
Using (\ref{eq:bound of nabla v})-(\ref{eq:bound of n 2k-2}) again, direct calculation  shows
\ben
&&|\int_{\R^2}[(\nabla_\beta\mv\cdot\nabla)n^i
+(\mv\cdot\nabla)\nabla_{\beta}n^{i}]\cdot\nabla_{\beta}\Delta^{2k-1} h^{i}dx|\nonumber\\
&&\leq \delta\int_{\R^2}|\nabla^{2k+2}\mn|^2dx+C(\delta,k,\epsilon_1,E_0,\tau, T,  \frac{T}{R^2})(\int_{\R^2}|\nabla^{2k}\mv|^2dx+\int_{\R^2}|\nabla^{2k+1}\mn|^2dx+1),\nonumber\\
\een
and
\ben
B_6''&=&\int_{\R^2}(\Delta\momega\cdot\mn)\cdot\Delta^{2k-1}\mh dx+\int_{\R^2}(\Delta(\momega\cdot\mn)-\Delta\momega\cdot\mn)\cdot\Delta^{2k-1}\mh dx\nonumber\\
&\le& \int_{\R^2}\Delta^k\mh\cdot(\Delta^k\momega\cdot\mn) dx+ \delta\int_{\R^2}|\nabla^{2k+2}\mn|^2dx\nonumber\\
&&+C(\delta,k,\epsilon_1,E_0,\tau,  T, \frac{T}{R^2})(\int_{\R^2}|\nabla^{2k}\mv|^2dx+\int_{\R^2}|\nabla^{2k+1}\mn|^2dx+1)
\een
By (\ref{eq:bound of n 2k-2}) and the fact $|\mn|=1$, we remark that
\beno
|\mn\cdot\Delta^{k+1}\mn|\leq C(k,\epsilon_1,E_0,\tau,  T, \frac{T}{R^2})(|\nabla^{2k+1}\mn|+|\nabla^{2k}\mn|+|\nabla^{2k-1}\mn|),
\eeno
and similar estimates for $B_7''$ as in (\ref{nhn}) imply that
\ben
B_7''\le -\mu_1 a\int_{\R^2}|\nabla^{2k+2}\mn|^2+C(k,\epsilon_1,E_0,\tau, T,  \frac{T}{R^2})(\int_{\R^2}|\nabla^{2k+1}\mn|^2dx+1).
\een
At last, we estimate $B_8''$. Noting the bound of $\Delta^k\mh$, we get
\ben\label{eq:4.1900}
B_8''&=&\mu_2\int_{\R^2}\{(\nabla_{\beta}\mn\times((\md\cdot\mn)\times\mn))\cdot\nabla_{\beta}\Delta^{2k-1}\mh
 + (\mn\times((\md\cdot\mn)\times\nabla_{\beta}\mn))\cdot\nabla_{\beta}\Delta^{2k-1}\mh\nonumber\\
 &&+(\mn\times((\md\cdot\nabla_{\beta}\mn)\times\mn))\cdot\nabla_{\beta}\Delta^{2k-1}\mh\} dx+\mu_2\int_{\R^2} (\mn\times((\nabla_{\beta}\md\cdot\mn)\times\mn))\cdot\nabla_{\beta}\Delta^{2k-1}\mh dx\nonumber\\
 &\le&\delta\int_{\R^2}|\nabla^{2k+2}\mn|^2dx+C(\delta,k,\epsilon_1,E_0,\tau,  T, \frac{T}{R^2})(\int_{\R^2}|\nabla^{2k}\mv|^2dx+1)\nonumber\\ &&-\mu_2\int_{\R^2}\nabla_{\beta}\Delta^{k-1}\big(\mn\times((\nabla_{\beta}\md\cdot\mn)\times\mn)\big)\cdot\Delta^{k}\mh dx\nonumber\\
  &\le&-\mu_2\int_{\R^2}((\Delta^k\md\cdot\mn)\cdot(\mn\times(\Delta^k\mh \times\mn)) dx+2\delta\int_{\R^2}|\nabla^{2k+2}\mn|^2dx\nonumber\\
 &&+C(\delta,k,\epsilon_1,E_0,\tau, T,  \frac{T}{R^2})(\int_{\R^2}|\nabla^{2k}\mv|^2dx+1).
 \een
Hence combining these estimates from (\ref{4.1300}) to (\ref{eq:4.1900}) and choosing $\delta$ small enough, we obtain
 \beno
&&\frac{d}{dt} \int_{\R^2}\frac{Re}{2(1-\gamma)}|\nabla^{2k}\mv|^2 +a|\nabla^{2k+1}\mn|^2dx
  +\int_{\R^2}\frac{\gamma}{4(1-\gamma)}|\nabla^{2k+1}\mv|^2 +\frac{a\mu_1}{4}|\nabla^{2k+2}\mn|^2 dx\nonumber\\
&&\quad+\frac{d}{dt}\int_{\R^2}(k_1-a)|\Delta^k{\rm div}\mn|^2
+(k_2-a)|\Delta^k(\mn\times\mcurl\mn)|^2+(k_3-a)|\Delta^k(\mn\cdot\mcurl\mn)|^2dx\nonumber\\
&&\le C(k,\epsilon_1,E_0,\tau,  T,  \frac{T}{R^2})(\int_{\R^2}|\nabla^{2k}\mv|^2dx+\int_{\R^2}|\nabla^{2k+1}\mn|^2dx+1)
\eeno
which implies (\ref{eq:l regularity}) for $l=2k$ by Gronwall's inequality.

{\bf Step 3.} Secondly, the case that $l_0$ is even may be obtained by the same step as above.

{\bf Step 4.} By arbitrarily of $\tau$,
the estimates in Step 2. and 3., we have for any $l>1$ and $(x,t)\in \R^2\times(0,T)$
\beno
|\nabla^l\mn|(x,t)+|\nabla^l\mv|(x,t)<\infty,
\eeno
which yields that $\mv, \mn$ are spatially smooth. Using the equation (\ref{eq:EL}), $\mv_t, \mn_t$ are also spatially smooth, and differentiating to the time $t$, hence we finally obtain $\mv, \mn$ are smooth in $\R^2\times(0,T)$.
\endproof

%Higher regularity estimates are derived by standard arguments as in \cite[Lemma 3.10]{St1} or \cite[Prop 5.2]{HX}, and we omitted the details.

\subsection{Global existence}
Now we'll complete the proof of Theorem 1.2. Indeed it's more or less standard since the local monotonic inequality , $\varepsilon$-regularity estimates in Section 4 and local existence of strong solutions for some regular data in Section 3 have been obtained. The following arguments are similar to \cite{St1} and \cite{LLW}, where the main difference is dealing with the Leslie coefficients, and we sketch its proof.

Let $E(t)=E(\mv,\mn)(t)=\int_{\R^2}e(\mv,\mn)(\cdot,t)dx$ denote the energy of $(\mv,\mn)$ at time $t$ . For any data $\mn_0\in H_{b}^{1}(\R^2;S^2)$, one can approximate it by a sequence of smooth maps $\mn_0^k$ in $H_{b}^{1}(\R^2;S^2)$, and we can assume that $\mn_0^k\in H_{b}^{4}(\R^2;S^2)$ (see \cite{SU}). For $\mv_0\in L^2(\R^2;\R^2)$, there exists a sequence of $\mv_0^k\in C_0^{\infty}(\R^2;\R^2)$ and $\mv_0^k\rightarrow \mv_0$ in $L^2(\R^2;\R^2)$.

Due to the absolute continuity property of the integral, for any $\epsilon_1>0$, there exists $R_0>0$ such that
$$
 \sup_{x\in \R^2}\int_{B_{R_0}(x)}|\nabla\mn_0|^2+|\mv_0|^2dx\le \epsilon_1,
$$
and by the strong convergence of $\mn_0^k$ and $\mv_0^k$
$$
 \sup_{x\in \R^2}\int_{B_{R_0}(x)}|\nabla\mn_0^k|^2+|\mv_0^k|^2dx\le 2\epsilon_1,
$$
for $k$ is large enough. Without loss of generality, we assume that it holds for all $k\geq 1.$

For the data $\mn_0^k$ and $\mv_0^k$, by Theorem 1.1 there exists a time $T^k$ and a strong solution $(\mn^k, \mv^k)$ with the pressure $p^k$ such that
$$
\mv^k\in C\left([0,T^k];H^{4}(\R^2)\right)\cap L^{2}(0,T;H^{5}(\R^2)),
\nabla\mn^k\in C\left([0,T^k];H^{4}(\R^2)\right).
$$
Hence there exists $T_0^k\leq T^k$ such that
$$
 \sup_{0<t<T_0^k,\\ x\in \R^2}\int_{B_{R}(x)}|\nabla\mn^k(y,t)|^2+|\mv^k(y,t)|^2dy\le 4\epsilon_1,
$$
where $R\leq R_0/2.$
However, by the local monotonic inequality in Lemma 4.3, we have $T_0^k\geq \frac{\epsilon_1^2R_0^2}{4C_2^2E_0^2}=T_0>0$  uniformly.
For any $0<\tau<T_0$, by the estimates in Corollary 4.6 for any $l\geq 1$ we get
\begin{eqnarray}\label{eq:5.1}
&&\sup_{\tau<t<T_0}\int_{\R^2}|\nabla^{l+1}\mn^k|^2(\cdot,t)+|\nabla^{l}\mv^k|^2(\cdot,t)dx + \int_{\tau}^{T_0}\int_{\R^2}|\nabla^{l+2}\mn^k(\cdot,s)|^2+|\nabla^{l+1}\mv^k|^2(\cdot,s)dxds\nonumber\\
&&\le C(l,\epsilon_1,E_0,\tau,  T, \frac{T}{R^2}).
\end{eqnarray}
Moreover, the energy inequality in Proposition 2.1, a priori estimates in Lemma 4.2 and the equations of the direction fields yield that
$$
E(\mv^k,\mn^k)(t)\leq E_0,$$
and
\begin{eqnarray}\label{eq:5.2}
\int_{\R^2\times [0,T_0^k]}\big(|\nabla^2 \mn^k|^2+|\nabla\mv^k|^2+|\partial_t \mn^k|^2+|\nabla \mn^k|^4+|\mv^k|^4\big)dxdt\leq C(\epsilon_1,C_2, E_0).
\end{eqnarray}

Now we estimate the pressure term of the velocity equations. In the distributional sense,
\begin{eqnarray*}
\triangle p^k=-\nabla^2(\mv^k\otimes\mv^k-\frac{1-\gamma}{Re}\sigma^k),
\end{eqnarray*}
hence by Calder\'{o}n-Zygmund estimates, we get
\begin{eqnarray}\label{eq:5.3}
\|p^k\|_{L^{2}(\R^2\times [0,T_0^k])}^2&\leq& C\||\mv^k|^2+|\nabla\mn^k|^2+|\nabla\mv^k|+|\partial_t\mn^k|\|_{L^{2}(\R^2\times [0,T_0^k])}^2\nonumber\\
&\leq& C(\epsilon_1,C_2, E_0).
\end{eqnarray}

At last, we estimate the term $\partial_t \mv^k$. For any $\phi\in C_0^{\infty}(\R^2\times (0,T_0^k);\R^2)$,
\begin{eqnarray*}
&&\int_0^{T_0^k}\int_{\R^2}\partial_t \mv^k\phi dxdt\\
&&= \int_0^{T_0^k}\int_{\R^2}\big(\mv^k\otimes\mv^k-\frac{\gamma}{Re}\nabla\mv^k-\frac{1-\gamma}{Re}\sigma^k\big):\nabla\phi+p^k {{\rm div}\phi} dxdt\\
&&\leq C\||\mv^k|^2+|\nabla\mn^k|^2+|\nabla\mv^k|+|p^k|+|\partial_t\mn^k|\|_{L^{2}(\R^2\times [0,T_0^k])}^2\|\phi\|_{L^2_tH^1_x}\\
&&\leq C(\epsilon_1,C_2) E_0\|\phi\|_{L^2_tH^1_x},
\end{eqnarray*}
that is, for any $k\geq 1,$
\begin{eqnarray}\label{eq:5.4}
\|\partial_t\mv^k\|_{L^2(0,T^k_0;H^{-1}(\R^2))}\leq C(\epsilon_1,C_2,E_0).
\end{eqnarray}

Hence the above estimates (\ref{eq:5.1})-(\ref{eq:5.4}) and Aubin-Lions Lemma yield that there exist a solution $(\mv,\mn-b)\in W^{1,0}_2(\R^2\times [0,T_0];\R^2)\times W^{2,1}_2 (\R^2\times [0,T_0];\R^3)$ with the pressure $p$ such that
(at most up to s subsequence)
\begin{eqnarray*}
&&\mv^k\rightarrow \mv,\quad {\rm locally\,\, in}\quad W^{2,1}_2(\R^2\times (0,T_0);\R^2);\\
&&\mn^k-b\rightarrow \mn-b,\quad {\rm locally\,\, in}\quad W^{3,1}_2(\R^2\times (0,T_0);\R^3);\\
&&\mn^k\rightarrow \mn,\quad {\rm in}\quad L^{\infty}(\R^2\times (0,T_0);\R^3);\\
&&\mv^k\rightharpoonup \mv,\quad {\rm in}\quad W^{1,0}_2(\R^2\times [0,T_0];\R^2);\\
&& \mn^k-b\rightharpoonup \mn-b,\quad {\rm in}\quad W^{2,1}_2(\R^2\times [0,T_0];\R^3);\\
&&p^k\rightharpoonup p,\quad {\rm in}\quad L^2(\R^2\times [0,T_0];\R);\\
&&\mv^k(t)\rightharpoonup \mv(t),\quad {\rm in}\quad H^2(\R^2;\R^2)\quad {\rm for}\quad{a.e.} \quad t\in(\tau,T_0);\\
&& \mn^k(t)-b\rightharpoonup \mn(t)-b,\quad {\rm in}\quad H^{3}(\R^2;\R^3)\quad {\rm for}\quad{a.e.}\quad  t\in(\tau,T_0).
\end{eqnarray*}
By (\ref{eq:5.2}) and (\ref{eq:5.4}), $(\mv(t),\nabla\mn(t))\rightharpoonup (\mv_0,\nabla\mn_0)$ weakly in $L^2(\R^2)$, thus
$$E(\mv_0,\mn_0)\leq\liminf_{t\rightarrow0} E(\mv(t),\mn(t))$$
On the other hand, by the energy estimates of $(\mv^k,\mn^k)$, we have
$$E(\mv_0,\mn_0)\geq\limsup_{t\rightarrow0} E(\mv(t),\mn(t)).$$
Hence,  $(\mv(t),\nabla\mn(t))\rightarrow (\mv_0,\nabla\mn_0)$ strongly in $L^2(\R^2)$ and $(\mv,\mn)$ is the solution of the equations (\ref{eq:EL}) with the indicated data $(\mv_0,\mn_0).$ From the weak limit of regular estimates (\ref{eq:5.1}), we know that $(\mv,\mn)\in C^{\infty} (\R^2\times(0,T_0])$, and $(\nabla^l\mv,\nabla^{l+1}\mn)(\cdot ,T_0)\in L^2(\R^2) $ for any $l\geq1$. By Theorem 1.1, there exists a unique smooth solution of (\ref{eq:EL}) with the initial data $(\mv,\mn)(\cdot,T_0)$, which is still written as $(\mv,\mn)$, and blow-up criterion yields that if $(\mv,\mn)$ blows up at time $T^*$, then
\beno
\|\nabla\times\mv\|_{L^{\infty}(\R^2)}(t)+\|\nabla\mn\|_{L^{\infty}(\R^2)}^2(t)\rightarrow\infty,\quad {\rm as} \quad t\rightarrow T^*,
\eeno
as a result,
\ben\label{eq:blow-up}
|\nabla^{4}\mn|(x,t)+|\nabla^{3}\mv|(x,t)\not\in L^{\infty}_tL^{2}_x((T_0,T^*)\times\R^2)
\een
We assume that $T_1$ is the first singular time of  $(\mv,\mn)$, then we have
\begin{eqnarray*}
(\mv,\mn)\in C^{\infty}(\R^2\times (0,T_1); \R^2\times S^2),\quad {\rm and}
\quad (\mv,\mn)\not\in C^{\infty}(\R^2\times (0,T_1]; \R^2\times S^2);
\end{eqnarray*}
and by Corollary 4.6 and (\ref{eq:blow-up}),
\begin{eqnarray*}
\lim\sup_{t\uparrow T_1}\max_{x\in \R^2}\int_{B_R(x)}(|\mv|^2+|\nabla\mn|^2)(\cdot,t)\geq \epsilon_1,\quad \forall R>0.
\end{eqnarray*}
Finally, we can prove that $(\mv,\mn-b)\in C^0([0,T_1], L^2(\R^2))$ by similar arguments as (\ref{eq:5.4}) (also see P330, \cite{LLW}). Hence, we can define
$$(\mv(T_1), \mn(T_1)-b)=\lim_{t\uparrow T_1}(\mv(t), \mn(t)-b) \quad {\rm in}\quad  L^2(\R^2).$$
On the other hand, by the energy inequality $\nabla\mn\in L^{\infty}(0,T_1;L^2(R^2))$, hence
$\nabla\mn(t)\rightharpoonup \nabla\mn(T_1)$. Similarly we can extend $T_1$ to $T_2$ and so on. It's easy to check that the energy loss at every singular time $T_i$ for $i\geq 1$ is at least $\epsilon_1$, thus the number of the singular time is finite as $L$, and for $1\leq i\leq L$ we have
\begin{eqnarray*}
\lim\sup_{t\uparrow T_i}\max_{x\in \R^2}\int_{B_R(x)}(|\mv|^2+|\nabla\mn|^2)(\cdot,t)\geq \epsilon_1,\quad \forall R>0.
\end{eqnarray*}
The proof is complete.
\endproof
\begin{remark}\label{rem:hongxin}
Formally, if the Leslie coefficients $\alpha_1,\cdots,\alpha_6$, $\momega$ and $\md$ vanish in (\ref{eq:EL}), then we arrive at the Oseen-Frank model as in \cite{HX}, and all the computations and arguments above also work, which seem to be more simple.
\end{remark}

\bigskip

\noindent {\bf Acknowledgments.}
The authors would like to thank Professor Zhifei ZHANG for his suggestion to consider the problem and his some valuable discussions with them. Part of this work is carried out when the first author is visiting Math. Department of Princeton University. Meng is partially supported by  NSFC  10931001£¬and Chen-su star project by Zhejiang University. Wendong is supported by "the Fundamental Research Funds for the Central Universities" and
the Institute of Mathematical Sciences of Chinese University of Hong Kong.

\end{document}